\documentclass[preprint,12pt,english]{elsarticle}

\usepackage[margin=2cm]{geometry}
\usepackage{mathrsfs}
\usepackage{amsfonts}
\usepackage{dsfont}
\usepackage{multirow}
\usepackage{bbold}
\usepackage{amssymb}
\usepackage{float}
\usepackage{placeins}
\usepackage[utf8]{inputenc}
\usepackage{hyperref}  
\usepackage{tikz}
\usetikzlibrary{backgrounds,arrows, shapes,decorations.markings,positioning}
\usetikzlibrary{plotmarks,calc,fadings,decorations.pathreplacing,decorations.pathmorphing}
\tikzset{%
  >=latex,
  inner sep=0pt,%
  outer sep=2pt,%
  mark coordinate/.style={inner sep=0pt,outer sep=0pt,minimum size=3pt,
    fill=black,circle}%
}

\tikzset{-dot-/.style={decoration={
  markings,
  mark=at position #1 with {\fill circle (2.0pt);}},postaction={decorate}}} 

\usepackage{pgfplots}
\usepackage{pgfplotstable}
\pgfplotsset{compat=1.13}
\usepackage{amsmath}
\usepackage{subfigure}
\usepackage{epstopdf}
\usepackage{graphicx}
\usepackage{color}
\usepackage{cleveref}

\newcommand{\Om}{\Omega}
\newcommand{\R}{\mathbb{R}}

\newcommand{\I}{\mathbf{I}}
\newcommand{\pa}{\partial}
\newcommand{\diff}[1]{{\mathrm{d}{#1}}}

\newtheorem{theorem}{Theorem}

\newtheorem{remark}[theorem]{Remark}

\journal{Elsevier}

\begin{document}

\begin{frontmatter}


\title{High order treatment of moving curved boundaries: \\ 
Arbitrary-Lagrangian-Eulerian methods with\\a shifted boundary polynomials correction}

\author[lab1,dmi]{Walter Boscheri}
\author[lab2]{Mirco Ciallella\corref{cor1}}

\address[lab1]{Laboratoire de Math\'ematiques, Universit\'e Savoie Mont Blanc, CNRS, UMR 5127, Chamb\'ery, France}
\address[lab2]{Laboratoire Jacques-Louis Lions, Universit\'e Paris Cit\'e, CNRS, UMR 7598, Paris, France}
\address[dmi]{Department of Mathematics and Computer Science, University of Ferrara, Ferrara, Italy}

\cortext[cor1]{Corresponding author (\href{mailto:mirco.ciallella@u-paris.fr}{mirco.ciallella@u-paris.fr})}

\begin{abstract}

  In this paper we present a novel approach for the prescription of high order boundary conditions 
  when approximating the solution of the Euler equations for compressible gas dynamics on curved moving domains.
  When dealing with curved boundaries, the consistency of boundary conditions is a real challenge, 
  and it becomes even more challenging in the context of moving domains discretized with high order Arbitrary-Lagrangian-Eulerian (ALE) schemes. 
  The ALE formulation is particularly well-suited for handling moving and deforming domains, thus allowing for the simulation of 
  complex fluid-structure interaction problems. However, if not properly treated, the imposition of boundary conditions  
  can lead to significant errors in the numerical solution, which can spoil the high order discretization of the underlying mathematical model.
  In order to tackle this issue, we propose a new method based on the recently developed shifted boundary polynomial correction,
  which was originally proposed in the discontinuous Galerkin (DG) framework on fixed meshes. 
  The new method is integrated into the space-time corrector step of a direct ALE finite volume method to account for the local 
  curvature of the moving boundary by only exploiting the high order reconstruction polynomial of the finite volume control volume.
  It relies on a correction based on the extrapolated value of the cell polynomial evaluated at the true geometry, 
  thus not requiring  the explicit evaluation of high order Taylor series.
  This greatly simplifies the treatment of moving curved boundaries, as it allows for the use of standard simplicial meshes,
  which are much easier to generate and move than curvilinear ones, especially for 3D time-dependent problems.
  Several numerical experiments are presented demonstrating the high order convergence properties of the new method in the context 
  of compressible flows in moving curved domains, which remain approximated by piecewise linear elements.
  
\end{abstract}




\end{frontmatter}

\section{Introduction}\label{section:introduction}

The Arbitrary-Lagrangian-Eulerian (ALE) framework is a powerful tool for solving partial differential equations (PDEs)
on moving and deforming domains. This approach is particularly useful in the context of fluid dynamics and 
fluid-structure interaction problems, where the mesh motion is deformed to follow the displacement of 
the fluid and/or immersed bodies. Problems involving moving boundaries are common in many engineering applications,
and their numerical resolution opens up a wide range of possibilities for simulating complex physical phenomena.

In order to ensure conservation and accuracy in numerical solutions of moving domains, the direct ALE framework oringally introduced in \cite{boscheri2013arbitrary,boscheri2014direct} is a good candidate
for handling moving meshes while achieving higher convergence rates of discretization errors, and it allows a rather notable flexibility when compared to pure Lagrangian schemes \cite{munz1994godunov,carre2009cell,despres2003symmetrization,atallah2024weak,atallah2025high,kincl2025semi,Maire2007,Maire2009}. Indeed, in the ALE context \cite{kucharik2011hybrid,kucharik2012one,shashkov2008closure,berndt2011two} the mesh velocity is independent of the local fluid velocity, and in the direct ALE approach the governing equations are directly integrated over the moving space-time domain, hence satisfying by construction conservation of physical quantities and the compliance with the Geometric Conservation Law (GCL). In the direct ALE method forwarded in \cite{boscheri2013arbitrary,boscheri2014direct}, a piecewise linear approximation of the geometry is retained, while achieving arbitrary order of accuracy in space and time for the numerical solution.

However, ensuring consistency in boundary conditions presents a major challenge in high order discretization techniques. 
While these methods have demonstrated their ability to deliver highly accurate results with fewer degrees of 
freedom \cite{Wang_et_al_ijnmf13,cockburn1998runge,jiang1996efficient,dumbser2008unified,ar2017,Deconinck3}, handling boundaries remains a complex and unresolved issue.
This difficulty becomes even more pronounced when curved boundaries are involved, as accurately representing 
geometric features and maintaining simulation precision are crucial considerations.
Specifically, the high order approximation of boundary conditions is essential to achieve the expected convergence properties that
are affected by both the discretization of the PDE problem, and that of the domain geometry. 

A widely used strategy for addressing these challenges is the isoparametric approach \cite{bassi1997high}.
This method relies on high order polynomial reconstructions \cite{zienkiewic} of the same degree as the discretization scheme 
to approximate the targeted geometry. By doing so, it allows for a much more accurate representation of curved 
boundaries, enabling the use of larger elements compared to traditional linear approximations. 
Despite notable progress in this area, generating curved meshes remains significantly more complex than producing 
linear meshes, which has already reached a high level of maturity also thanks to its widespread development in both research and industry \cite{GMSH,arpaia2022h}.
The process of generating curvilinear meshes introduces several challenges \cite{luo2004automatic,FORTUNATO20161,sahni2010curved,dey2001towards}, 
including nonlinear mappings and the need for complex quadrature rules. 
These factors contribute to increased computational costs and complexity, particularly due to the higher number of required degrees of freedom and the more sophisticated algorithms involved.
Moreover, even with significant advancements in the development of numerical methods for handling curvilinear elements, 
producing high-quality curved meshes remains a difficult and multifaceted task, especially for realistic geometric configurations \cite{moxey2015isoparametric,coppeans2025aerodynamic}.

An alternative and computationally less expensive approach involves improving boundary conditions on simplicial meshes,
meaning that the discretization of the domain geometry is performed through simple linear elements, like triangles or tetrahedrons, 
but the boundary conditions are modified to achieve high order precision. 
However, this requires careful consideration of the local geometric properties in the affected regions.
In \cite{WangSun,krivodonova2006high} the first works in this direction were presented, where the authors proposed 
curvature corrected boundary conditions and their applications to high order schemes. 
These methods have been formulated only for slip wall boundary conditions and 2D geometries.

Some more recent works on the subject, to tackle general boundary conditions and 3D geometries, can be found in 
\cite{costa2018very,costa2019very,fernandez2020very,clain2021very,costa2021efficient,costa2022very,costa2023imposing,ciallella2024very}
where the modification in the prescription of boundary condition consists in replacing the polynomial reconstruction in each boundary cell with a constrained least squares problem.
In the finite volume framework \cite{costa2018very}, the minimization problem  extracts information from neighboring cells and imposes constraints tied to the physical boundary,
while in the discontinuous Galerkin framework \cite{ciallella2024very} the minimization problem is solved locally in each element using information coming from the internal degrees of freedom.
Differently, in \cite{ciallella2023shifted}, a direct approach is introduced to design boundary 
conditions through the extrapolation of the solution based on truncated Taylor expansions, 
similarly to  
\cite{Scovazzi1,Scovazzi3,carlier2023enriched,burman2018cut,ciallella2022extrapolated,visbech2025spectral,ciallella2020extrapolated,assonitis2022extrapolated},
that exploits the high order discontinuous Galerkin (DG) polynomial of the simplicial element to avoid the computation of local high order Taylor series.

In this paper, we investigate the extension of the shifted boundary polynomial correction initially developed in the DG framework on fixed grids \cite{ciallella2023shifted}
to the simulation of moving curved boundaries in the context of a direct high order ALE finite volume (FV) method \cite{boscheri2013arbitrary}. 
Specifically, we focus on the treatment of boundary conditions for the compressible Euler equations
using unstructured linear meshes to overcome the limitations given by standard boundary treatments, while greatly simplifying the 
computational complexity of treating 2D and 3D moving high order curvilinear meshes \cite{boscheri2016high}. 
This remarkable gain in terms of computational efficiency and algorithmic complexity is exploited to simplify the simulation 
of high order fluid-structure interaction problems \cite{gao2023high}. For the sake of simplicity in the presentation of the new method, 
we will focus on the interaction between a compressible inviscid fluid and geometries moving with a dictated velocity, 
but the proposed method has very high potential to be extended to more complex scenarios involving the interaction 
with deformable structures. With the new methods, we are able to show convergence rates of the isentropic Kidder problem \cite{kidder1976laser} up to fourth order of accuracy on moving tetrahedral meshes.

The paper is organized as follows. 
In \cref{section:model}, we present the mathematical model,
specifically the hyperbolic system of the Euler equations of gas dynamics. In \cref{section:ALE-FV}, we recall the high order ALE-FV method \cite{boscheri2013arbitrary,boscheri2014direct} developed for 2D and 3D unstructured linear meshes,
which achieves high order accuracy in time thanks to the space-time predictor described in \cref{subsection:predictor}.
In \cref{subsection:BC}, we explain how standard boundary conditions are implemented in our ALE-FV framework.
Then, in \cref{section:SBM}, we introduce the shifted boundary polynomial correction method used to improve the accuracy of boundary
conditions on moving meshes. In particular, we present the general methodology and its application in the context of moving slip wall
boundary conditions. In \cref{section:Results}, several numerical tests are shown to validate the proposed method
and to show its performance in terms of accuracy. In \cref{section:Conclusions}, we draw some conclusions and outline future perspectives.

\section{Mathematical model}\label{section:model}

The mathematical model considered in this work is the system of Euler equations for compressible 
gas dynamics in $d=2,3$ space dimensions reading:
\begin{align}\label{eq:euler0}
\frac{\partial \mathbf{Q}}{\partial t}+\nabla\cdot\mathbf{F}(\mathbf{Q}) =0, \quad &\text{in}\quad \Omega_T=\Omega\times[0,T]\subset\mathbb{R}^d\times\mathbb{R}^+, \\
\mathcal{B}(\mathbf{Q},\,\mathbf{Q}^{BC})=0, \quad &\text{on}\quad \partial\Omega_T=\partial\Omega\times[0,T]\subset\mathbb{R}^{d-1}\times\mathbb{R}^+, \nonumber
\end{align}
where $\mathcal{B}$ is the boundary condition operator, 
$\mathbf{Q}:\R^+\times\R^d\rightarrow\R^{d+2}$ is the vector of conserved variables and $\mathbf{F}:\R^{d+2}\rightarrow\R^{d+2}\times\R^{d}$ the nonlinear flux tensor, respectively defined as
\begin{equation}\label{eq:euler0a}
\mathbf{Q}=\left(
\begin{array}{c}
\rho \\ \rho \mathbf{u}  \\ \rho E
\end{array}
\right)\;,\;\;
\mathbf{F}(\mathbf{Q})= \left(
\begin{array}{c}
\rho \mathbf{u}\\ \rho \mathbf{u} \otimes \mathbf{u} + p \I \\ \rho H \mathbf{u} 
\end{array}
\right), 
\end{equation}
having denoted by $\rho$ the mass density, by $\mathbf{u}$ the velocity vector, by $p$ the pressure, 
and with $E=e+\mathbf{u}\cdot \mathbf{u}/2$ the specific total energy, $e$ being the specific internal 
energy. Finally, the specific total enthalpy is  $H=h+\mathbf{u}\cdot \mathbf{u}/2$ with  
$h=e+p/\rho$ the specific enthalpy. The relation between the pressure and the internal energy is 
given by the perfect gas equation of state
\begin{equation}\label{eq:EOS}
p=(\gamma-1)\rho e,
\end{equation}
where $\gamma$ is the constant ratio of specific heats ($\gamma=1.4$ for air).\\

The validation of the numerical methods are also performed by means of manufactured solutions 
on moving unstructured meshes which imply the discretization of an additional source term 
in~\cref{eq:euler0} which thus reads as follows
\begin{equation}\label{eq:euler1}
 \frac{\partial\mathbf{Q}}{\partial t}+\nabla\cdot\mathbf{F}\left(\mathbf{Q}\right) = \mathbf{S}(\mathbf{Q}).
\end{equation}
For this reason, we are going to include the algebraic source term $\mathbf{S}$ in the model discretization presented 
in~\cref{section:ALE-FV}.

\section{Numerical method}\label{section:ALE-FV}

\subsection{Computational domain and notation}\label{subsection:domain}
The time-dependent computational domain $\Om_t$ is discretized at the current time level $t^n$ 
by a total number $N_E$ of non-overlapping simplicial elements $T_i^n$. 
To simplify the description of the method, we describe here its implementation for 2D triangles, 
although the algorithm is also implemented and validated for 3D tetrahedrons. 
More details about the 3D implementation are given in \cite{boscheri2014direct}.
The union of all elements is referred to as the current tessellation $\mathcal{T}^n_\Om$ of $\Om_{t^n}=\Om^n$. It should be noticed that, in general, $\Om^n\neq\Om$ and in 
particular  $\pa\Om^n\neq\pa\Om$, except for very simple geometries with no curvature.

In a Lagrangian framework, the mesh moves and deforms as the solution evolves over time. 
To simplify the analysis, it is convenient to introduce a local reference coordinate system 
$\boldsymbol{\xi}=(\xi,\eta)$ with $(\xi,\eta)\in[0,1]$, where the reference element $T_e$ is defined. 
Each triangular element $T^n_i$ at the current time $t_n$ is mapped from the reference 
coordinates $(\xi, \eta)$ to the physical coordinates $(x, y)$ by a transformation relation: 
\begin{align}\label{eq:referencesystem}
x &= X^n_{1,i} + (X^n_{2,i} - X^n_{1,i})\xi + (X^n_{3,i} - X^n_{1,i})\eta, \nonumber\\
y &= Y^n_{1,i} + (Y^n_{2,i} - Y^n_{1,i})\xi + (Y^n_{3,i} - Y^n_{1,i})\eta,
\end{align}
where $\mathbf{X}^n_{k,i}=(X^n_{k,i},Y^n_{k,i})$ are the spatial coordinates of the $k$-th vertex 
of triangle $T^n_i$ at time $t^n$. The reference element $T_e$ is the triangle defined by
$\boldsymbol{\xi}_{1,e} = (0,0)$, $\boldsymbol{\xi}_{2,e} = (1,0)$, and $\boldsymbol{\xi}_{3,e} = (0,1)$.  

In the finite volume framework, the solution $\mathbf{Q}$ within each triangle $T^n_i$ is represented 
by spatial cell averages, which are given at time $t^n$ by
\begin{equation}\label{eq:cellaverage}
  \mathbf{Q}^n_i = \frac{1}{|T^n_i|} \int_{T^n_i} \mathbf{Q}(\mathbf{x},t^n) \, \diff{V},
\end{equation}
where $|T^n_i|$ denotes the volume of triangle $T^n_i$. 
High order accuracy is then achieved by reconstructing the solution using piecewise high order polynomials $\mathbf{w}_h(\mathbf{x},t^n)$
from a stencil of cell averages $\mathbf{Q}^n_i$, and approximating multi-dimensional integrals using consistent quadrature formulas.

\subsection{Piecewise WENO reconstruction on unstructured meshes}\label{subsection:WENO}

Herein we are going to use the WENO reconstruction technique in the polynomial formulation~\cite{dumbser2007arbitrary,dumbser2007quadrature}, 
rather than the classical pointwise method~\cite{jiang1999high,jiang1996efficient}. As we will see below, this will simplify notably the imposition
of the high order boundary conditions described in \cref{section:SBM}.
We limit us to provide the main points of the reconstruction procedure, and we refer the reader to the aforementioned references for more details on this subject.

The reconstruction polynomial of degree $M$ is built by using basis functions defined 
in the reference system $(\xi, \eta)$, and considering a reconstruction stencil $S^s_i$ of $n_s$ elements.
The stencil is defined to contain a greater number of elements than the optimal number of degrees of freedom 
$D=\prod_{\ell=1}^d (M+\ell)/\ell$ needed to reach the order of accuracy $M+1$ 
(in general, $n_s$ is chosen equal to $d \, D$).

The reconstruction polynomial is built involving a set of high order spatial basis functions 
$\psi_k(\boldsymbol{\xi})$ defined in the reference space for each candidate stencil $s$, and for each triangle $T^n_i$:
\begin{equation}
\mathbf{w}^s_h(\mathbf{x},t^n) = \sum_{k=1}^D \psi_k(\boldsymbol{\xi})  \mathbf{\hat w}^{n,s}_{k,i}(t) 
:= \psi_k(\boldsymbol{\xi})  \mathbf{\hat w}^{n,s}_{k,i},
\end{equation}
where the last definition is introduced to shorten the notation. 
In the following, we will use the tensor index notation with the Einstein summation convention, 
which implies summation over two equal indices. 
As commonly done in the context of finite volume methods, the polynomial coefficients are found by imposing
the integral conservation relation on each element $T^n_j \in S^s_i$, meaning that
\begin{equation}\label{eq:reconstructionsystem}
\frac{1}{|T_j^n|} \int_{T_j^n} \psi_k(\boldsymbol{\xi})  \mathbf{\hat w}^{n,s}_{k,i} \, \diff{V} = \mathbf{Q}^n_j, \qquad \forall T_j^n\in S^s_i .
\end{equation}
Since the size of the stencil is larger than the number of polynomial coefficients $D$, the system of equations \eqref{eq:reconstructionsystem}
is an overdetermined linear system and can be solved to find $\mathbf{\hat w}^{n,s}_{k,i}$ using a least-squares approach.
It should be noticed that, since the shape of the triangles evolves in time, the linear system \eqref{eq:reconstructionsystem} 
must be solved at every time step, while the choice of the stencils $S^s_i$ stays the same because we do not allow any topology change in the computational mesh.

The WENO polynomial is devised from the reconstruction polynomials computed in each stencil $S^s_i$.
To make the scheme nonlinear we introduce the classical oscillation indicators $\sigma_s$, as reported in \cite{jiang1996efficient}
\begin{equation}
\sigma_s = \Sigma_{km}  \mathbf{\hat w}^{n,s}_{k,i} \mathbf{\hat w}^{n,s}_{m,i}, \quad \text{with}\quad
\Sigma_{km} = \sum_{\alpha+\beta \leq M} \int_{T_e} \frac{\partial^{\alpha+\beta} \psi_k(\xi,\eta)}{\partial \xi^{\alpha} \partial \eta^{\beta}} \frac{\partial^{\alpha+\beta} \psi_m(\xi,\eta)}{\partial \xi^{\alpha} \partial \eta^{\beta}} \, \diff{\xi} \, \diff{\eta},
\end{equation}
and the nonlinear weights $\omega_s$, defined as
\begin{equation}
\tilde\omega = \frac{\lambda_s}{(\sigma_s+\varepsilon)^r}, \qquad \omega = \frac{\tilde\omega}{\sum_q \tilde\omega_q},
\end{equation}
where $\varepsilon=10^{-14}$, $r=8$, $\lambda_s =1$ for one-sided stencils and $\lambda_0 =10^5$ for the central stencil.
Finally, the WENO polynomial is given by
\begin{equation}
\mathbf{w}^s_h(\mathbf{x},t^n) = \psi_k(\boldsymbol{\xi})  \mathbf{\hat w}^{n}_{k,i}, \qquad  \mathbf{\hat w}^{n}_{k,i} = \sum_s \omega_s \mathbf{\hat w}^{n,s}_{k,i}.
\end{equation}

\subsection{Local space-time predictor on moving meshes}\label{subsection:predictor}

The high order accuracy in time for moving meshes is achieved through a local space-time predictor $\mathbf{q}_h(\mathbf{x},t)$ 
within each element $T_i(t)$ in the temporal slab $[t^n,t^{n+1}]$, which does not require any neighbor information. 
Notice that, due to the moving mesh, the space-time element will be characterized by different triangle shapes at time $t^n$ and $t^{n+1}$. 
The computation of $\mathbf{q}_h(\mathbf{x},t)$ is performed through the local evolution in time of the polynomial $\mathbf{w}^s_h(\mathbf{x},t^n)$
with a weak space-time formulation of \cref{eq:euler0}. This idea finds his roots in \cite{dumbser2008unified} and was then
developed in the context of moving meshes in 
\cite{dumbser2013arbitrary,boscheri2013arbitrary,boscheri2014direct,boscheri2015direct,gaburro2020high}.
For this reason, we here describe the main steps while further details can be found in the
aforementioned references.

The space-time formulation is introduced through the modified reference systems that also
include the time variable. In particular, $\mathbf{\tilde x}=(x,y,t)$ is the physical coordinate vector,
while $\boldsymbol{\tilde \xi}=(\xi,\eta,\tau)$ is the  reference one, with $\tau \in [0,1]$.
The predictor is then defined through local space-time basis functions $\theta(\boldsymbol{\tilde \xi})$ in each element $T_i(t)$, that is
\begin{equation}
\mathbf{q}_h(\mathbf{x},t) = \sum_{k=1}^Q \theta_k(\boldsymbol{\tilde \xi}) \mathbf{\hat q}_{k,i} 
:=\theta_k(\boldsymbol{\tilde \xi}) \mathbf{\hat q}_{k,i}, \quad \text{where} \quad Q =\prod_{\ell=1}^{d+1} (M+\ell)/\ell.
\end{equation}
The same space-time representation is also used for fluxes $\mathbf{F}=(\mathbf{f},\mathbf{g})$ and source terms $\mathbf{S}$ by exploiting the interpolation property of the chosen nodal basis:
$$ \mathbf{f}_h(\mathbf{x},t) = \theta_k(\boldsymbol{\tilde \xi}) \mathbf{\hat f}_{k,i}, \quad 
\mathbf{g}_h(\mathbf{x},t) = \theta_k(\boldsymbol{\tilde \xi}) \mathbf{\hat g}_{k,i}, \quad 
\mathbf{S}_h(\mathbf{x},t) = \theta_k(\boldsymbol{\tilde \xi}) \mathbf{\hat S}_{k,i}. $$
By introducing the mapping in time $t=t^n + \tau \Delta t$, the space-time weak formulation
can be written in the reference space-time element $T_e\times[0,1]$.
The space-time transformation Jacobian and its inverse read
\begin{equation}
  J = \frac{\pa \mathbf{\tilde x}}{\pa \boldsymbol{\tilde \xi}} = \left(
    \begin{array}{ccc}
    x_\xi & x_\eta & x_\tau \\ y_\xi & y_\eta & y_\tau  \\ 0 & 0 & \Delta t
    \end{array}
    \right), \quad 
  J^{-1} = \frac{\pa \boldsymbol{\tilde \xi}}{\pa \mathbf{\tilde x}} = \left(
    \begin{array}{ccc}
    \xi_x & \xi_y & \xi_t \\ \eta_x & \eta_y & \eta_t  \\ 0 & 0 & \frac{1}{\Delta t}
    \end{array}
    \right).
\end{equation}
The local reference system with the space-time Jacobian $J$ are used to rewrite 
\cref{eq:euler0} as
\begin{equation}\label{eq:pdereferencespacetime}
  \frac{\pa \mathbf{Q}}{\pa \tau} + \Delta t \left(\frac{\pa \mathbf{Q}}{\pa \xi}\xi_t + \frac{\pa \mathbf{Q}}{\pa \eta}\eta_t + \frac{\pa \mathbf{f}}{\pa \xi}\xi_x + \frac{\pa \mathbf{f}}{\pa \eta}\eta_x + \frac{\pa \mathbf{g}}{\pa \xi}\xi_y + \frac{\pa \mathbf{g}}{\pa \eta}\eta_y \right) = \Delta t \mathbf{S}(\mathbf{Q}).
\end{equation}
To simplify the notation, we introduce the following operator
$$\langle f,g \rangle = \int_0^1 \int_{T_e} f(\xi,\eta,\tau) g(\xi,\eta,\tau) \, \diff{\xi} \, \diff{\eta} \, \diff{\tau}.$$
By plugging $\mathbf{q}_h$, $\mathbf{f}_h$, $\mathbf{g}_h$, and $\mathbf{S}_h$ into \cref{eq:pdereferencespacetime}, 
and then integrating it over the space-time reference element  $T_e\times[0,1]$, we obtain the following compact weak formulation:
\begin{equation}
\mathbf{K}_\tau \mathbf{\hat q}_{k,i} + \Delta t \left( \mathbf{K}_t \mathbf{\hat q}_{k,i} + \mathbf{K}_x \mathbf{\hat f}_{k,i} + \mathbf{K}_y \mathbf{\hat g}_{k,i}  \right) = \Delta t \mathbf{M}_t \mathbf{\hat S}_{k,i}, \label{eqn.pdeweak}
\end{equation}
where
\begin{align*}
\mathbf{K}_\tau = \bigg\langle \theta_l, \frac{\pa \theta_k}{\pa \tau} \bigg\rangle, 
\quad \mathbf{M} =  \langle \theta_l, \theta_k \rangle, 
\quad \mathbf{K}_t = \bigg\langle \theta_l, \frac{\pa \theta_k}{\pa \xi} \xi_t \bigg\rangle + \bigg\langle \theta_l, \frac{\pa \theta_k}{\pa \eta} \eta_t \bigg\rangle, \\ 
\quad \mathbf{K}_x = \bigg\langle \theta_l, \frac{\pa \theta_k}{\pa \xi} \xi_x \bigg\rangle + \bigg\langle \theta_l, \frac{\pa \theta_k}{\pa \eta} \eta_x \bigg\rangle,
\quad \mathbf{K}_y = \bigg\langle \theta_l, \frac{\pa \theta_k}{\pa \xi} \xi_y \bigg\rangle + \bigg\langle \theta_l, \frac{\pa \theta_k}{\pa \eta} \eta_y \bigg\rangle.
\end{align*}
It should be noticed that all the above matrices can be computed on the reference space-time element and stored once and for all in the pre-processing step. Only the inverse of the space-time Jacobian $J^{-1}$ must be evaluated at each time step in order to take into account the mesh motion. The coefficients of the predictor $\mathbf{\hat q}_{k,i}$ in each element $T_i(t)$ 
can be found independently for each cell through an iterative procedure, since the resulting system given by \cref{eqn.pdeweak} is nonlinear. More details are given in \cite{boscheri2013arbitrary}.

Since the mesh is also moving, we have to consider the evolution of the vertex coordinates of the space-time element,
whose motion is described by the trajectory equation
\begin{equation}\label{eq:meshspeed}
\frac{\diff{\mathbf{x}}}{\diff{t}} = \mathbf{V}(\mathbf{x},t),
\end{equation}
where $\mathbf{V}(\mathbf{x},t)$ is the local mesh velocity, which can be independent from the local fluid velocity.
Following the same space-time nodal expansion, the mesh position and velocity in each space-time element $T_i(t)$ are expressed as
$\mathbf{x}=\theta_k(\boldsymbol{\tilde \xi}) \mathbf{\hat x}_{k,i}$ and $\mathbf{V}_h=\theta_k(\boldsymbol{\tilde \xi}) \mathbf{\hat V}_{k,i} $.

As done in \cite{dumbser2013arbitrary}, \cref{eq:meshspeed} can be solved in a space-time finite element fashion:
\begin{equation}
  \mathbf{K}_\tau \mathbf{\hat x}_{k,i} = \Delta t \mathbf{M} \mathbf{\hat V}_{k,i}.
\end{equation}
This formulation is employed to obtain a local space-time predictor for the mesh nodal coordinates.

To guarantee continuity of the discretized geometry, the mesh velocity must be uniquely determined at vertex $\nu$, hence requiring a so-called nodal solver, since the mesh velocity is defined by the local predictor velocity within all elements surrounding the vertex. Here, a unique vertex velocity is computed by simply averaging all the element contributions obtaining $\mathbf{\bar V}^n_\nu$. 
Finally, vertex $\nu$ can be moved according to $\mathbf{X}^{n+1}_\nu = \mathbf{X}^{n}_\nu + \Delta t \mathbf{\bar V}^n_\nu $.

\subsection{Mesh motion}\label{subsection:meshmotion}

In this work, we consider that the computational domain is moving according to a prescribed velocity at the boundary,
which is a common assumption in many applications related to fluid-structure interaction.
In order to expand coherently the mesh motion from the Lagrangian boundary to a fluid grid, 
there are many possible techniques in the literature 
\cite{lohner1996improved,helenbrook2003mesh,dwight2009robust,stein2003mesh}.
One of the most common and straightforward approaches to implement is the Harmonic equation, which 
consists in solving the following Laplacian problem for the mesh velocity:
\begin{equation}\label{eq:laplace}
  \begin{cases}
  -\Delta \mathbf{V}(\mathbf{x},t) &= 0, \quad \text{in}\quad \Omega, \quad \\
  \mathbf{V}(\mathbf{x},t) &= \mathbf{V}^{b}(\mathbf{x},t), \quad \text{on}\quad \partial\Omega,
  \end{cases}
\end{equation}
where $\mathbf{V}^{b}(\mathbf{x},t)$ is the prescribed velocity at the boundary $\partial\Omega$.
Herein, the solution of \cref{eq:laplace} is computed by means of a standard $\mathcal{P}_1$ finite element method
on triangles and tetrahedrons.

\subsection{High order Finite Volume scheme}\label{subsection:FV}

Let \ref{eq:euler0} be written in a space-time divergence form \cite{boscheri2013arbitrary} as
\begin{equation}
\tilde\nabla \cdot \tilde{\mathbf{Q}} = \mathbf{S}(\mathbf{Q}), \quad \text{where} \quad 
\tilde\nabla = \left(\frac{\pa}{\pa t},\frac{\pa}{\pa x},\frac{\pa}{\pa y}\right)^T, \quad 
\tilde{\mathbf{Q}} = \left( \mathbf{Q}, \mathbf{F}\right).  
\end{equation}
Integrating over the space-time control volume $C_i^n = T_i(t)\times[t^n,t^{n+1}]$, and applying the divergence theorem gives
\begin{equation}
\int_{\pa C_i^n} \tilde{\mathbf{Q}} \cdot \mathbf{\tilde n} \, \diff{S} = \int_{C_i^n} \mathbf{S}(\mathbf{Q}) \, \diff{V},
\end{equation}
where $\mathbf{\tilde n} = (\tilde n_x,\tilde n_y,\tilde n_t)$ is the outward pointing space-time normal on the surface $\pa C_i^n$.
The space-time surface $\pa C_i^n$ includes $T_i^n$, $T_i^{n+1}$ and the three shared surfaces with neighboring triangles. 
Given that, for $T_i^n$ and $T_i^{n+1}$, the space-time normals respectively read $\mathbf{\tilde n} = (0,0,-1)$ and $\mathbf{\tilde n} = (0,0,1)$,
the ALE one-step finite volume scheme can be written as
\begin{equation}\label{eq:alefv}
|T_i^{n+1}| \mathbf{Q}_i^{n+1} = |T_i^{n}| \mathbf{Q}_i^{n} - \sum_{j\in\mathcal N_i}\int_{S_{ij}^n} \tilde{\mathbf{G}}_{ij} \cdot \mathbf{\tilde n}_{ij} \, \diff{S} + \int_{C_i^n} \mathbf{S}(\mathbf{q}_h)\, \diff{V} ,
\end{equation}
where $\bigcup_{j\in\mathcal N_i} S^n_{ij} = \pa C_i^n\setminus (T_i^{n} \cup T_i^{n+1})$, and $\mathcal N_i$ denotes the
neighborhood of triangle $T_i(t)$.

The term inside the surface integral in \cref{eq:alefv} represents the ALE numerical flux at the space-time surface $S_{ij}^n$.
In this work, we consider the Osher ALE flux developed before in the Eulerian \cite{dumbser2011universal} and then in the ALE framework \cite{dumbser2013arbitrary}:
\begin{equation}\label{eq:osherale}
  \tilde{\mathbf{G}}_{ij} \cdot \mathbf{\tilde n}_{ij}(\mathbf{q}^-,\mathbf{q}^+) = \frac12 \left( \tilde{\mathbf{Q}}(\mathbf{q}^+) + \tilde{\mathbf{Q}}(\mathbf{q}^-) \right)\cdot \mathbf{\tilde n}_{ij}  - \frac12 \left(\int_0^1 |\mathbf{A}_{\mathbf{n}}^{\mathbf{V}}(\boldsymbol\Psi(s))| \diff{s}\right) \left(\mathbf{q}^+-\mathbf{q}^-\right) ,
\end{equation}
with the linear path connecting the left and right states, as $\boldsymbol\Psi(s) = \mathbf{q}^+ + s \left(\mathbf{q}^+-\mathbf{q}^-\right)$ with $0\leq s\leq 1$.  
The matrix $\mathbf{A}_{\mathbf{n}}^{\mathbf{V}}$ represents the ALE Jacobian in the spatial normal direction, defined as
\begin{equation}
  \mathbf{A}_{\mathbf{n}}^{\mathbf{V}}(\mathbf{Q}) = \frac{\pa(\mathbf{F}\cdot\mathbf{n})}{\pa\mathbf{Q}} - (\mathbf{V}\cdot\mathbf{n})\I, 
  \quad \text{where}\quad \mathbf{n} = \frac{(\tilde n_x,\tilde n_y)^T}{\sqrt{\tilde n_x^2+\tilde n_y^2}}.
\end{equation}
In \cref{eq:osherale}, we used the common definition of matrix absolute value:
$$|\mathbf{A}| = \mathbf{R}|\mathbf{\Lambda}|\mathbf{R}^{-1}, \quad\text{where}\quad |\mathbf{\Lambda}| = \text{diag}(|\lambda_1|,...,|\lambda_{d+2}|),$$
where $\mathbf{R}$ is the matrix of right eigenvectors. Notice that the path integral in \cref{eq:osherale} is approximated using Gaussian quadrature rules.

\subsection{Boundary conditions}\label{subsection:BC}

When the boundary $\partial T_i$ of element $T_i$ belongs to $\partial\Omega^n$, the numerical flux must  
account for the appropriate  boundary conditions. The space-time flux consistent with such conditions will be denoted 
by $\tilde{\mathbf{G}}^{BC}_{ij} \cdot \mathbf{\tilde n}_{ij}$ and will be computed through a ghost state $\mathbf{Q}^{BC}$. 
Meaning that the neighboring reconstructed state is replaced by the ghost state in \cref{eq:osherale}, as 
$\tilde{\mathbf{G}}^{BC}_{ij} \cdot \mathbf{\tilde n}_{ij}(\mathbf{q}^-,\mathbf{Q}^{BC}).$
In this work, we focus on two different boundary conditions to be enforced: 
\begin{itemize}
\item Dirichlet-type BC can be enforced weakly through fluxes when the whole state $\mathbf{Q}^{BC}$ is set to prescribed values;
\item slip walls with non-zero speed $\mathbf{w}$ can be imposed through the constraint that 
$\mathbf{u}\cdot\mathbf{n}=\mathbf{w}\cdot\mathbf{n}$ at the boundary. 
In the simple situation where the wall does not move, i.e.\ $\mathbf{u}\cdot\mathbf{n}\equiv 0$, it is possible to show that $\mathbf{Q}^{BC}$ has the same density, 
internal energy and tangential velocity of $\mathbf{q}^-$, and the opposite normal velocity component.
In particular, considering the primitive variables, we can set $\rho^{BC}=\rho^-$ and $p^{BC}=p^-$, while 
\begin{equation}\label{eq:slipwall}
\mathbf{u}^{BC} = \mathbf{u}^- - 2 (\mathbf{u}^-\cdot\mathbf{n})\mathbf{n}.  
\end{equation}
Note that the last relation can be recovered from the more general one that consider a moving wall, that reads 
\begin{equation}\label{eq:slipwall2}
\mathbf{u}^{BC} = 2 \mathbf{u}^b - \mathbf{u}^- .
\end{equation}
where $\mathbf{u}^b$ represents the velocity one would impose at the boundary. 
Therefore, in the specific case of fixed walls, we can easily recover 
equation \eqref{eq:slipwall} from \eqref{eq:slipwall2}, 
given that we would like to have $\mathbf{u}^b=(\mathbf{u}^-\cdot\boldsymbol{\tau})\boldsymbol{\tau}$, and
\begin{align*}
  \mathbf{u}^{BC} &= 2 (\mathbf{u}^-\cdot\boldsymbol{\tau})\boldsymbol{\tau} - \mathbf{u}^- \\
                  &= (\mathbf{u}^-\cdot\boldsymbol{\tau})\boldsymbol{\tau} - (\mathbf{u}^-\cdot\mathbf{n})\mathbf{n} \\
                  &= \mathbf{u}^- - 2 (\mathbf{u}^-\cdot\mathbf{n})\mathbf{n} ,
\end{align*}
where we considered that $\mathbf{u}^-=(\mathbf{u}^-\cdot\mathbf{n})\mathbf{n} + (\mathbf{u}^-\cdot\boldsymbol{\tau})\boldsymbol{\tau}$,
with $\boldsymbol{\tau}$ being the unit vector tangent to the wall.
In the case of moving slip walls, the velocity $\mathbf{u}^b$ can be set starting from the wall velocity $\mathbf{w}$;
in particular we should have
\begin{equation}\label{eq:slipwall3}
  \mathbf{u}^b = (\mathbf{w}\cdot\mathbf{n})\mathbf{n} + (\mathbf{u}^-\cdot\boldsymbol{\tau})\boldsymbol{\tau} ,
\end{equation}
which is then used in combination with \eqref{eq:slipwall2} to compute the ghost state $\mathbf{Q}^{BC}$.

\end{itemize}
It is clear that the computation of the ghost state $\mathbf{Q}^{BC}$ strongly depends on the position and normal at the discretized boundary,
which only represents the discretized counterpart of the exact data.
For high order methods, achieving a truly high level of accuracy relies on the simultaneous control of errors in 
both geometry representation and flow variables. This also involves the use of consistent quadrature rules for 
approximating boundary integrals. A straightforward yet computationally demanding approach is to employ a 
curved high order approximation of the domain boundary. This typically requires an isoparametric representation of 
the boundary and the generation of a valid curved volume mesh~\cite{Wang_et_al_ijnmf13}. Such an approach has been pursued in the context of direct ALE schemes in \cite{boscheri2016high}.

\begin{figure}[!h]
  \centering
  \includegraphics[width=0.6\textwidth]{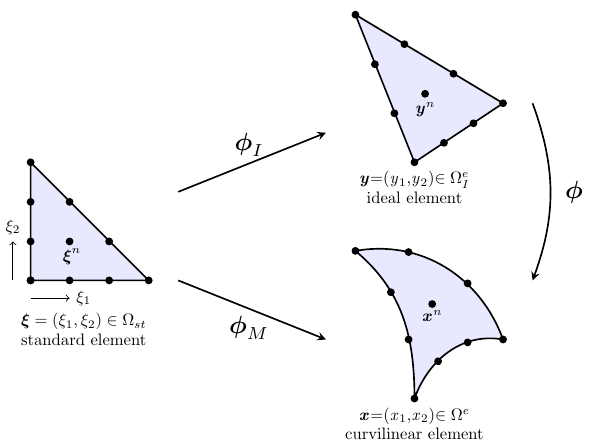}
  \caption{This illustration provides a visual representation of the necessary transformations when working with curved meshes.
  First, a mapping is required from the standard reference element, $\Omega_{st}$ (left), to the straight-sided element, $\Omega^e_{I}$ (top right), denoted as $\boldsymbol\phi_I:\Omega_{st}\rightarrow\Omega^e_I$.
  Similarly, a mapping to the curvilinear element (bottom right) is needed, represented as $\boldsymbol\phi_M:\Omega_{st}\rightarrow\Omega^e$.
  Finally, the deformation mapping $\boldsymbol\phi:\Omega^e_I\rightarrow\Omega^e$ is obtained by composing these transformations as $\boldsymbol\phi=\boldsymbol\phi_M\circ\boldsymbol\phi_I^{-1}$.}
  \label{fig:curvedElement}
\end{figure}

Curvilinear grids provide significantly more precise geometric boundaries, enabling the use of larger elements 
compared to linear ones. Various techniques exist for generating high order meshes, including curving pre-existing 
linear meshes~\cite{dey2001towards,luo2004automatic,sahni2010curved,FORTUNATO20161,MOXEY2016130} or employing 
optimization and variational methods~\cite{https://doi.org/10.1002/nme.4888,TOULORGE20138,TURNER201873}. 
Regardless of the approach, defining and evaluating mappings between curvilinear and reference elements 
(see \cref{fig:curvedElement}) remains necessary. While advancements have significantly refined linear mesh 
generation for complex geometries, the robust creation of curved meshes continues to be a challenging task.

In the following section, we present an alternative approach given by the boundary polynomial correction method, which allows us to
bypass the need for generating curved meshes, and still recover high order accuracy with moving curved domains.

\section{High order boundary conditions on moving linear meshes}\label{section:SBM}

In this work, we seek to overcome the need for generating curved meshes by directly utilizing conformal 
linear meshes. Our approach is based on the boundary polynomial correction introduced in \cite{ciallella2023shifted}. 
By leveraging this method, we aim to compensate for geometric errors while preserving the intended accuracy.

For clarity, we briefly introduce the notation used for boundaries. 
Let $\Omega_h$ denote a linear conformal mesh discretizing the physical domain $\Omega$, 
and let $\partial \Omega_h$ be the piecewise affine approximation of the exact curved boundary $\partial \Omega$. 
We refer to $\partial \Omega_h$ as the discretized boundary. 
Additionally, we assume that each point on $\partial \Omega_h$ can be mapped to a unique corresponding point 
on the true boundary $\partial \Omega$, such that:  
\begin{equation*}
\forall \mathbf{\tilde{x}} \in \partial \Omega_h, \quad \exists \mathbf{x} \in \partial \Omega \quad 
\text{such that} \quad \mathbf{x} = \mathcal{M}(\mathbf{\tilde{x}}),	
\end{equation*}
where $\mathcal{M}$ represents the mapping function, for instance in this work we rely on a signed distance function defined using distances
along the exact boundary normals $\mathbf{n}$:
\begin{equation}
  \mathbf{d}(\mathbf{\tilde{x}}) = \mathbf{x} - \mathbf{\tilde{x}} = (\mathcal{M}-\I) (\mathbf{\tilde{x}}),\quad\text{with}\quad
  \mathbf{x} = \mathbf{\tilde{x}} + d(\mathbf{\tilde x}) \mathbf{n},
\end{equation}
where $d(\mathbf{\tilde x})=\|\mathbf{d}(\mathbf{\tilde{x}})\|$.

We provide a detailed explanation of the boundary correction for handling Dirichlet boundary conditions. 
Let $\phi_D$ denote the prescribed variable value (analogous expressions can be formulated for all variables).
In the same spirit of \cite{huffman1993practical,burman2018cut,Scovazzi1}, the considered boundary condition 
on $\pa \Omega$ can be approximated on $\pa \Omega_h$, under some smoothness hypothesis, by the following Taylor expansion
\begin{equation}\label{eq:taylorexpansion}
\phi(\mathbf{x}) = \phi(\mathbf{\tilde{x}}) + 
\sum_{k=1}^M d^k n_{i_1} n_{i_2} \dots n_{i_k} \frac{\partial^k \phi}{\partial x_{i_1} \partial x_{i_2} \dots \partial x_{i_k}} \Big|_{\mathbf{\tilde{x}}},
\end{equation}
where $n_{i_1}, n_{i_2}, \dots, n_{i_k}$ are the components of the normal vector $\mathbf{n}$,
and $\frac{\partial^k \phi}{\partial x_{i_1} \partial x_{i_2} \dots \partial x_{i_k}}$ 
represents the $k$-th order partial derivative of $\phi$, evaluated at $\mathbf{\tilde{x}}$.

The approach involves adjusting the boundary condition on $\pa\Om_h$ to incorporate all corrective terms. 
This effectively leads to using a modified prescribed value given by
\begin{equation}\label{eq:boundarycorrection}
\phi^\star(\mathbf{\tilde{x}}) = \phi_D(\mathbf{x}) - \sum_{k=1}^M d^k n_{i_1} n_{i_2} \dots n_{i_k} \frac{\partial^k \phi}{\partial x_{i_1} \partial x_{i_2} \dots \partial x_{i_k}} \Big|_{\mathbf{\tilde{x}}}.
\end{equation}
It is important to note that all derivative terms are computed based on the available polynomial approximation of $\phi$
and are evaluated at the appropriate quadrature points.  
As the order of accuracy increases, the correction terms in \cref{eq:boundarycorrection} become increasingly complex and 
computationally expensive, particularly in three space dimensions. To overcome this problem, we consider here an alternative formulation 
\cite{ciallella2023shifted} that skips the explicit evaluation of these terms. 

Starting from the Taylor expansion \cref{eq:taylorexpansion}, we can deduce that 
\begin{equation}
  \phi(\mathbf{x}) - \phi(\mathbf{\tilde{x}}) = 
  \sum_{k=1}^M d^k n_{i_1} n_{i_2} \dots n_{i_k} \frac{\partial^k \phi}{\partial x_{i_1} \partial x_{i_2} \dots \partial x_{i_k}} \Big|_{\mathbf{\tilde{x}}}.  
\end{equation}
In practice, this means that by simply evaluating the polynomial of $\phi$ at $\mathbf{x}=\mathbf{\tilde{x}}+d(\mathbf{\tilde{x}})\mathbf{n}$,
and taking the difference from its value at $\mathbf{\tilde{x}}$, all correction terms can be efficiently obtained in a 
single step. This remarkably simplifies the computation. 
Moreover, since the required data at quadrature points can always be derived from the same basis used in the approximation of the numerical solution, 
the approach remains universally applicable regardless of the chosen basis function. 
The use of linear meshes further enhances this process by ensuring that the mapping is inherently linear, 
eliminating any ambiguity when transitioning from reference to physical space, regardless of the approximation degree.

The modified boundary condition given in \cref{eq:boundarycorrection} can be recast more elegantly as
\begin{equation}\label{eq:boundarycorrection2}
  \phi^\star(\mathbf{\tilde{x}}) = \phi_D(\mathbf{x}) - \left[ \phi(\mathbf{x}) - \phi(\mathbf{\tilde{x}}) \right] .
\end{equation}
This approach drastically simplifies the computation, requiring just one additional polynomial evaluation. 
Its implementation is particularly straightforward for straight-sided simplex elements, 
where basis functions are defined in either reference or physical space. 
As discussed earlier, the extrapolated variables obtained through this process are then used to define the ghost state 
$\mathbf{Q}^{BC}$. 
Equation \ref{eq:boundarycorrection2}, applied to either conservative or primitive variables, can be straightforwardly used 
to compute the ghost state $\mathbf{Q}^{BC}$ for Dirichlet-type boundary conditions.

\begin{remark}[Polynomial correction evaluation in the ALE framework]
In order to simplify the computation of integrals in the ALE formulation, everything is evaluated in the reference space-time element $T_e\times[0,1]$.
Therefore, when computing the boundary integrals, the polynomial correction must be also evaluated for each quadrature point 
of the reference space-time element. This can be easily done by finding the point $\mathbf{x}$ in the physical space, as mentioned above,
and then using the inverse Jacobian transformation to find the corresponding point in the reference space-time element, according to the mapping \cref{eq:referencesystem}, which is linear, hence allowing this transformation to be carried out very easily.
\end{remark}

\subsection{Moving slip walls}\label{subsection:slipwalls}

\begin{figure}%
  \centering
  \subfigure[Boundaries]{\includegraphics[width=0.45\textwidth]{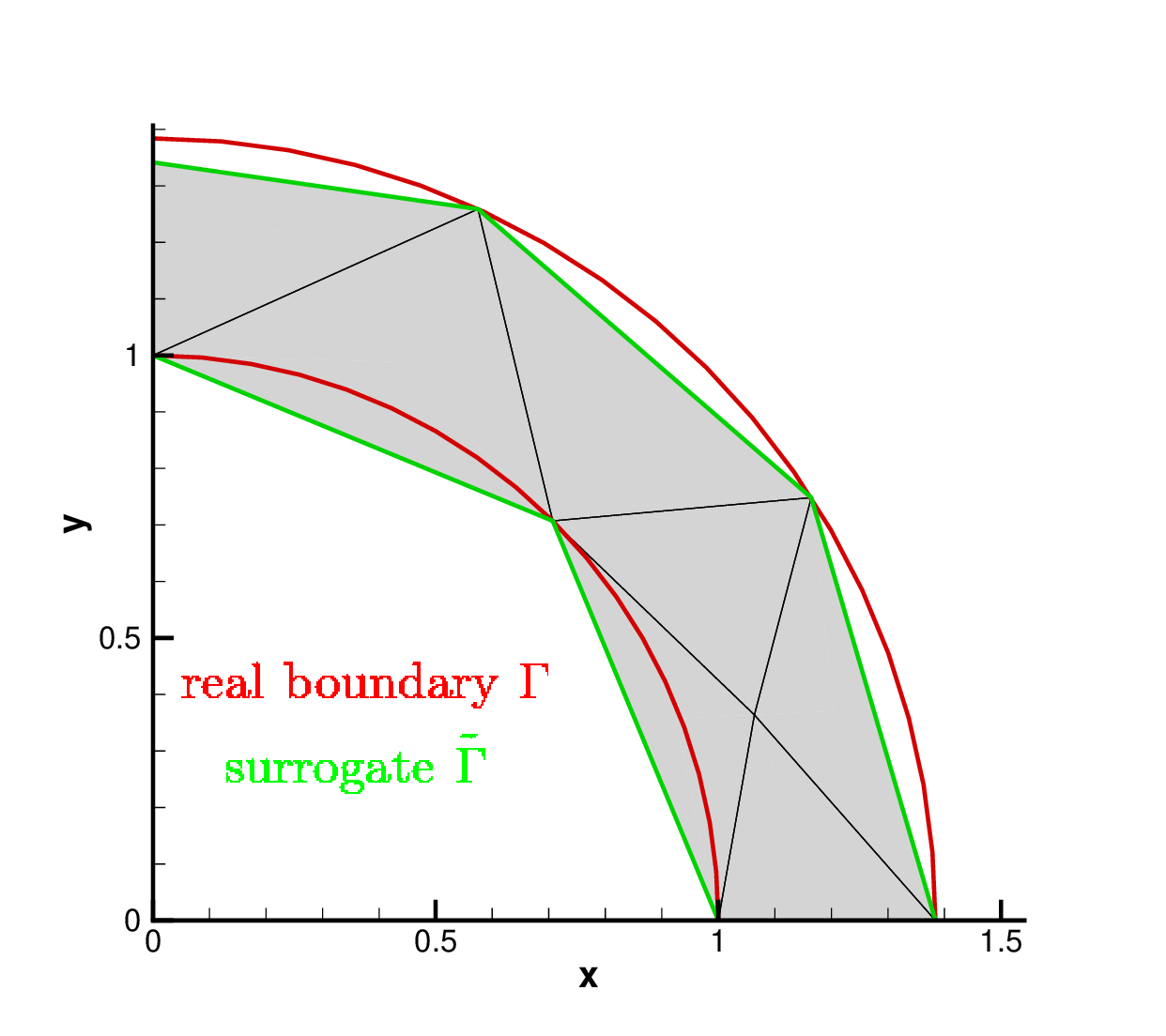}}\quad
  \subfigure[Normals]{\includegraphics[width=0.35\textwidth]{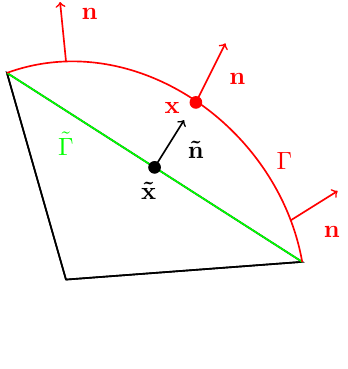}}
  \caption{Visual representation of the real and discretized (surrogate) boundaries for curved geometries, with correspondent normals.}
  \label{fig:boundarynormals}%
  \end{figure}

In the context of slip walls, when the boundary is discretized with linear meshes, the degradation of the order of accuracy
also experienced in previous works \cite{krivodonova2006high,ciallella2023shifted}, is related to the low accuracy of 
boundary normals, which are simply taken as the edge normals of the linear boundary mesh (see figure \ref{fig:boundarynormals}).
Indeed, when no boundary correction is applied, the velocity at the boundary is prescribed as
$  \mathbf{u}^b = (\mathbf{w}\cdot\mathbf{\tilde n})\mathbf{\tilde n} + (\mathbf{u}^-\cdot\boldsymbol{\tilde\tau})\boldsymbol{\tilde\tau} $, 
where $\mathbf{\tilde n}$  and $\boldsymbol{\tilde\tau}$ are the normal and tangent vectors, respectively, of the linear boundary mesh.

Contrary to what was done in \cite{ciallella2023shifted}, in this work, we propose to apply the polynomial correction \eqref{eq:boundarycorrection2} 
to improve the accuracy of moving slip walls, discretized through linear meshes, by directly modifying the ghost state to take 
into account the real geometry and correspondent normals.
Therefore, the correction \eqref{eq:boundarycorrection2} can be applied starting from \eqref{eq:slipwall2} and \eqref{eq:slipwall3}.
In this case, the ghost state $\mathbf{Q}^{BC}$ is determined by applying the polynomial correction to the normal velocity field to compute
$  \mathbf{u}^b = (\mathbf{w}^{\star}\cdot\mathbf{n})\mathbf{n} + (\mathbf{u}^-\cdot\boldsymbol{\tau})\boldsymbol{\tau} $, 
where the corrected term is evaluated through the following correction:
\begin{equation}\label{eq:slipwallSBM}
  (\mathbf{w}^{\star}\cdot\mathbf{n})(\mathbf{\tilde x}) = (\mathbf{w}\cdot\mathbf{n})(\mathbf{x}) - [\mathbf{u}(\mathbf{x}) - \mathbf{u}(\mathbf{\tilde x})] \cdot\mathbf{n}. 
\end{equation}
Notice that in this work, we consider a prescribed velocity at the boundary, and therefore the
term $(\mathbf{w}\cdot\mathbf{n})(\mathbf{x})$ can be computed explicitly given the local boundary normal $\mathbf{n}$.
  
\section{Numerical experiments}\label{section:Results}

In this section, we will assess the performance of the proposed high order boundary conditions on 2D and 3D moving linear meshes.
All the simulations will be performed using the same high order ALE-FV scheme described in \cref{section:ALE-FV},
while only modifying the imposition of boundary conditions to take into account the curved boundary.
To verify the high order accuracy of the method, we will consider particular manufactured solutions, for both 2D and 3D meshes,
that are characterized by a moving curved domain, whose exact solution is known.
Then, we will simulate a more complex test case often studied in the literature of Lagrangian hydrodynamics, 
which is the Kidder problem \cite{kidder1976laser}, for both 2D and 3D meshes, for which exact solutions are also available.
In particular, high order convergence analysis for the Kidder problem  are here presented in both 2D and 3D computational domains.
We would like to stress the fact that, up to our knowledge, this is the first time that a high order convergence analysis
(higher than second order) is presented for the Kidder problem in the literature.
We will conclude with two qualitative test cases, where the high order boundary conditions will be used to enforce moving slip
wall condition for the fluid-structure interaction of oscillating cylinders.  
For all cases, the mesh motion is prescribed at the boundaries, and the mesh velocity at the internal nodes 
is computed as described in \cref{subsection:meshmotion}. 

\subsection{Manufactured solution on 2D moving meshes}\label{subsection:Manufactured2D}

\begin{figure}
\centering
\subfigure[]{\includegraphics[width=0.48\textwidth]{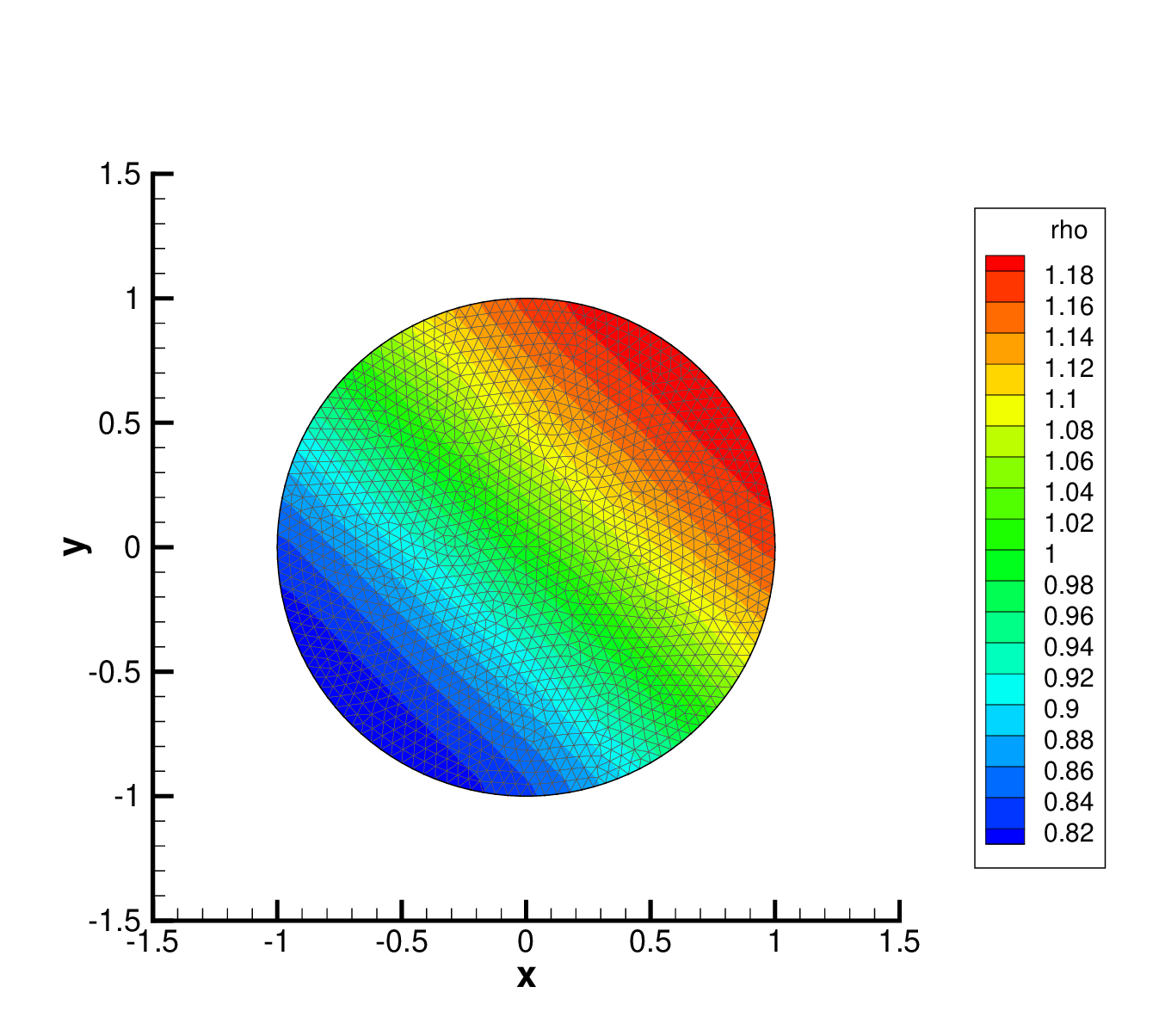}}
\subfigure[]{\includegraphics[width=0.48\textwidth]{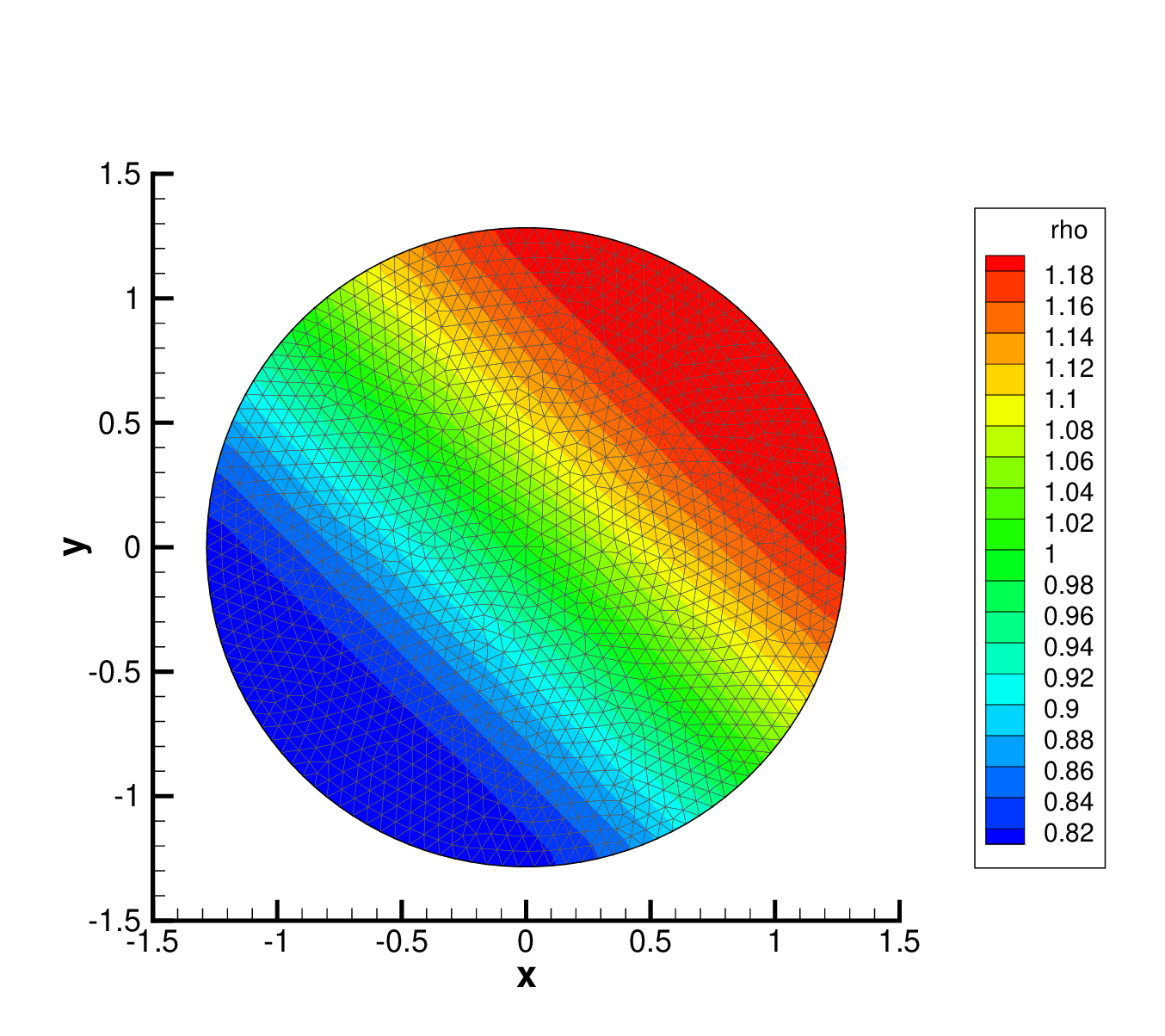}}
\caption{Manufactured solution on 2D moving meshes: initial (left) and final (right) mesh configuration.
The mesh is moving radially outward with a velocity field that depends on the distance from the center of the domain.
The color map represents the density field.}
\label{fig:manuf2D}
\end{figure}

In order to test the accuracy of the new boundary conditions on 2D moving meshes, we set up a smooth manufactured solution,
in a circular domain, by considering the 2D inhomogeneous Euler equations with source term
\begin{equation*}\label{eq:manufactured}
\mathbf{S} = \left(
    \begin{array}{c}
    0.4\cos(x+y) \\ 0.6\cos(x+y)  \\ 0.6\cos(x+y) \\ 1.8\cos(x+y) 
      \end{array}
    \right).
\end{equation*}
The exact steady state solution reads as
$$ \rho = 1 + 0.2 \sin(x+y),\quad u=1,\quad v=1, \quad p=1 + 0.2 \sin(x+y),$$
which is imposed on the boundary as a ghost state.
To study the impact of the boundary correction on the moving mesh, we impose the following radial velocity field
on each boundary node $\nu$:
$$ \mathbf{V}_{\nu} = u_0 \|\mathbf{x}_{\nu}\| \left(
  \begin{array}{c}
   \cos(\alpha_{\nu}) \\ \sin(\alpha_{\nu}) 
    \end{array}
  \right),$$
where $u_0=0.1$ is the maximum velocity, and $\alpha_{\nu}=\arctan(y_{\nu},x_{\nu})$.
It can be noticed that since the domain is evolving the solution, which is space-dependent, the boundary condition
is also changing in time, providing a more challenging test case for the high order boundary conditions.
We show in figure \ref{fig:manuf2D} the initial and final mesh configurations, 
where the color map represents the density field.
The simulations are performed on a set of four meshes with and without the polynomial correction until the final time $t_f=0.5$, 
and the convergence analysis is presented in \cref{tab:manuf2D}.

Due to the curvature of the external boundary, approximated through a polygonal representation, 
we expect to achieve second order accuracy at best for the ALE-FV with trivial imposition of boundary conditions.
Indeed, the results in \cref{tab:manuf2D} confirm this expectation, with the convergence rates of primitive variables 
being around 2, no matter the order of the polynomials used for the reconstruction of the numerical solution.
On the other hand, the results obtained with the new correction demonstrate that the novel method allows us to recover 
the expected high order accuracy of the scheme. In particular, the convergence rates are close to the expected 
order of accuracy for all variables and polynomial orders, confirming the effectiveness of the proposed methodology on moving triangular grids. 

\begin{table}
  \caption{Manufactured solution on 2D moving meshes: convergence analysis for the test case presented in~\cref{subsection:Manufactured2D}. 
  Numerical results obtained with Dirichlet-type boundary conditions imposed with and w/o the shifted polynomial correction 
  on conformal linear meshes.}\label{tab:manuf2D}
  \footnotesize
  \centering
  \begin{tabular}{ccccccccc}
          \hline\hline          &\multicolumn{2}{c}{$\rho$} &\multicolumn{2}{c}{$u$}   &\multicolumn{2}{c}{$\rho$} &\multicolumn{2}{c}{$u$}  \\[0.5mm]
          \cline{2-9}
          Grid size & $L_2$        & $\tilde{n}$ & $L_2$        & $\tilde{n}$ & $L_2$      & $\tilde{n}$  & $L_2$      & $\tilde{n}$ \\[0.5mm]\hline
          &\multicolumn{8}{c}{FV-$\mathcal{P}_1$}\\
          &\multicolumn{4}{c}{w/o correction}   &\multicolumn{4}{c}{with correction} \\[0.5mm]
          1.90E-01  &  4.47E-03  &   --   &  2.00E-03   &  --  &  4.65E-03 &  --   & 1.95E-03 &  --   \\
          9.76E-02  &  1.27E-03  &  1.89  &  7.32E-04   & 1.51 &  1.29E-03 &  1.92 & 7.30E-04 & 1.47   \\
          5.04E-02  &  3.59E-04  &  1.91  &  2.27E-04   & 1.77 &  3.62E-04 &  1.92 & 2.27E-04 & 1.77   \\
          2.55E-02  &  9.90E-05  &  1.89  &  6.78E-05   & 1.77 &  9.93E-05 &  1.89 & 6.79E-05 & 1.77   \\
          &\multicolumn{8}{c}{FV-$\mathcal{P}_2$}\\
          &\multicolumn{4}{c}{w/o correction}   &\multicolumn{4}{c}{with correction} \\[0.5mm]
          1.90E-01  &  6.83E-04  &   --   &  2.30E-04   & --   &  7.82E-04 &  --   & 2.53E-04 & --   \\
          9.76E-02  &  8.56E-05  &  3.11  &  4.32E-05   & 2.51 &  8.17E-05 &  3.39 & 2.49E-05 & 3.47  \\
          5.04E-02  &  1.73E-05  &  2.42  &  1.11E-05   & 2.05 &  1.02E-05 &  3.15 & 3.60E-06 & 2.92  \\
          2.55E-02  &  4.33E-06  &  2.03  &  2.84E-06   & 2.00 &  1.24E-06 &  3.08 & 4.64E-07 & 3.00  \\
          &\multicolumn{8}{c}{FV-$\mathcal{P}_3$}\\
          &\multicolumn{4}{c}{w/o correction}   &\multicolumn{4}{c}{with correction} \\[0.5mm]
          1.90E-01  & 2.66E-04  &   --   &  1.96E-04   & --    &  1.02E-04 &  --   & 5.30E-05 & --    \\
          9.76E-02  & 6.78E-05  &  2.05  &  4.56E-05   & 2.19  &  6.54E-06 &  4.12 & 4.06E-06 & 3.85  \\
          5.04E-02  & 1.80E-05  &  2.00  &  1.18E-05   & 2.05  &  5.29E-07 &  3.80 & 3.92E-07 & 3.53  \\
          2.55E-02  & 4.56E-06  &  2.01  &  2.97E-06   & 2.02  &  3.72E-08 &  3.89 & 3.32E-08 & 3.62  \\
          \hline\hline\\[1pt]
  \end{tabular}
  \end{table}

\subsection{Manufactured solution on 3D moving meshes}\label{subsection:Manufactured3D}

The goal of this section is again to assess the performance of the high order polynomial correction, but this time 
on 3D moving meshes. We consider a smooth manufactured solution in a sphere domain, by considering the following
source term in \eqref{eq:euler0}:
\begin{equation*}\label{eq:manufactured3D}
  \mathbf{S} = \left( \begin{array}{c}  0.6\cos(x+y+z) \\ 0.8\cos(x+y+z)  \\ 0.8\cos(x+y+z) \\ 0.8\cos(x+y+z) \\ 3\cos(x+y+z) \end{array} \right),
\end{equation*}
The exact steady state solution is given by
$$ \rho = 1 + 0.2 \sin(x+y+z),\quad u=1,\quad v=1, \quad w=1, \quad p=1 + 0.2 \sin(x+y+z),$$
which is imposed on the boundary as a ghost state.

\begin{figure}
  \centering
  \subfigure[]{\includegraphics[width=0.48\textwidth]{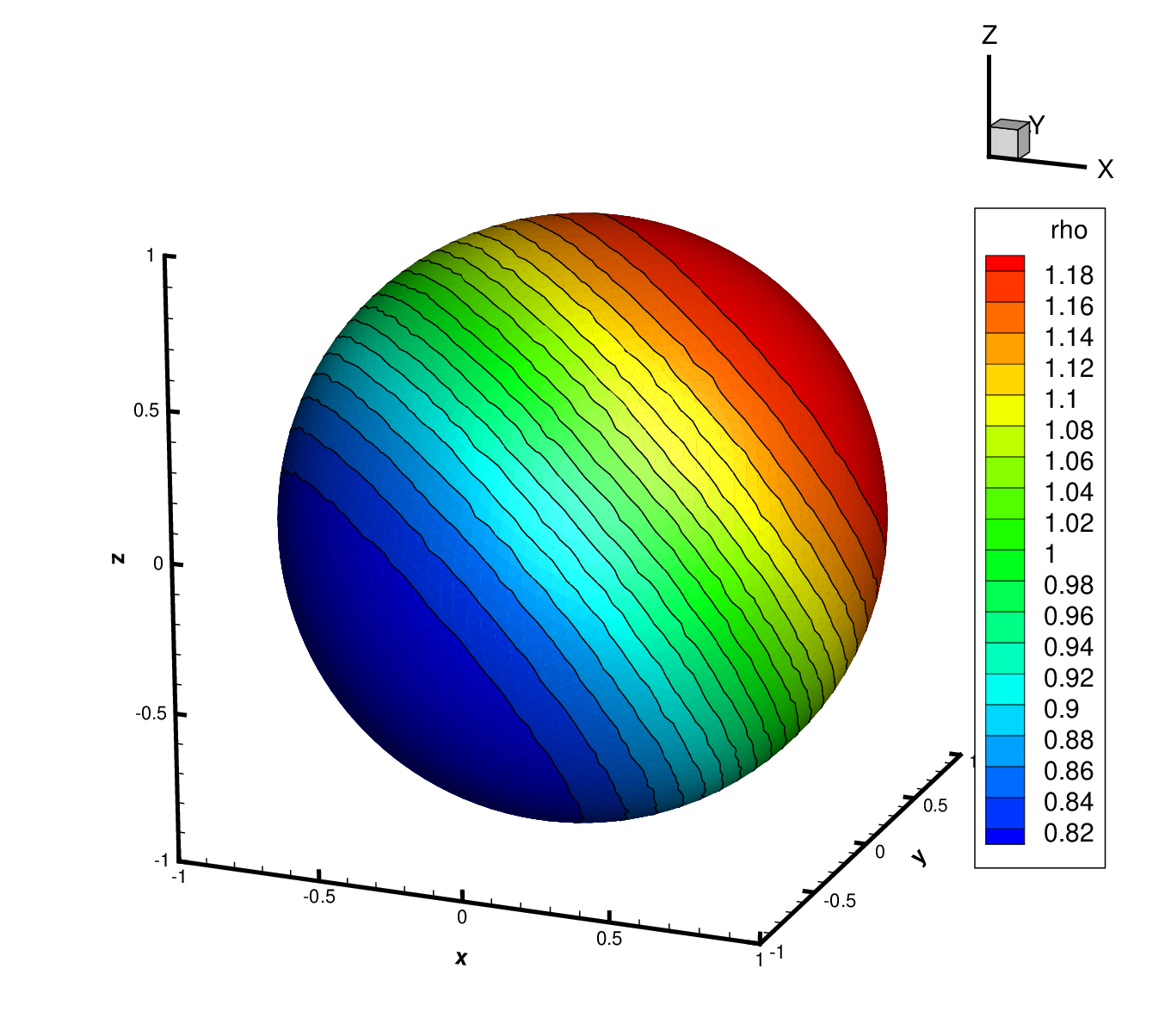}}
  \subfigure[]{\includegraphics[width=0.48\textwidth]{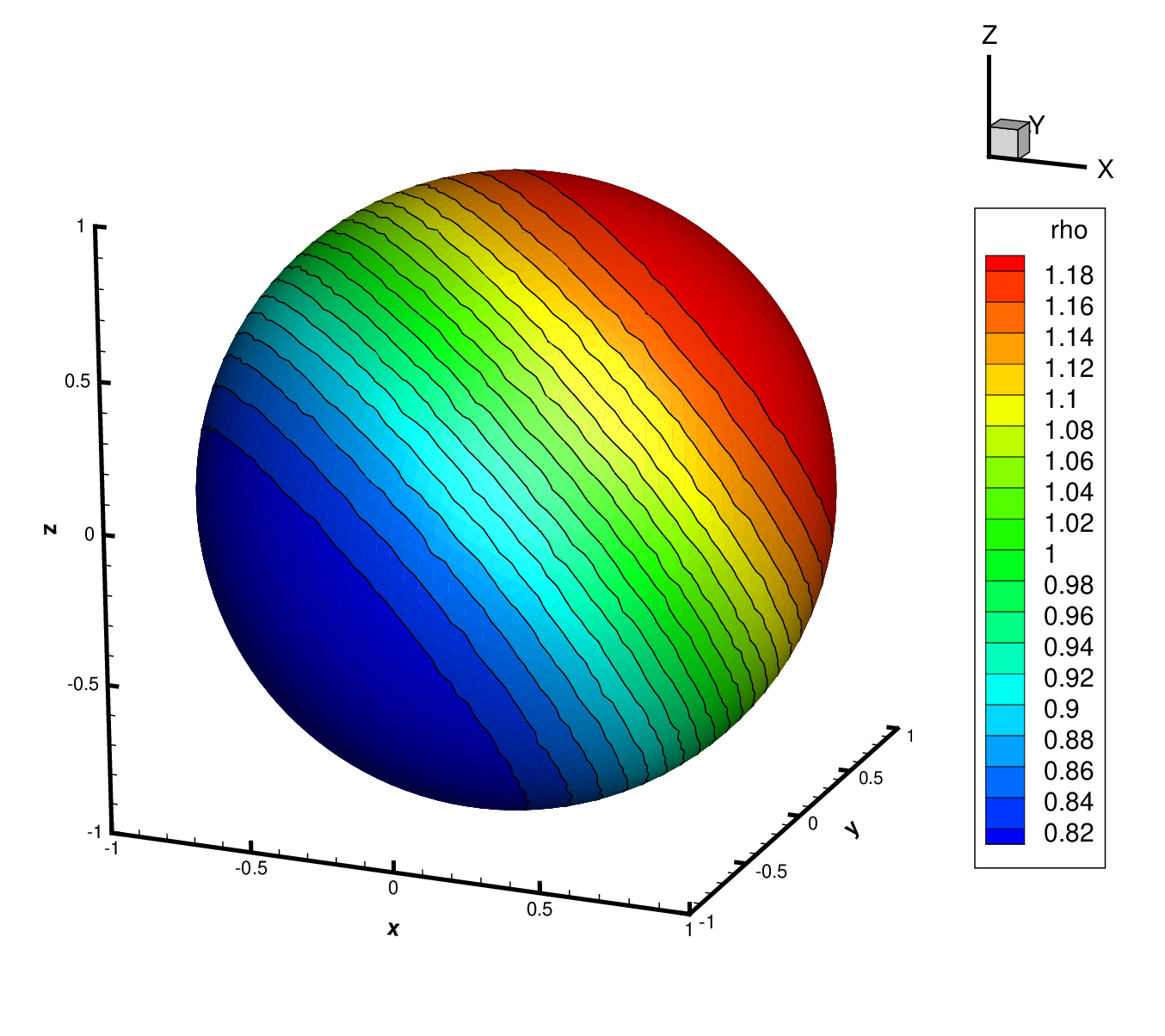}}
  \caption{Manufactured solution on 3D moving meshes: initial (left) and final (right) mesh configuration.
  The mesh is moving radially outward.
  The color map represents the density field.}
  \label{fig:manuf3D}
\end{figure}

To study the impact of the boundary correction on the moving mesh, we impose the following velocity field
on each boundary node $\nu$:
$$ \mathbf{V}_{\nu} = u_0 \frac{\mathbf{x}}{\|\mathbf{x}\|} ,$$
where $u_0=0.1$ is the maximum velocity.
Once again, the solution is space-dependent, and the domain is evolving, providing a moving test case for the proposed high order boundary conditions
in the considered ALE framework. In figure \ref{fig:manuf3D} we show the three-dimensional initial and final mesh configurations.
The simulations are always performed on a set of four meshes with and without the polynomial correction,
and the convergence analysis is presented in \cref{tab:manuf3D}.
As expected, the results in \cref{tab:manuf3D} show that the trivial imposition of boundary conditions leads to a loss of accuracy,
with the convergence rates being around 2 for all variables and polynomial orders.
On the other hand, the results obtained through the polynomial correction prove that the new method allows us to recover the expected high order accuracy of the scheme on moving tetrahedral meshes.
In particular, the convergence rates are close to the expected order of accuracy.

\begin{table}
  \caption{Manufactured solution on 3D moving meshes: convergence analysis for the test case presented in~\cref{subsection:Manufactured3D}. 
  Numerical results obtained with Dirichlet-type boundary conditions imposed with and w/o the shifted boundary polynomial correction 
  on conformal linear meshes.}\label{tab:manuf3D}
  \footnotesize
  \centering
  \begin{tabular}{ccccccccc}
          \hline\hline          &\multicolumn{2}{c}{$\rho$} &\multicolumn{2}{c}{$u$}   &\multicolumn{2}{c}{$\rho$} &\multicolumn{2}{c}{$u$}  \\[0.5mm]
          \cline{2-9}
          Grid size & $L_2$        & $\tilde{n}$ & $L_2$        & $\tilde{n}$ & $L_2$      & $\tilde{n}$  & $L_2$      & $\tilde{n}$ \\[0.5mm]\hline
          &\multicolumn{8}{c}{FV-$\mathcal{P}_1$}\\
          &\multicolumn{4}{c}{w/o correction}   &\multicolumn{4}{c}{with correction} \\[0.5mm]
          2.56E-01  &  3.44E-03  &   --  &  1.13E-03   &  --  & 4.27E-03  &  --   &  8.75E-04   &  --    \\
          1.45E-01  &  1.27E-03  & 1.75  &  4.69E-04   & 1.55 & 1.40E-03  & 1.96  &  4.19E-04   & 1.29    \\
          7.75E-02  &  3.52E-04  & 2.06  &  1.56E-04   & 1.76 & 3.76E-04  & 2.10  &  1.48E-04   & 1.67    \\
          3.91E-02  &  8.99E-05  & 2.00  &  4.54E-05   & 1.80 & 9.48E-05  & 2.02  &  4.42E-05   & 1.76    \\
          &\multicolumn{8}{c}{FV-$\mathcal{P}_2$}\\
          &\multicolumn{4}{c}{w/o correction}   &\multicolumn{4}{c}{with correction} \\[0.5mm]
          2.56E-01  &  2.02E-03  &  --   &  4.40E-04   &   --  & 2.13E-03  &   --  &  2.37E-04   &  --   \\
          1.45E-01  &  4.09E-04  & 2.80  &  1.21E-04   &  2.27 & 3.39E-04  &  3.22 &  4.87E-05   & 2.78  \\    
          7.75E-02  &  8.38E-05  & 2.54  &  3.24E-05   &  2.11 & 4.42E-05  &  3.26 &  8.14E-06   & 2.86  \\  
          3.91E-02  &  1.94E-05  & 2.14  &  8.28E-06   &  2.00 & 5.07E-06  &  3.17 &  1.15E-06   & 2.87  \\       
          &\multicolumn{8}{c}{FV-$\mathcal{P}_3$}\\
          &\multicolumn{4}{c}{w/o correction}   &\multicolumn{4}{c}{with correction} \\[0.5mm]
          2.56E-01  &  1.41E-03  &   --  &  4.12E-04   &  --   &  2.91E-04  &   --  &  4.63E-05   &   --  \\
          1.45E-01  &  3.30E-04  & 2.55  &  1.28E-04   & 2.06  &  2.69E-05  &  4.18 &  6.64E-06   &  3.41 \\
          7.75E-02  &  7.95E-05  & 2.28  &  3.38E-05   & 2.13  &  1.96E-06  &  4.19 &  6.17E-07   &  3.81 \\
          3.91E-02  &  1.95E-05  & 2.06  &  8.53E-06   & 2.02  &  1.16E-07  &  4.14 &  4.70E-08   &  3.77 \\
          \hline\hline\\[1pt]
  \end{tabular}
  \end{table}

\subsection{Kidder problem in 2D}\label{subsection:Kidder2D}

The so-called Kidder problem is an isentropic compression of a shell filled with an ideal gas. 
This problem, along with its exact analytical solution was proposed by Kidder in \cite{kidder1976laser}.
The problem is often used as a benchmark for Lagrangian hydrodynamics codes to assure that there is no production
of spurious entropy. The shell has a time–dependent internal radius $r_i(t)$ and an external radius $r_e(t)$.
The initial values for the internal and external radius are $r_i(0)=r_{i,0}=0.9$ and $r_e(0)=r_{e,0}=1$, respectively.
The ratio of specific heats is $\gamma=2$, and the initial density distribution is given by
$$ \rho_0(r) = \rho(r,0) = \left( \frac{r^2_{e,0}-r^2}{r^2_{e,0}-r^2_{i,0}}\rho_{i,0}^{\gamma-1} + 
\frac{r^2-r^2_{i,0}}{r^2_{e,0}-r^2_{i,0}}\rho_{e,0}^{\gamma-1} \right)^{\frac{1}{\gamma-1}},\qquad\text{with}\quad r_i(t)\leq r \leq r_e(t), $$
where we set $\rho_{i,0}=1$ and $\rho_{e,0}=2$,  which are the initial density values defined at the internal and 
external boundary of the shell, respectively.
The initial entropy $s_0=\frac{p_0}{\rho_0^\gamma}=1$ is uniform such that the initial pressure distribution
is given by $p_0(r)= s_0 \rho_0(r)^{\gamma}$.
Initially the fluid is set at rest ($u=v\equiv 0$).

The time-varying solution of the Kidder problem can be written as a solution of the form $R(r,t)=h(t) r$, where
$R(r,t)$ denotes the radius at time $t>0$ of a fluid particle initially located at $r$. Hence, the self-similar solution
for $t\in[0,\tau]$ is given by
\begin{align*}
\rho(R(r,t),t) &=   h(t)^{-\frac{2}{\gamma-1}} \rho_0\left(\frac{R(r,t)}{h(t)}\right),\\
u_r(R(r,t),t)  &= \frac{R(r,t)}{h(t)} \frac{\diff{}}{\diff{t}}h(t), \\
p(R(r,t),t)    &=  h(t)^{-\frac{2\gamma}{\gamma-1}}p_0\left(\frac{R(r,t)}{h(t)}\right),
\end{align*}
where $u_r$ is the radial velocity component.\\
The homothety rate reads
$$ h(t) = \sqrt{1-\frac{t^2}{\tau^2}}, $$
with $\tau$ representing the focalisation time
$$ \tau = \sqrt{\frac{\gamma-1}{2} \frac{r^2_{e,0}-r^2_{i,0}}{c^2_{e,0}-c^2_{i,0}}}, $$
while the internal and external sound speeds are given by
$$ c_{i,0} = \sqrt{\gamma\frac{ p_{i,0}}{\rho_{i,0}}}\quad \text{ and } \quad c_{e,0} = \sqrt{\gamma\frac{p_{e,0}}{\rho_{e,0}}}. $$

The boundary conditions are imposed on the internal and external boundaries of the shell, and the exact solution is imposed as a ghost state. The exact mesh velocity is prescribed at the boundary nodes.

\begin{figure}
  \centering
  \subfigure[]{\includegraphics[width=0.48\textwidth]{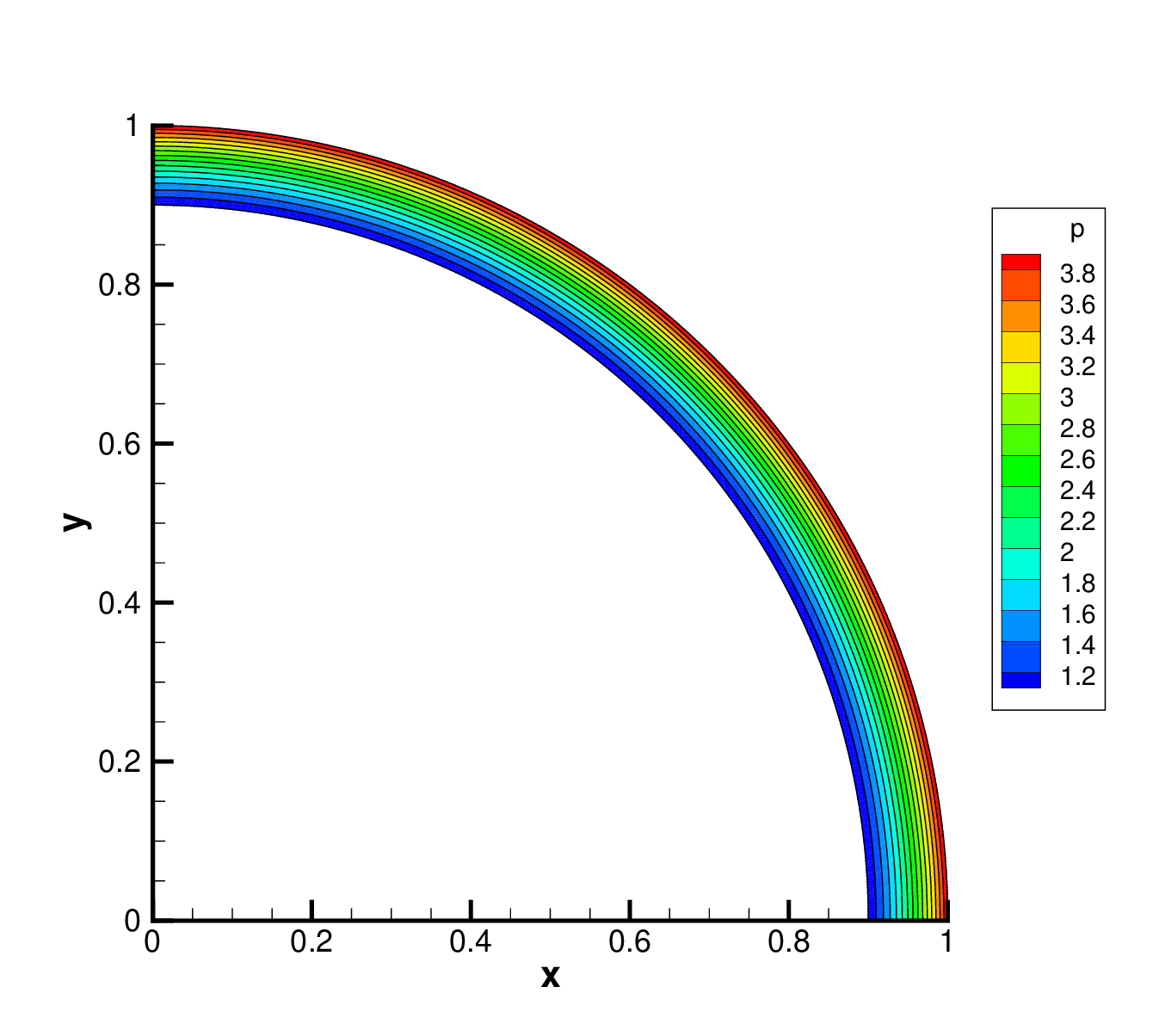}}
  \subfigure[]{\includegraphics[width=0.48\textwidth]{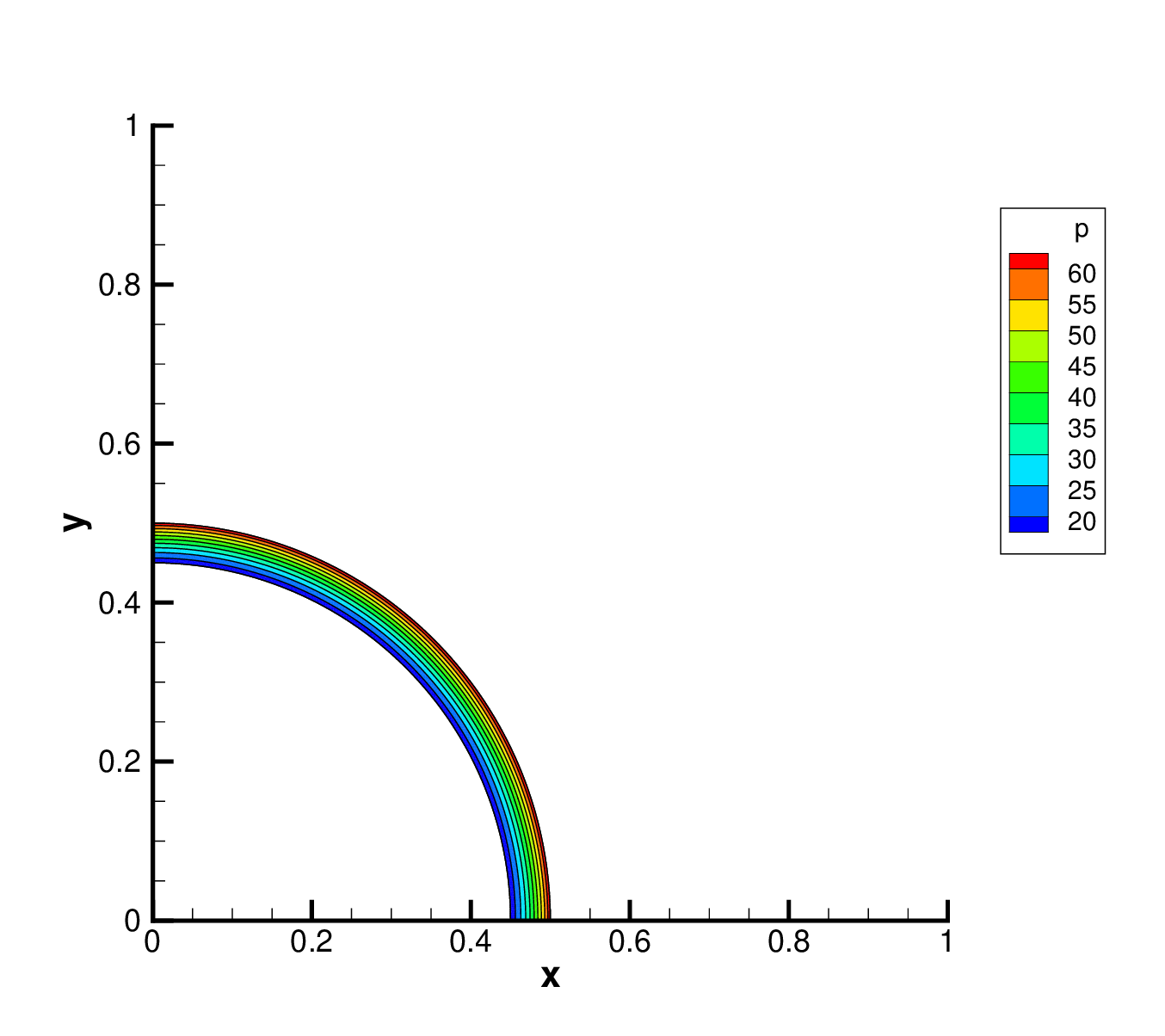}}
  \caption{Kidder problem in 2D: initial (left) and final (right) mesh configuration.
  The mesh is moving radially with a velocity field that depends on the Kidder exact solution.
  The color map represents the pressure field.}
  \label{fig:Kidder2D}
\end{figure}

Following \cite{carre2009cell}, the final time of the simulation is chosen as $t_f=\frac{\sqrt{3}}{2}\tau$, 
such that the resulting compression rate at final time is $h(t_f)=1/2$, and the exact location of the final shell
is bounded between 0.45 and 0.5.
In figure \ref{fig:Kidder2D} we show the initial and final mesh configurations.
Although the simulations are performed on an entire circle domain, we only show its top right quarter 
for visualization purposes. 

In \cref{tab:Kidder2D}, we present the convergence analysis for the Kidder problem in 2D.
The results show that the high order boundary conditions allow us to recover the expected high order accuracy of the scheme,
while the trivial imposition of boundary conditions leads to a loss of accuracy.
In particular, it should be noticed that without any boundary correction, the convergence rates are around 2 for all variables and polynomial orders.
This leads to numerical solutions that are almost as accurate as the low order scheme, which makes the use of high order approximation
not worth the effort. On the other hand, by using the polynomial correction we achieve a huge improvement in terms of accuracy
with discretization errors on the finest mesh being several orders of magnitude smaller than the ones obtained without the correction. 

\begin{table}
  \caption{Kidder problem in 2D: convergence analysis for the test case presented in~\cref{subsection:Kidder2D}. 
  Numerical results obtained with Dirichlet-type boundary conditions imposed with and w/o the shifted polynomial correction 
  on conformal linear meshes.}\label{tab:Kidder2D}
  \footnotesize
  \centering
  \begin{tabular}{ccccccccc}
          \hline\hline          &\multicolumn{2}{c}{$\rho$} &\multicolumn{2}{c}{$u$}   &\multicolumn{2}{c}{$\rho$} &\multicolumn{2}{c}{$u$}  \\[0.5mm]
          \cline{2-9}
          Grid size & $L_2$        & $\tilde{n}$ & $L_2$        & $\tilde{n}$ & $L_2$      & $\tilde{n}$  & $L_2$      & $\tilde{n}$ \\[0.5mm]\hline
          &\multicolumn{8}{c}{FV-$\mathcal{P}_1$}\\
          &\multicolumn{4}{c}{w/o correction}   &\multicolumn{4}{c}{with correction} \\[0.5mm]
          1.28E-02  &  3.56E-02 & --   & 5.54E-03 & --   & 3.54E-02 & --   & 5.54E-03 &  --  \\
          5.99E-03  &  1.28E-02 & 1.35 & 1.57E-03 & 1.66 & 1.27E-02 & 1.34 & 1.57E-03 & 1.66  \\
          3.06E-03  &  2.97E-03 & 2.17 & 3.48E-04 & 2.24 & 2.98E-03 & 2.16 & 3.48E-04 & 2.24  \\
          1.49E-03  &  7.35E-04 & 1.95 & 7.93E-05 & 2.06 & 7.38E-04 & 1.95 & 7.95E-05 & 2.06  \\
          &\multicolumn{8}{c}{FV-$\mathcal{P}_2$}\\
          &\multicolumn{4}{c}{w/o correction}   &\multicolumn{4}{c}{with correction} \\[0.5mm]
          1.28E-02  &  3.46E-03 & --   & 3.33E-04 & --   & 2.73E-03 &  --  & 3.20E-04 &  --  \\
          5.99E-03  &  6.35E-04 & 2.23 & 5.48E-05 & 2.38 & 3.84E-04 & 2.59 & 4.07E-05 & 2.72  \\
          3.06E-03  &  1.20E-04 & 2.48 & 1.09E-05 & 2.40 & 4.05E-05 & 3.34 & 4.08E-06 & 3.42  \\
          1.49E-03  &  2.54E-05 & 2.16 & 2.62E-06 & 1.99 & 4.34E-06 & 3.11 & 3.93E-07 & 3.26  \\
          &\multicolumn{8}{c}{FV-$\mathcal{P}_3$}\\
          &\multicolumn{4}{c}{w/o correction}   &\multicolumn{4}{c}{with correction} \\[0.5mm]
          1.28E-02  &  1.39E-03 & --   & 1.68E-04 & --   & 1.74E-04 &  --  & 4.03E-05 &  -- \\
          5.99E-03  &  3.47E-04 & 1.82 & 4.22E-05 & 1.82 & 9.65E-06 & 3.81 & 2.36E-06 & 3.74  \\
          3.06E-03  &  8.73E-05 & 2.05 & 1.06E-05 & 2.06 & 4.26E-07 & 4.64 & 1.23E-07 & 4.39  \\
          1.49E-03  &  2.17E-05 & 1.94 & 2.63E-06 & 1.94 & 1.92E-08 & 4.32 & 5.50E-09 & 4.34  \\
          \hline\hline\\[1pt]
  \end{tabular}
  \end{table}

  \subsection{Kidder problem in 3D}\label{subsection:Kidder3D}

\begin{figure}
  \centering
  \subfigure[]{\includegraphics[width=0.48\textwidth]{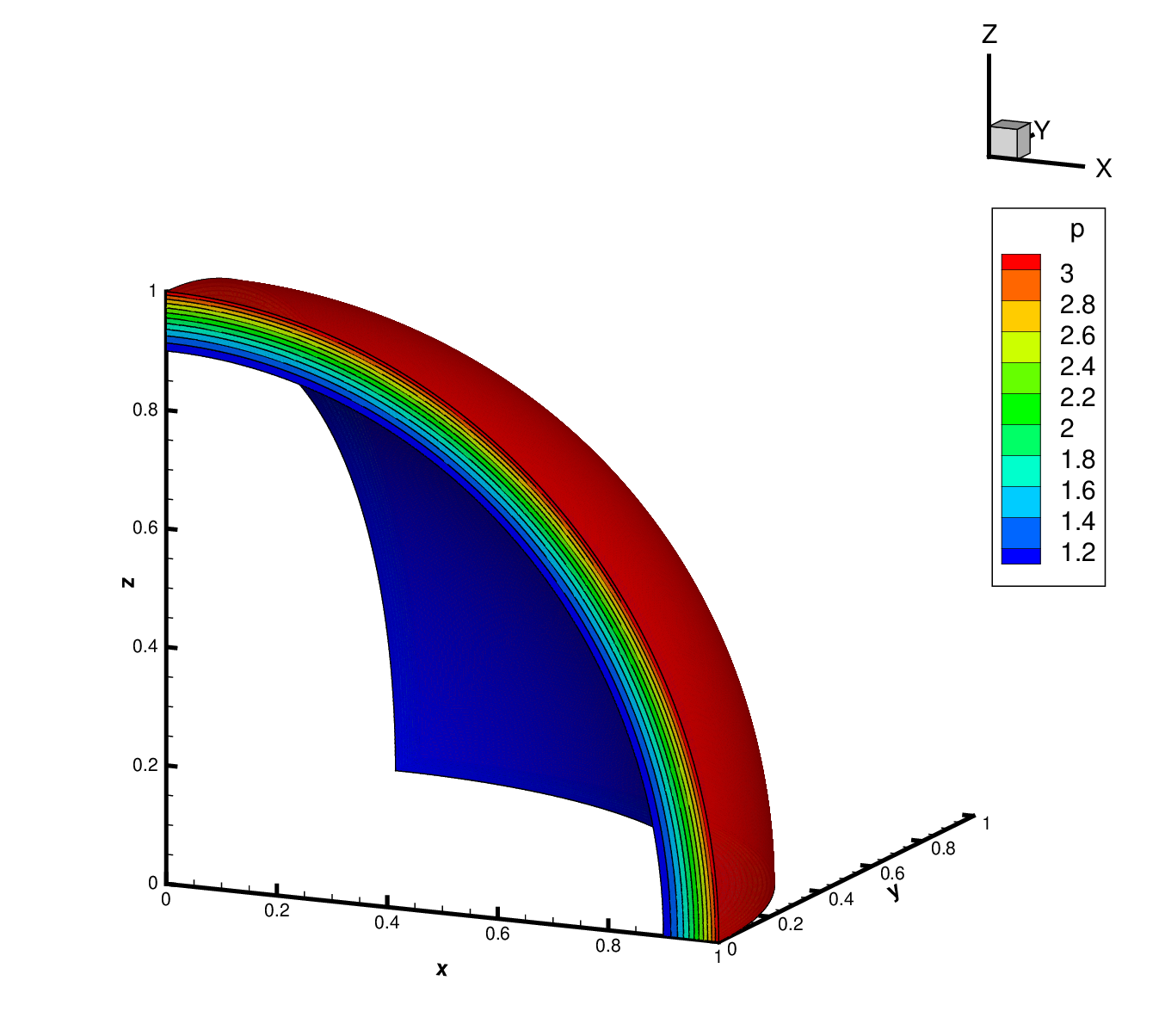}}
  \subfigure[]{\includegraphics[width=0.48\textwidth]{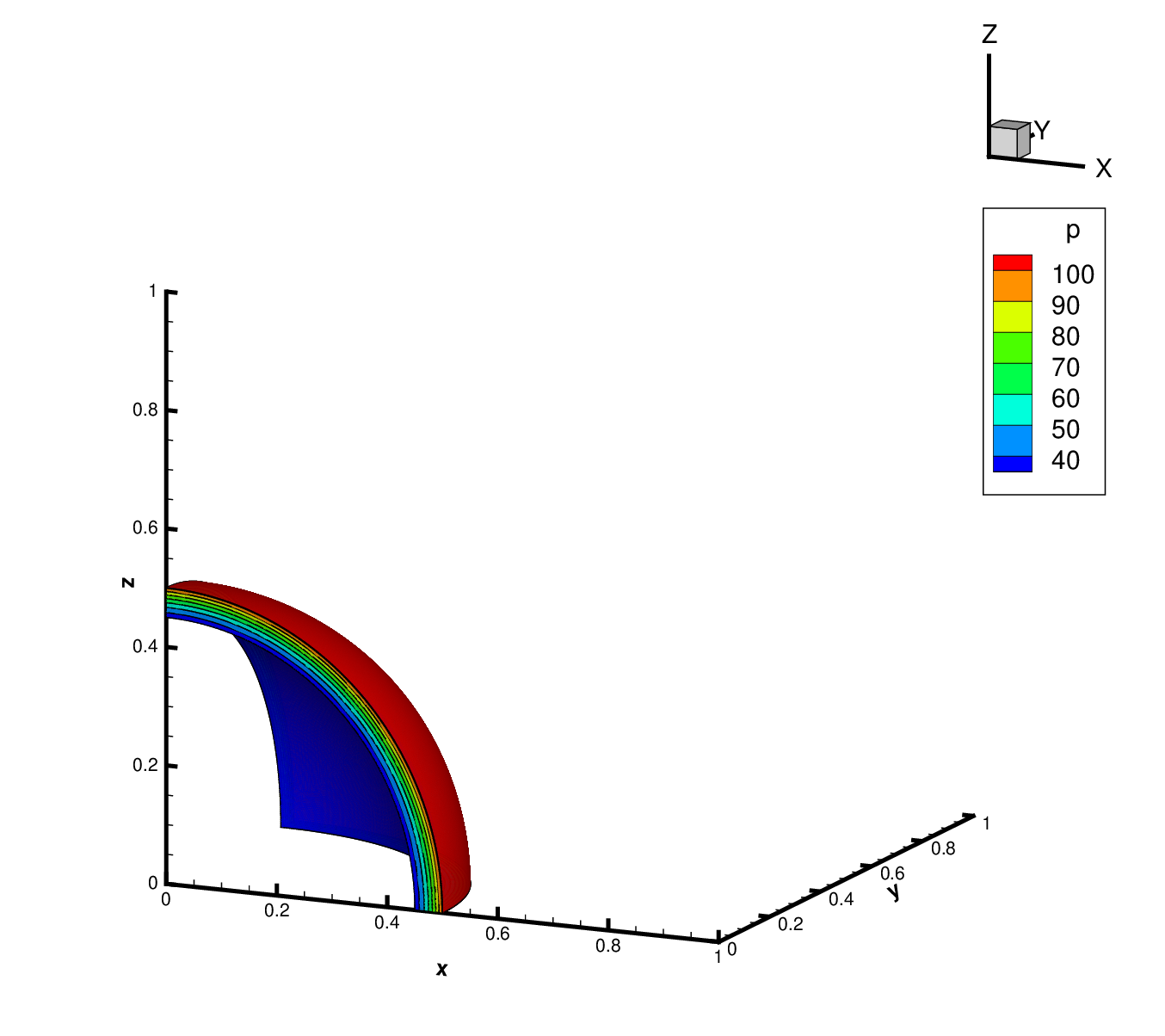}}
  \caption{Kidder problem in 3D: initial (left) and final (right) mesh configuration.
  The mesh is moving radially with a velocity field that depends on the Kidder exact solution.
  The color map represents the pressure field.}
  \label{fig:Kidder3D}
\end{figure}

In this section, we aim at testing the Kidder problem on three-dimensional unstructured moving meshes.
Due to the much higher computational cost required to run 3D simulations, we consider as initial computational domain
one eighth of the entire shell. On lateral faces, classical boundary conditions are imposed since they are simple 
straight faces with no curvature. In figure \ref{fig:Kidder3D} we show the mesh configurations at the initial and final time 
for the part of the shell considered for the simulations. The setup of the 3D version of the Kidder problem is slightly different compared to the two-dimensional setting, see \cite{boscheri2014direct} for further details.

In \cref{tab:Kidder3D}, we report the convergence analysis for the Kidder problem in 3D, where the discretization
errors are presented.
Again, the results demonstrate that the high order boundary conditions allow us to recover the expected high order accuracy of the scheme even 
in this 3D case, with discretization errors that are much lower than those obtained without the correction.
In particular, on the finest mesh for $\mathcal{P}_3$, the discretization errors are almost four orders of magnitude smaller than the ones 
obtained without the correction.
For the simulations performed without the correction we obtain convergence rates that are extremely degraded.
Since the only difference between the simulations is the boundary treatment, we can conclude that
this order-of-accuracy degradation might be due to the strong anisotropy of the 3D mesh that, 
combined with the low order boundary treatment, provides such a spoiled convergence. 

\begin{figure}
  \centering
  \subfigure[]{\includegraphics[width=0.48\textwidth]{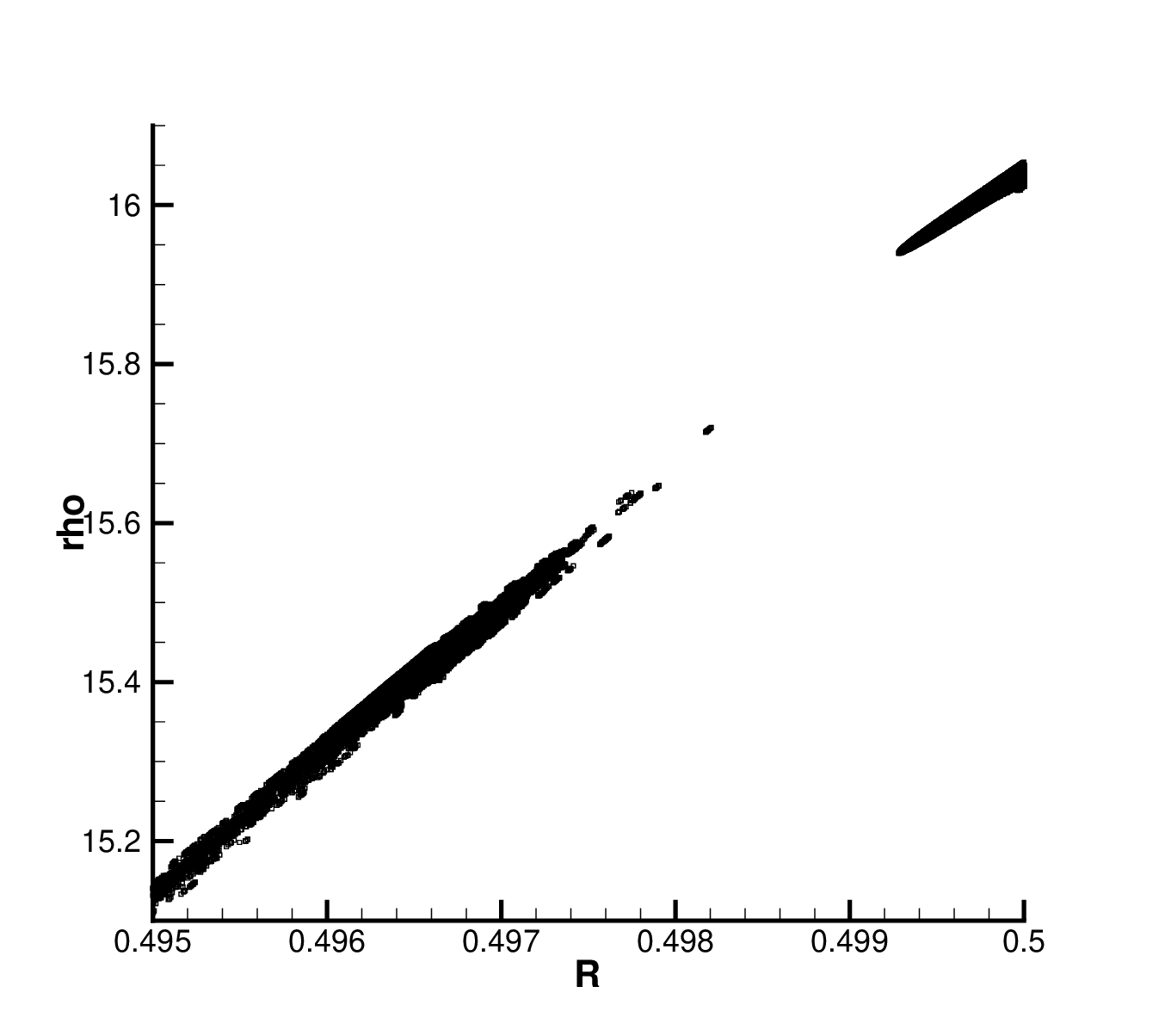}}
  \subfigure[]{\includegraphics[width=0.48\textwidth]{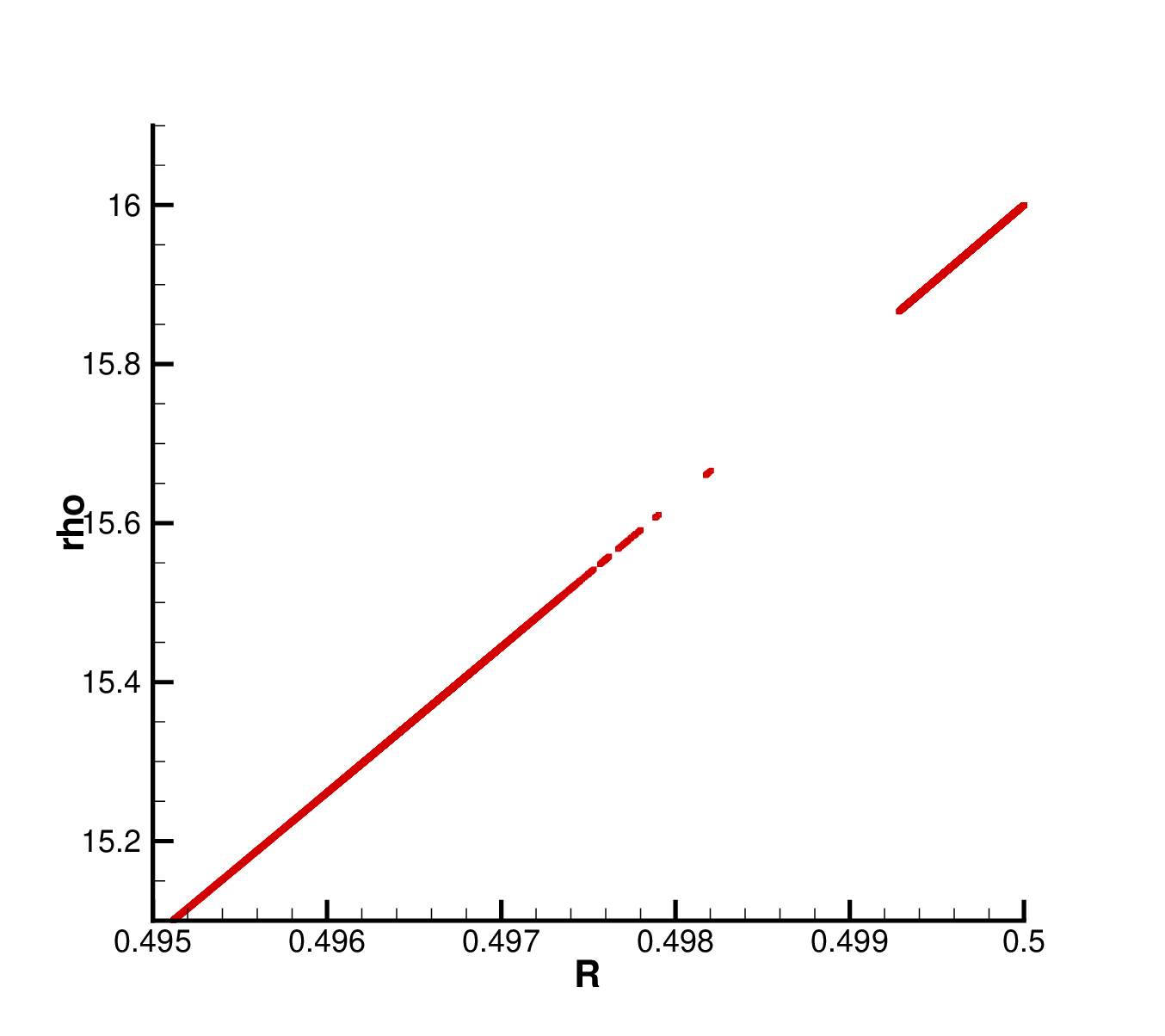}}
  \caption{Kidder problem in 3D: scatter plot of the density with respect to the 
  radius without correction (left) and with correction (right).}
  \label{fig:Kidder3Dscatter}
\end{figure}

In figure \ref{fig:Kidder3Dscatter} we also show the scatter plot of the density with respect to the radius for the two cases,
with and without the polynomial correction.
It is extremely clear from the scatter plot that the polynomial correction allows us to recover the expected density distribution, 
which should be linear in the radial direction. On the contrary, the simulation performed without the correction
shows a very irregular distribution of the density. Let us also remark that a more radially symmetric solution is attained even in the interior part of the domain when the correct boundary treatment is carried out. 
We would like to stress again the fact that since the mesh motion is prescribed by the exact solution of the problem 
for both simulations, the only difference between the two cases is the treatment of the boundary conditions.

\begin{table}
    \caption{Kidder problem in 3D: convergence analysis for the test case presented in~\cref{subsection:Kidder3D}. 
    Numerical results obtained with Dirichlet-type boundary conditions imposed with and w/o the shifted polynomial correction 
    on conformal linear meshes.}\label{tab:Kidder3D}
    \footnotesize
    \centering
    \begin{tabular}{ccccccccc}
            \hline\hline          &\multicolumn{2}{c}{$\rho$} &\multicolumn{2}{c}{$u$}   &\multicolumn{2}{c}{$\rho$} &\multicolumn{2}{c}{$u$}  \\[0.5mm]
            \cline{2-9}
            Grid size & $L_2$        & $\tilde{n}$ & $L_2$        & $\tilde{n}$ & $L_2$      & $\tilde{n}$  & $L_2$      & $\tilde{n}$ \\[0.5mm]\hline
            &\multicolumn{8}{c}{FV-$\mathcal{P}_1$}\\
            &\multicolumn{4}{c}{w/o correction}   &\multicolumn{4}{c}{with correction} \\[0.5mm]
            0.05 & 5.14E-02  &  --  & 4.93E-03 &  --  & 4.86E-02 &  --  & 5.73E-03 &  --   \\
            0.04 & 4.46E-02  & 1.91 & 5.89E-03 &  --  & 4.24E-02 & 1.83 & 4.69E-03 & 2.71  \\
            0.03 & 2.15E-02  & 2.28 & 3.09E-03 & 2.02 & 1.88E-02 & 2.56 & 2.75E-03 & 1.68  \\
            0.02 & 1.06E-02  & 1.46 & 1.78E-03 & 1.15 & 7.07E-03 & 2.02 & 1.05E-03 & 1.99  \\
            &\multicolumn{8}{c}{FV-$\mathcal{P}_2$}\\
            &\multicolumn{4}{c}{w/o correction}   &\multicolumn{4}{c}{with correction} \\[0.5mm]
            0.05 & 1.50E-03  &  --  & 1.64E-03 &  --  & 3.78E-03 &  --  & 3.00E-04 &  --   \\
            0.04 & 1.03E-02  &  --  & 1.56E-03 & 0.67 & 2.93E-03 & 3.43 & 2.45E-04 & 2.70  \\
            0.03 & 8.27E-03  & 0.68 & 1.50E-03 & 0.12 & 1.17E-03 & 2.87 & 1.09E-04 & 2.55  \\
            0.02 & 6.97E-03  & 0.35 & 1.46E-03 & 0.05 & 3.05E-04 & 2.78 & 2.55E-05 & 3.00  \\
            &\multicolumn{8}{c}{FV-$\mathcal{P}_3$}\\
            &\multicolumn{4}{c}{w/o correction}   &\multicolumn{4}{c}{with correction} \\[0.5mm]
            0.05 & 9.78E-03  &  --  & 1.64E-03 &  --  & 3.62E-04 &  --  & 3.26E-05 &  -- \\
            0.04 & 8.47E-03  & 1.95 & 1.55E-03 & 0.70 & 2.01E-04 & 7.94 & 2.47E-05 & 3.78 \\
            0.03 & 7.49E-03  & 0.39 & 1.50E-03 & 0.11 & 4.18E-05 & 4.93 & 5.95E-06 & 4.46 \\
            0.02 & 6.75E-03  & 0.21 & 1.46E-03 & 0.06 & 6.57E-06 & 3.82 & 9.95E-07 & 3.69 \\
            \hline\hline\\[1pt]
    \end{tabular}
\end{table}

\subsection{Fluid-structure interaction: horizontal and vertical oscillating cylinders}\label{subsection:FSIcylinder}

\begin{figure}
  \centering
  \subfigure[]{\includegraphics[width=0.4\textwidth]{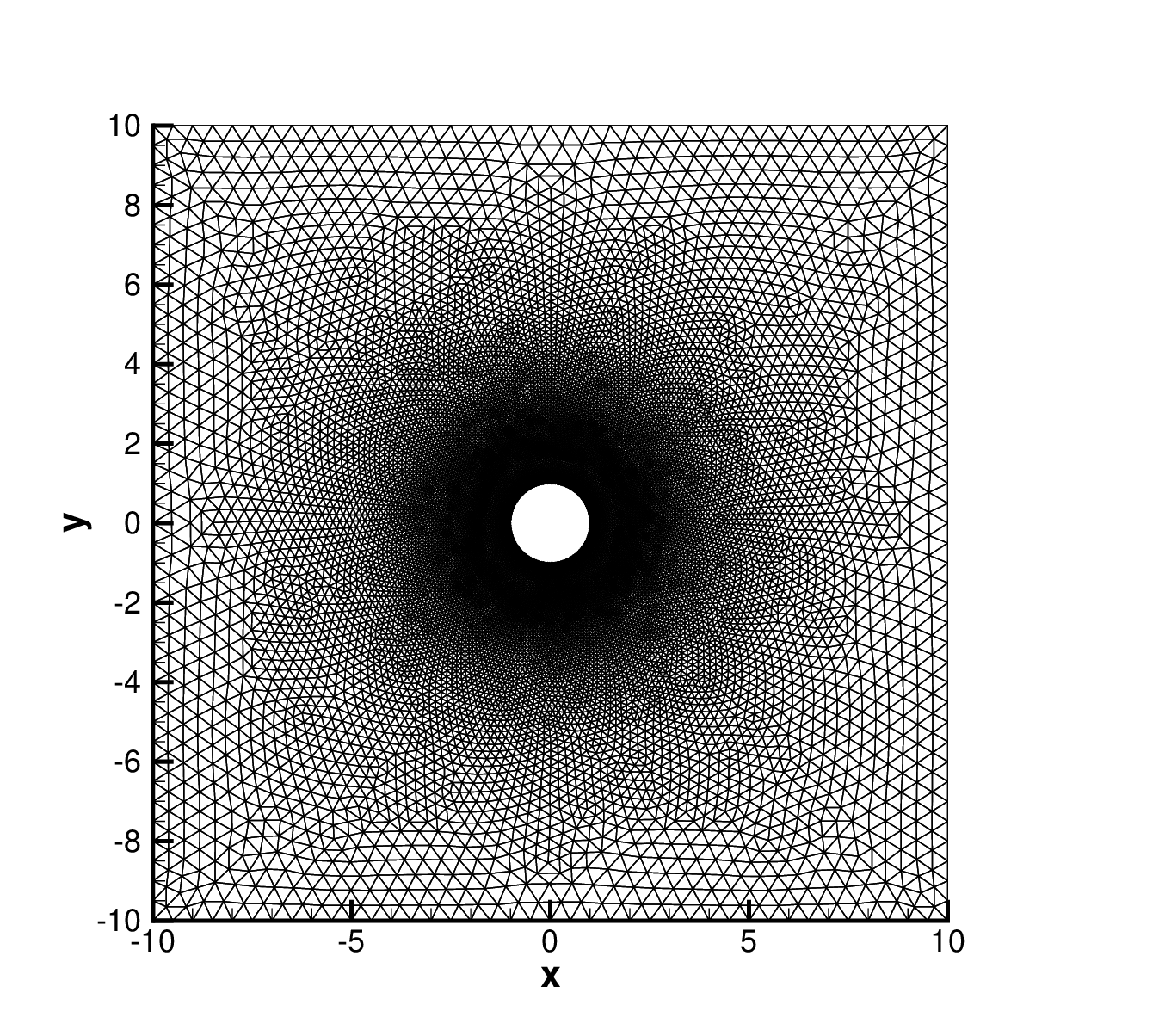}}
  \subfigure[]{\includegraphics[width=0.4\textwidth]{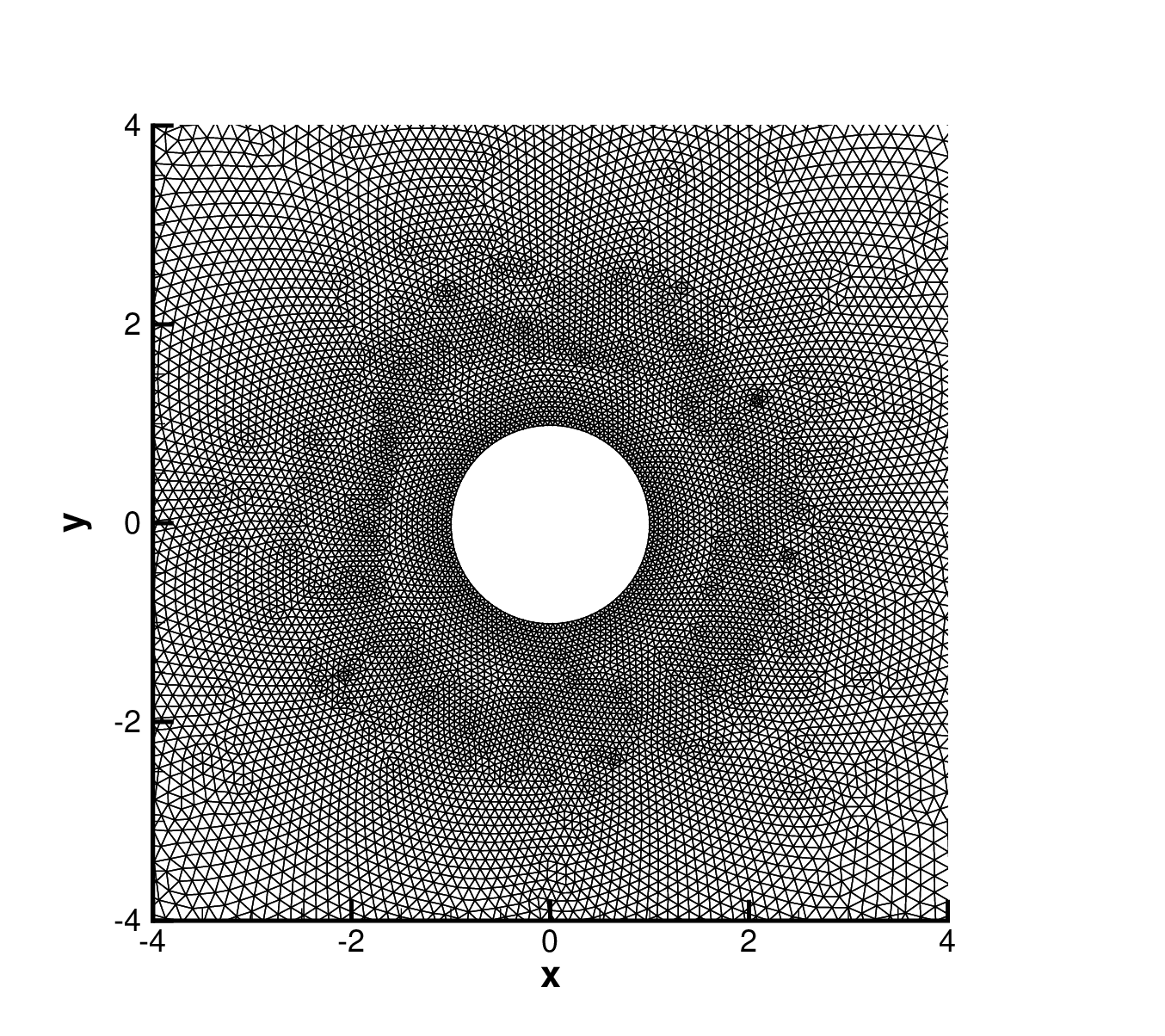}}
  \caption{Horizontal and vertical oscillating cylinders: initial mesh configuration (left) and a zoom close to the cylinder (right).}
  \label{fig:OC-mesh}
\end{figure}

Finally, we conclude the numerical experiments with a fluid-structure interaction problem, where we consider two different test cases of a cylinder
oscillating in a fluid. In this section, we present the results obtained using the polynomial correction for the imposition of moving slip wall boundary conditions,
and compare them qualitatively to those computed without any correction.
Both problems are set in a domain of size $[-10,10]\times[-10,10]$ with a cylinder of radius $R=1$ 
and the same initial mesh configuration shown in figure \ref{fig:OC-mesh} which counts a total number of $N_E=7145$ triangles.

The first case considers the motion of the circular cylinder, in a fluid initially at rest, described in the form of a simple harmonic horizontal oscillation as follows:
\begin{equation}
  x(t) = A \sin(2\pi f t),
\end{equation}
where $A$ denotes the amplitude of the oscillation and $f$ the frequency.
In this case, we set $A=0.1$ and $f=0.1$, and the final time is $t_f=10$ such that a whole period is completed.

\begin{figure}
  \centering
  \subfigure[$t=0.5$ s]{\includegraphics[width=0.4\textwidth]{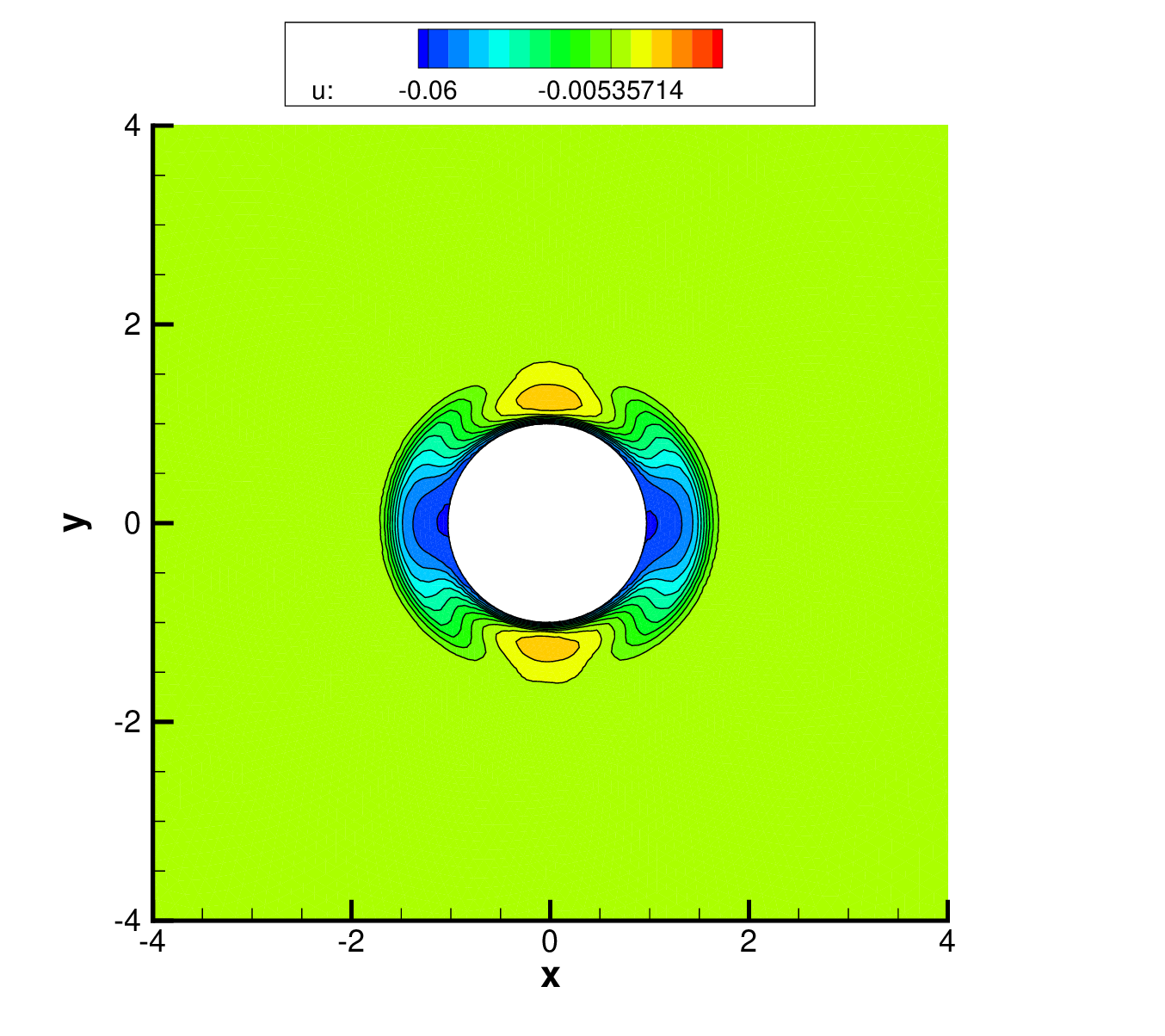}}
  \subfigure[$t=0.5$ s]{\includegraphics[width=0.4\textwidth]{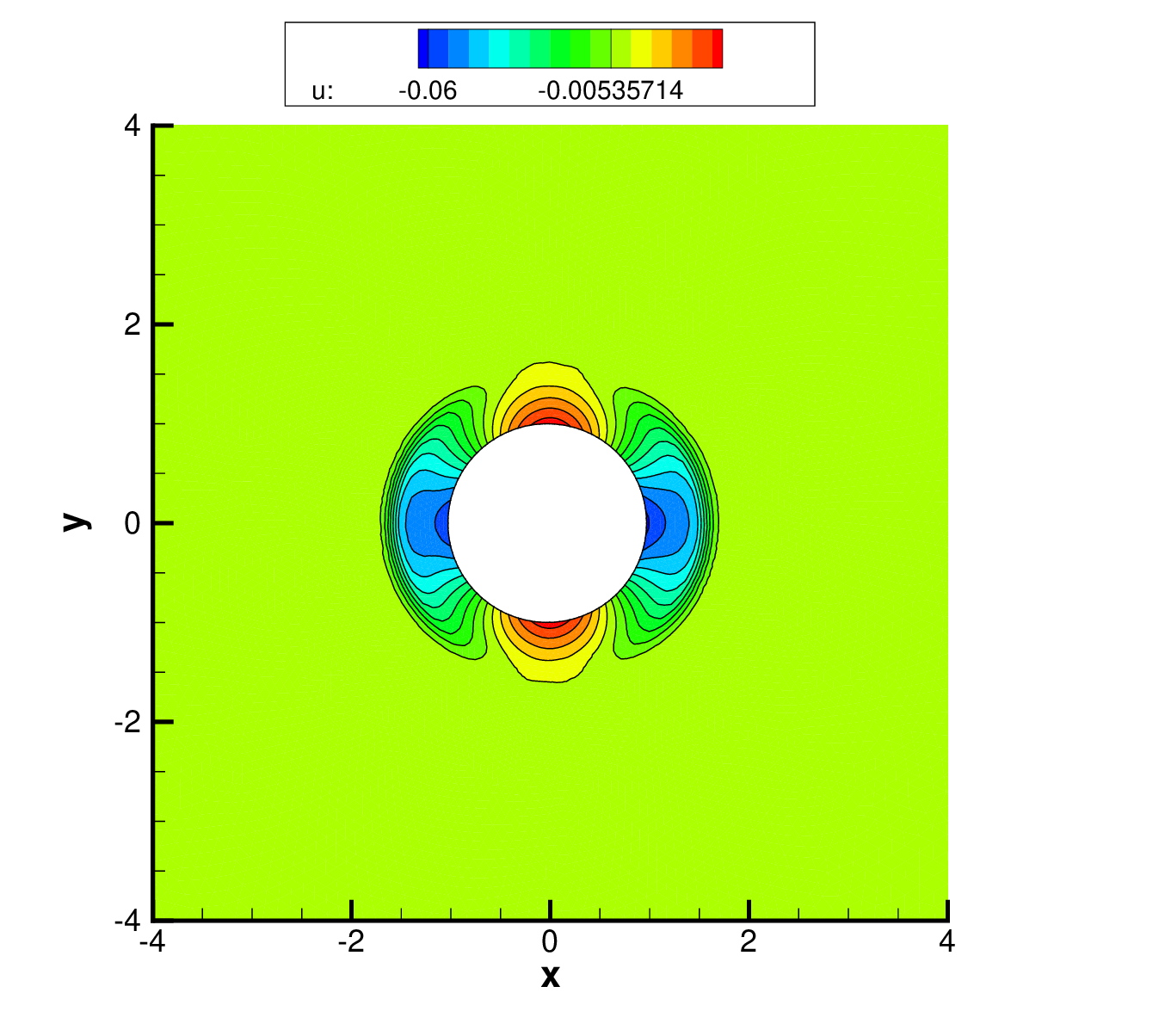}}
  \subfigure[$t=5$ s]{\includegraphics[width=0.4\textwidth]{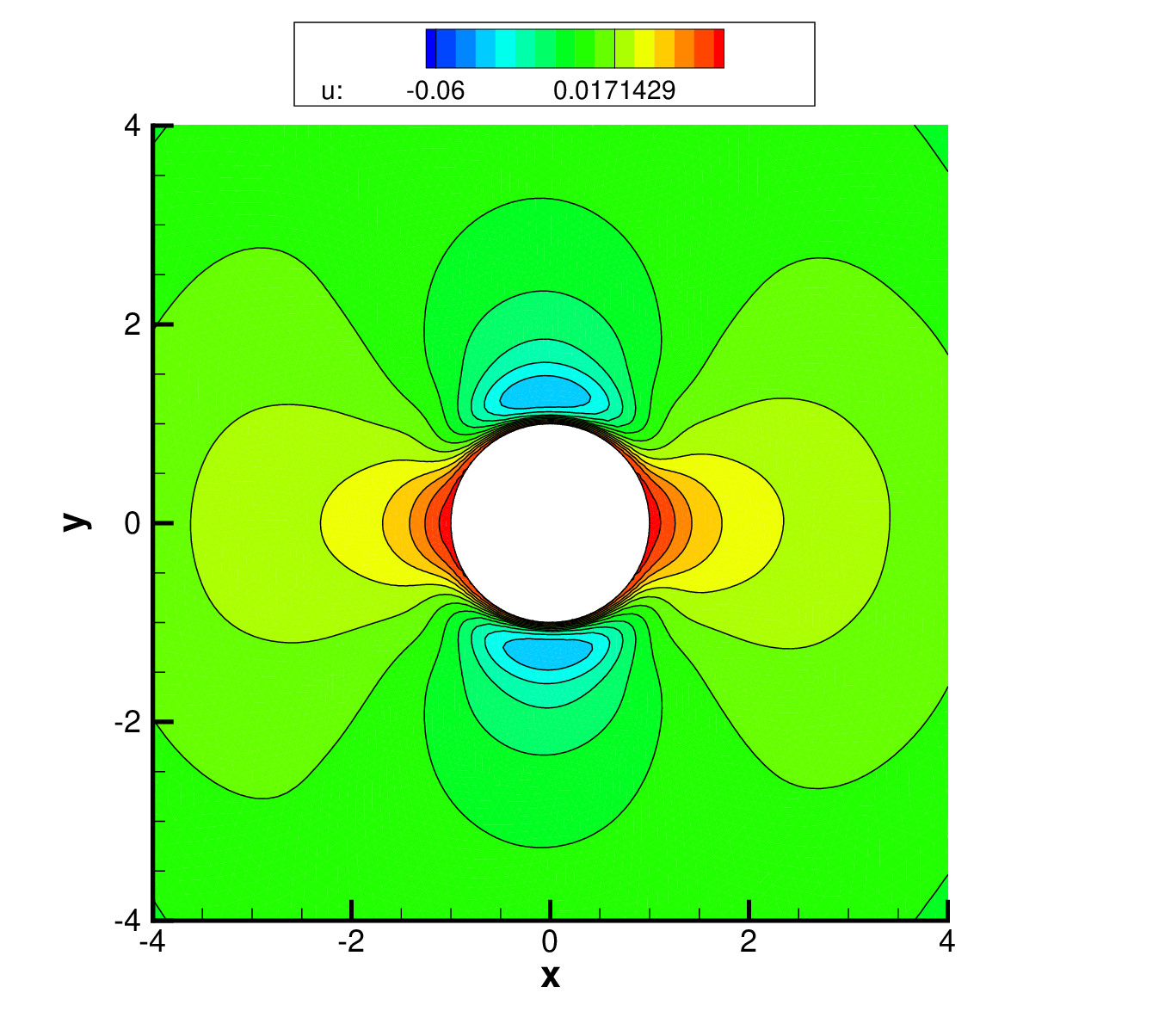}}
  \subfigure[$t=5$ s]{\includegraphics[width=0.4\textwidth]{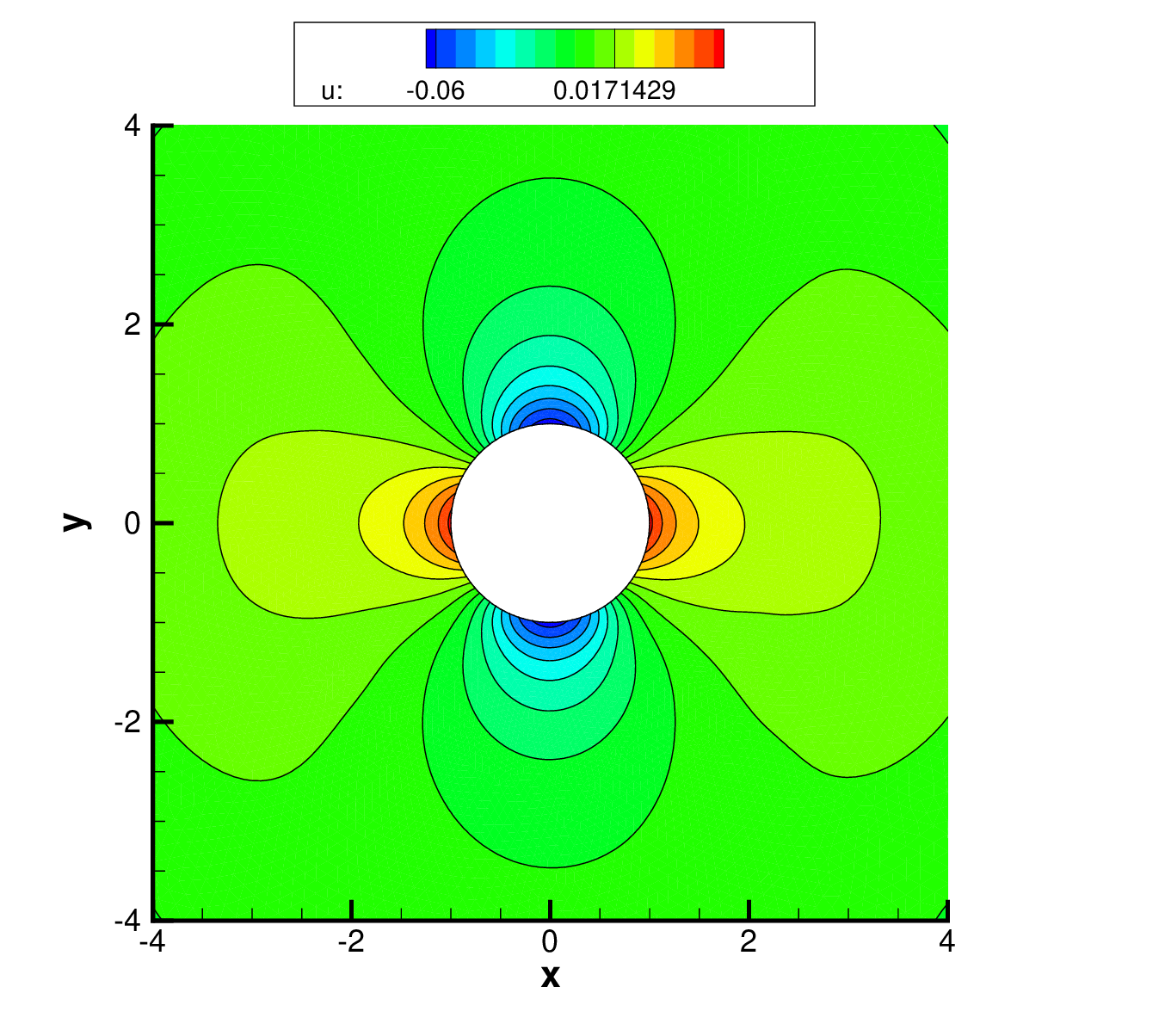}}
  \subfigure[$t=10$ s]{\includegraphics[width=0.4\textwidth]{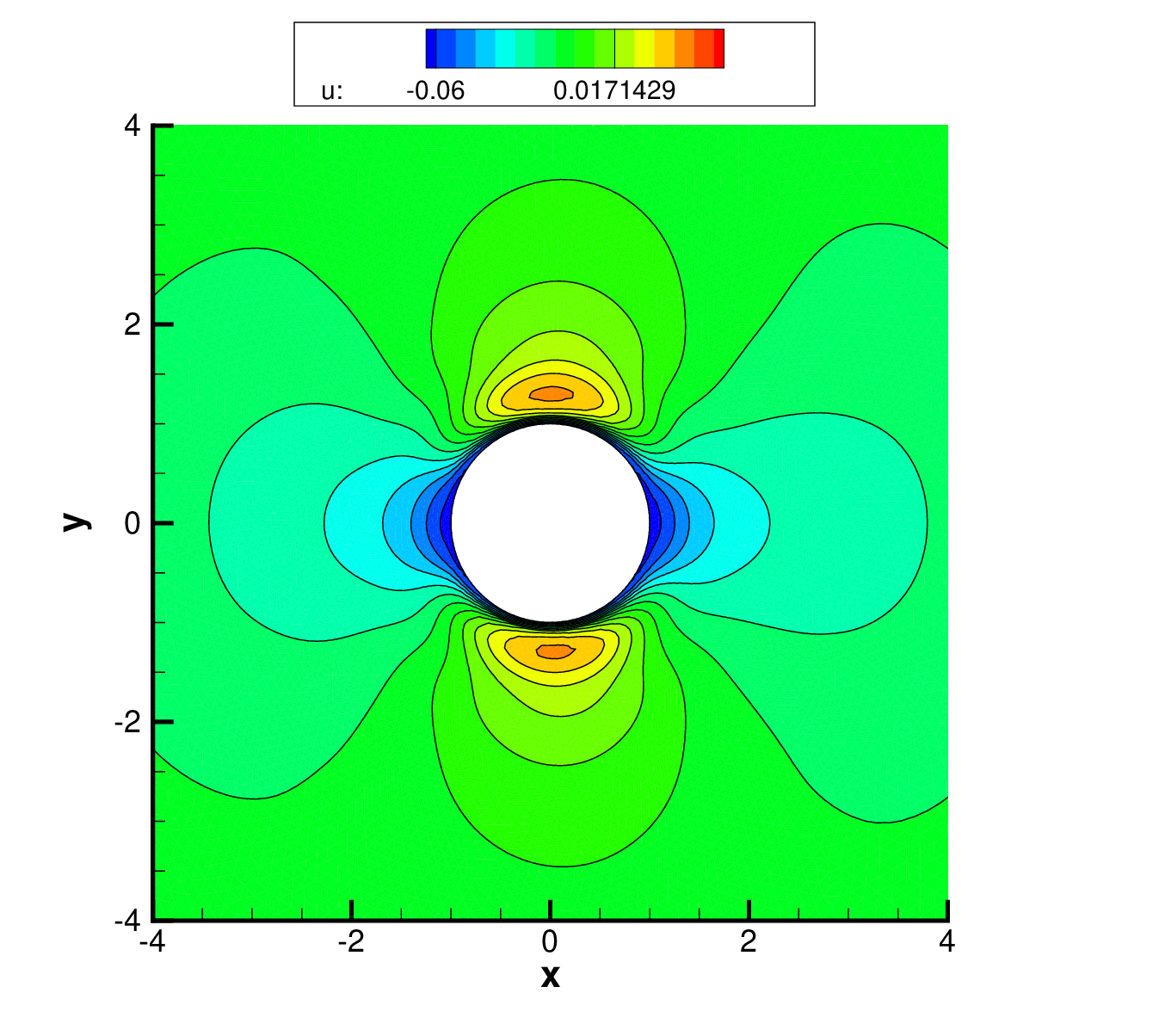}}
  \subfigure[$t=10$ s]{\includegraphics[width=0.4\textwidth]{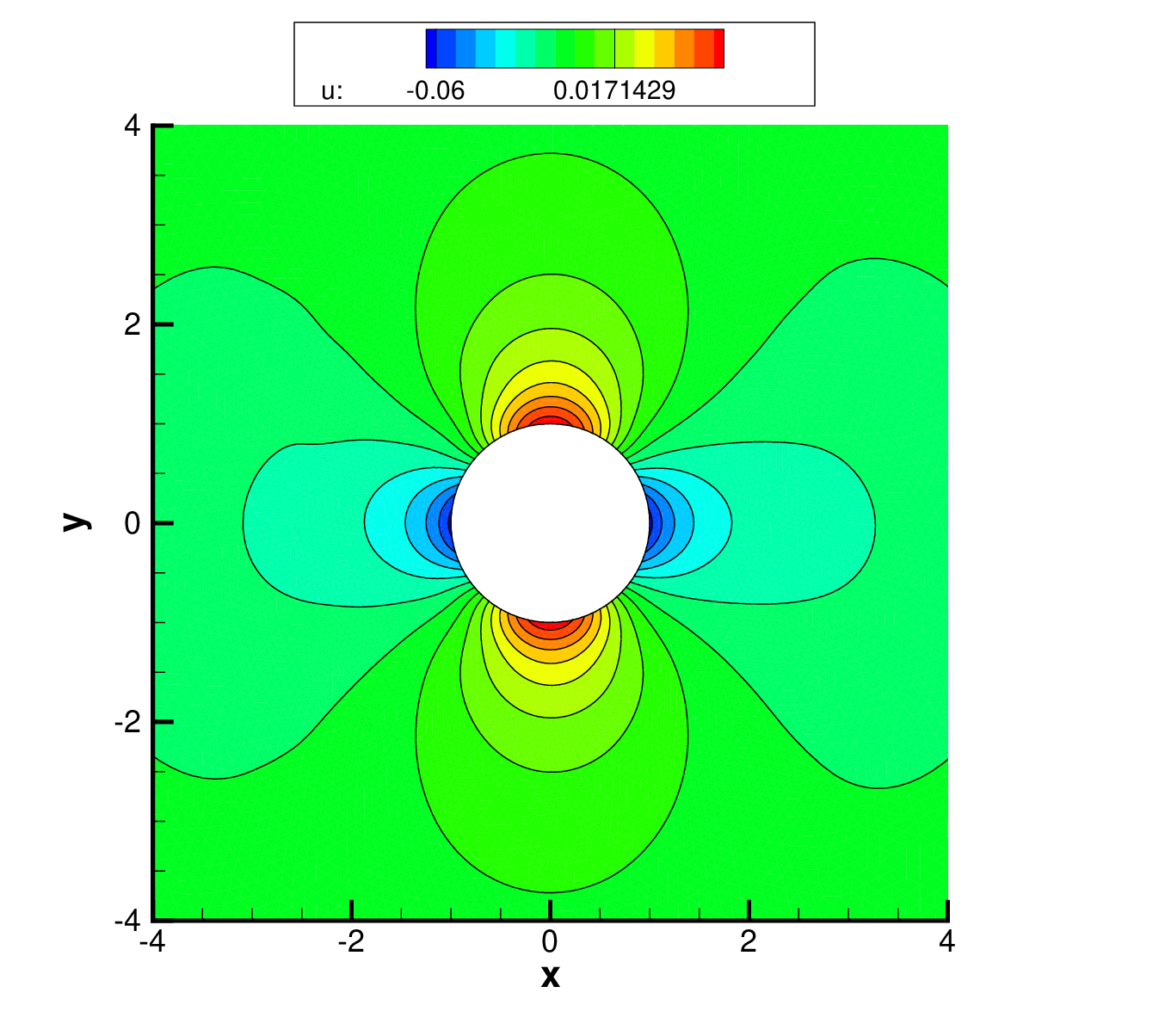}}
  \caption{Horizontal oscillating cylinder: velocity component $u$ distribution without correction (left) and with correction (right).
  The solution is plot at different times: $t=0.5$ s (top), $t=5$ s (middle) and $t=10$ s (bottom).}
  \label{fig:OC-H-velU}
\end{figure}

\begin{figure}
  \centering
  \subfigure[$t=0.5$ s]{\includegraphics[width=0.4\textwidth]{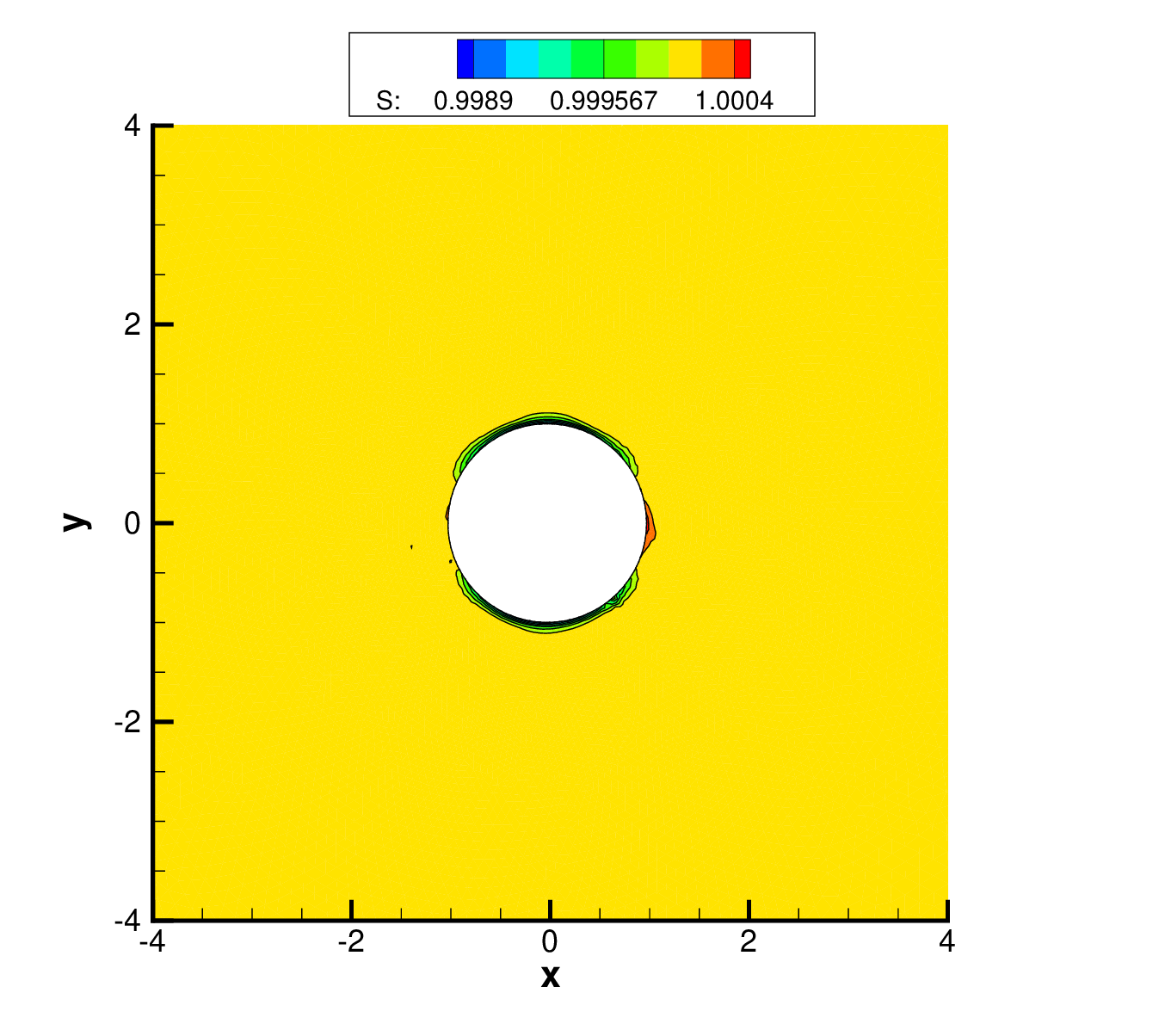}}
  \subfigure[$t=0.5$ s]{\includegraphics[width=0.4\textwidth]{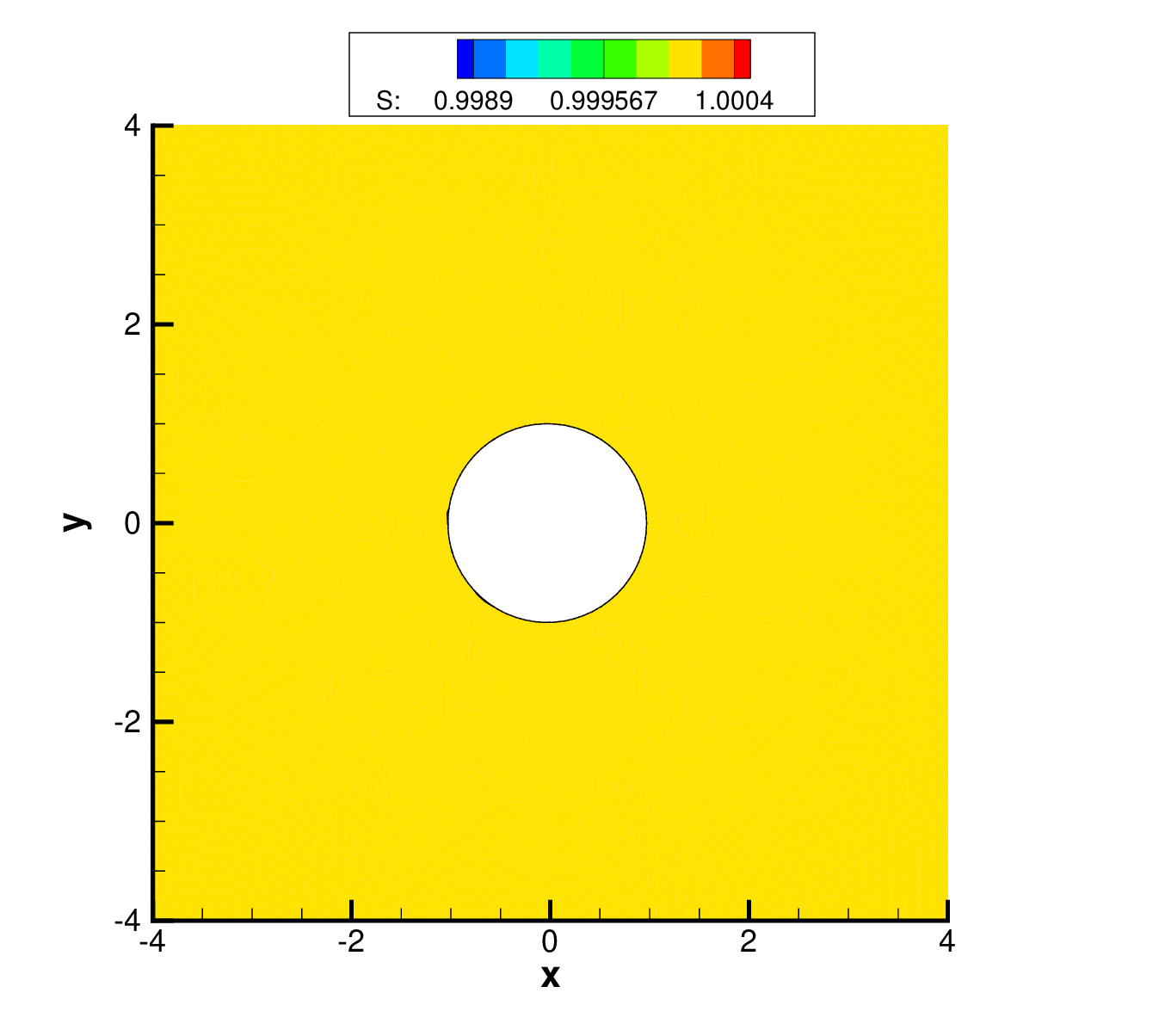}}
  \subfigure[$t=5$ s]{\includegraphics[width=0.4\textwidth]{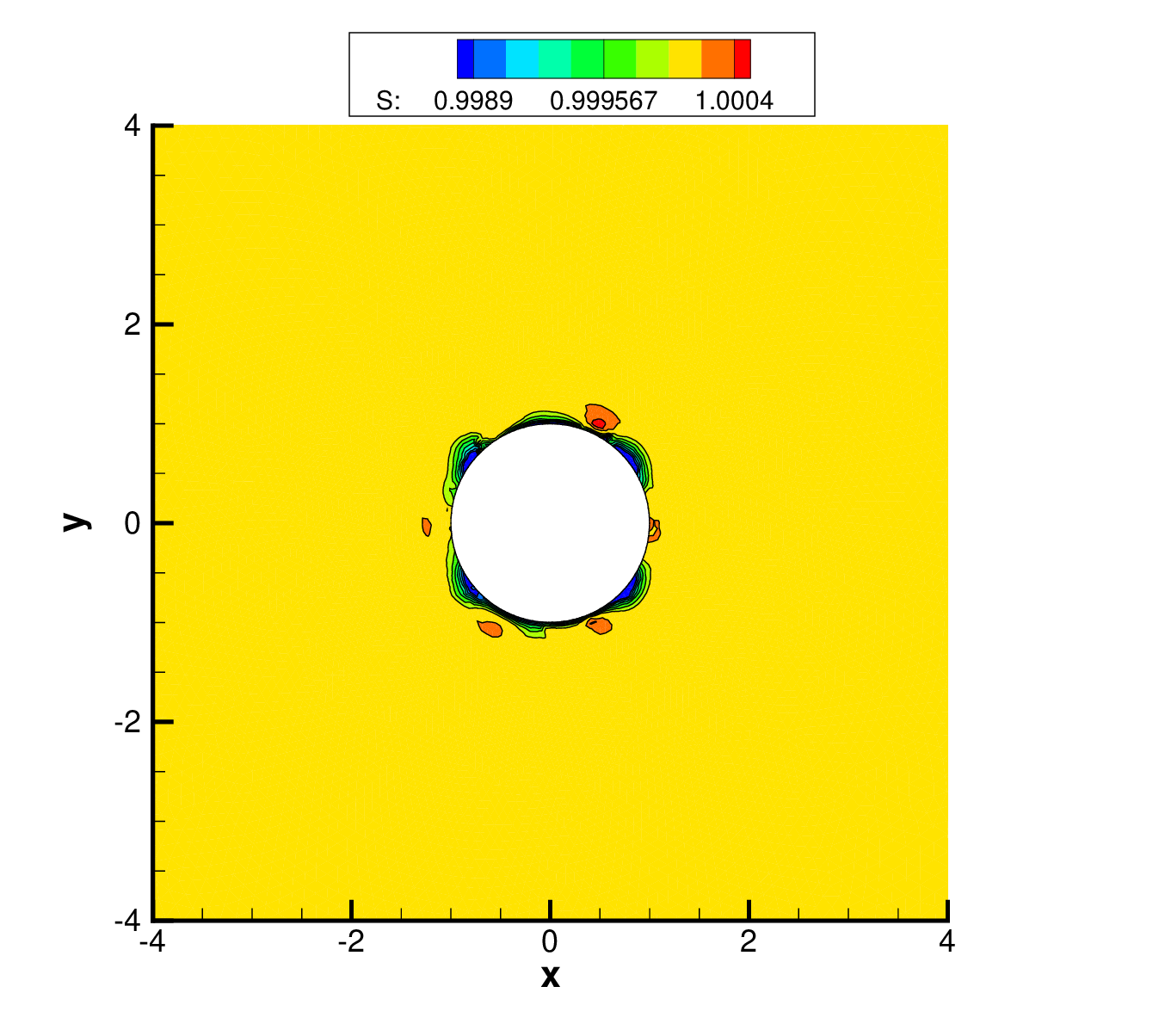}}
  \subfigure[$t=5$ s]{\includegraphics[width=0.4\textwidth]{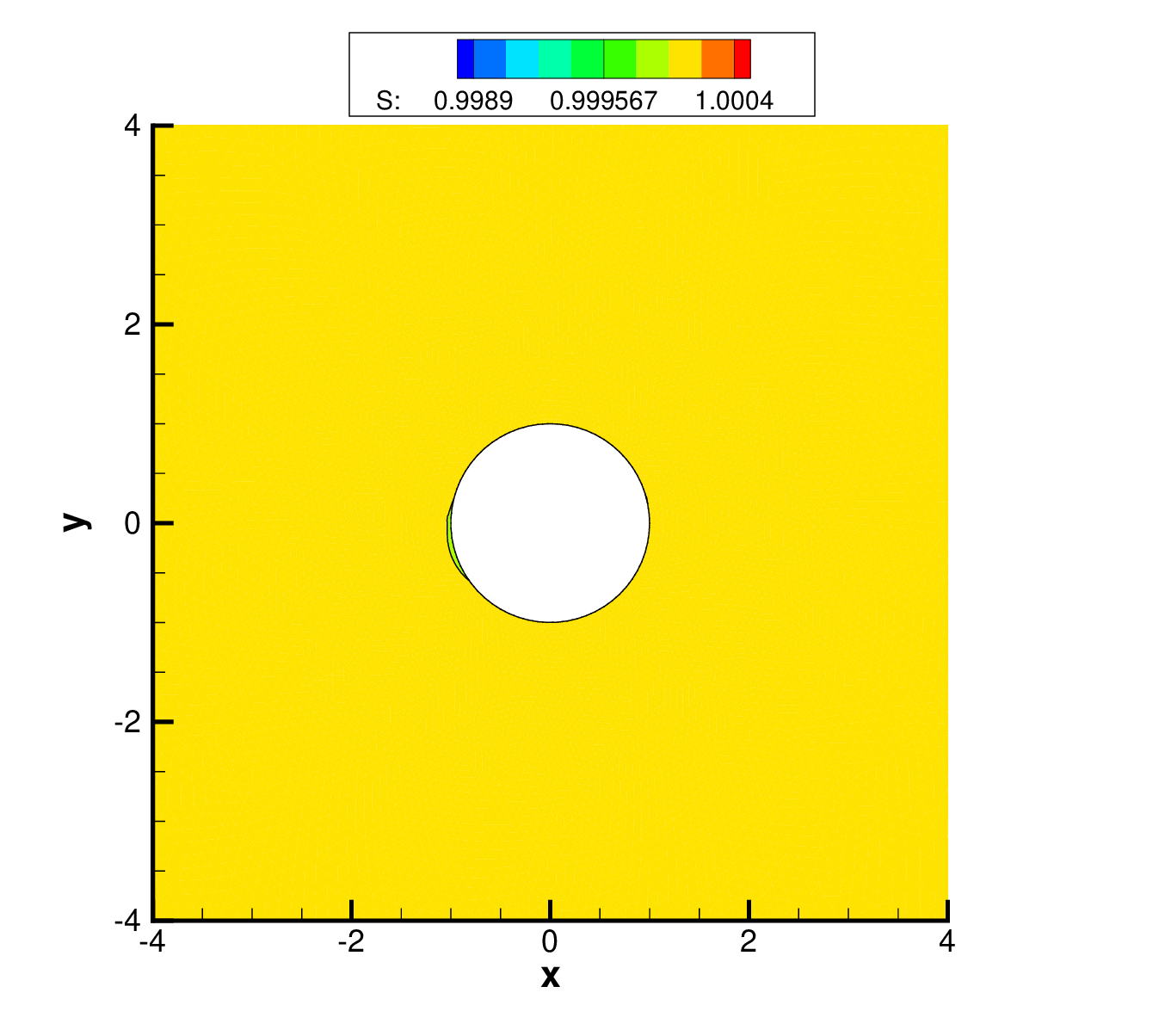}}
  \subfigure[$t=10$ s]{\includegraphics[width=0.4\textwidth]{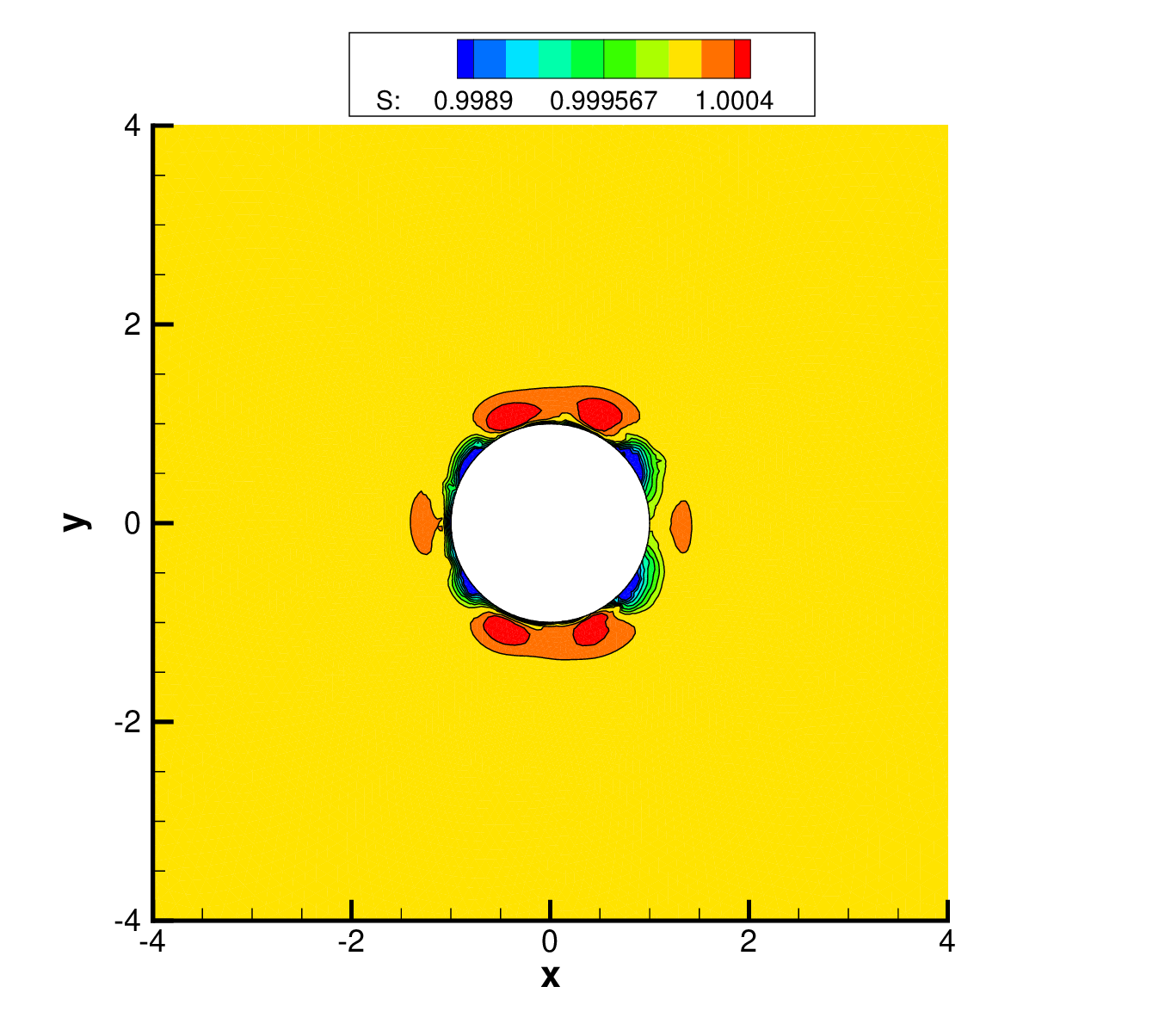}}
  \subfigure[$t=10$ s]{\includegraphics[width=0.4\textwidth]{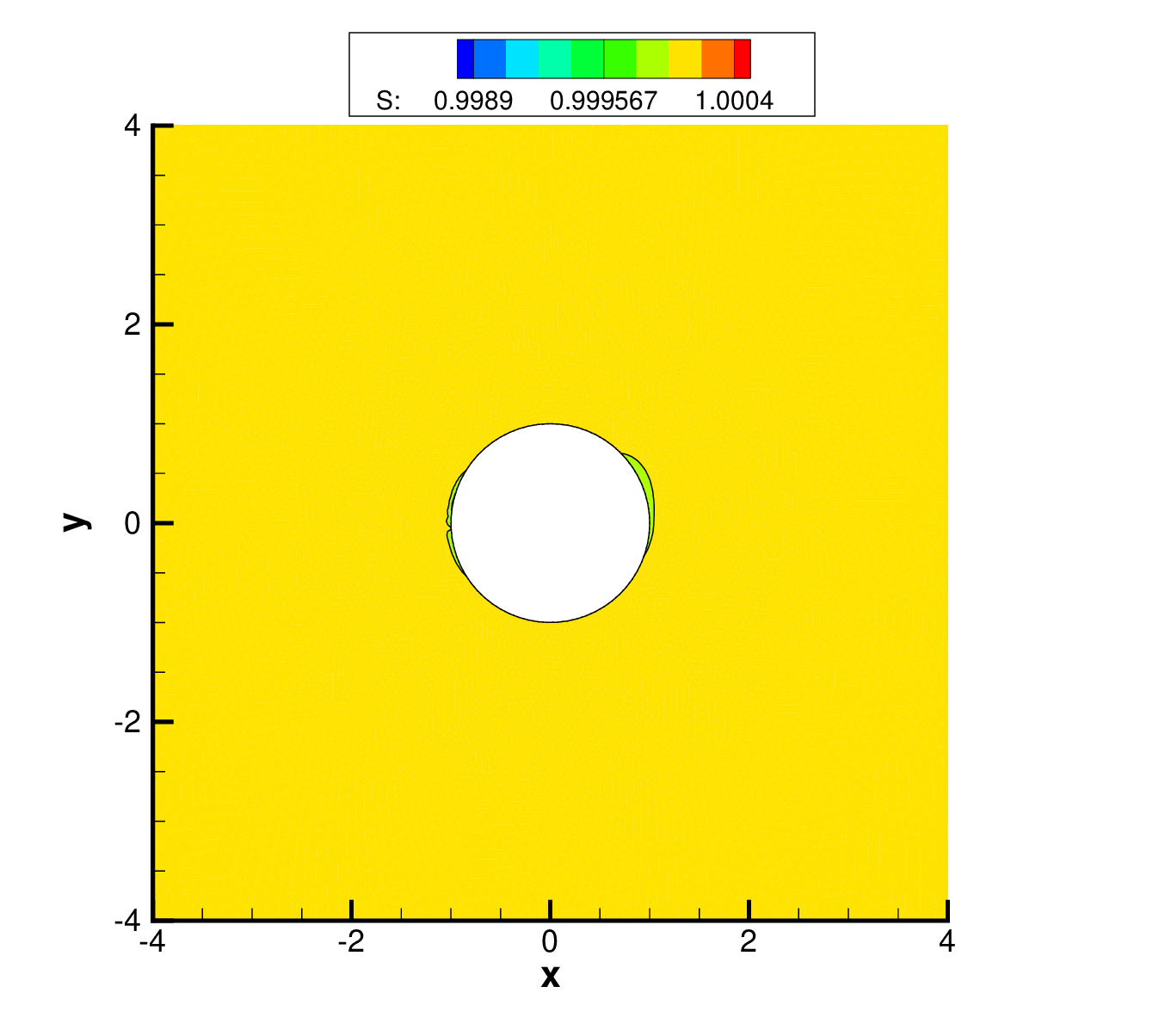}}
  \caption{Horizontal oscillating cylinder: entropy $S$ distribution without correction (left) and with correction (right).
  The solution is plot at different times: $t=0.5$ s (top), $t=5$ s (middle) and $t=10$ s (bottom).}
  \label{fig:OC-H-entropy}
\end{figure}

In figure \ref{fig:OC-H-velU} we depict the velocity field $u$ at different times for the horizontal oscillating cylinder, and
compare the results obtained with and without the polynomial correction.
The left column shows the results obtained without the correction, while the right column shows the results obtained with the correction.
The results demonstrate that the polynomial correction allows us to recover a much smoother and symmetric velocity field, 
while the results obtained without the correction present spurious errors close to the boundary.
In particular, the geometry errors dominates the numerical results that present a thick layer of errors close to the boundary.
The same conclusions can be drawn by looking at the entropy distribution $S=p/\rho^\gamma$ in figure \ref{fig:OC-H-entropy},
where the spurious entropy production is clearly visible in the left column.
It should be noticed that the results with and without the correction are compared using the same color map to highlight the differences between the two cases and to provide a fair comparison. 

The second test case is concerned with the motion of the cylinder, immersed in a horizontally-moving fluid with speed $u=0.25$, 
described by a harmonic vertical oscillation as follows:
\begin{equation}
  y(t) = A \cos(2\pi f t),
\end{equation}
where we set $A=0.05$ and $f=0.25$, and consider a simulation final time $t_f=8$.
Similar test cases have been performed in the literature \cite{wang2011improved} in the context of viscous flows, and
here we replicate them for inviscid flows for the sake of comparison of the new approach with the standard one in the context of moving slip wall conditions,
and to show the potential of the new approach to deal with the simulation of fluid-structure interaction problems with high order methods.

\begin{figure}
  \centering
  \subfigure[$t=0.5$ s]{\includegraphics[width=0.4\textwidth]{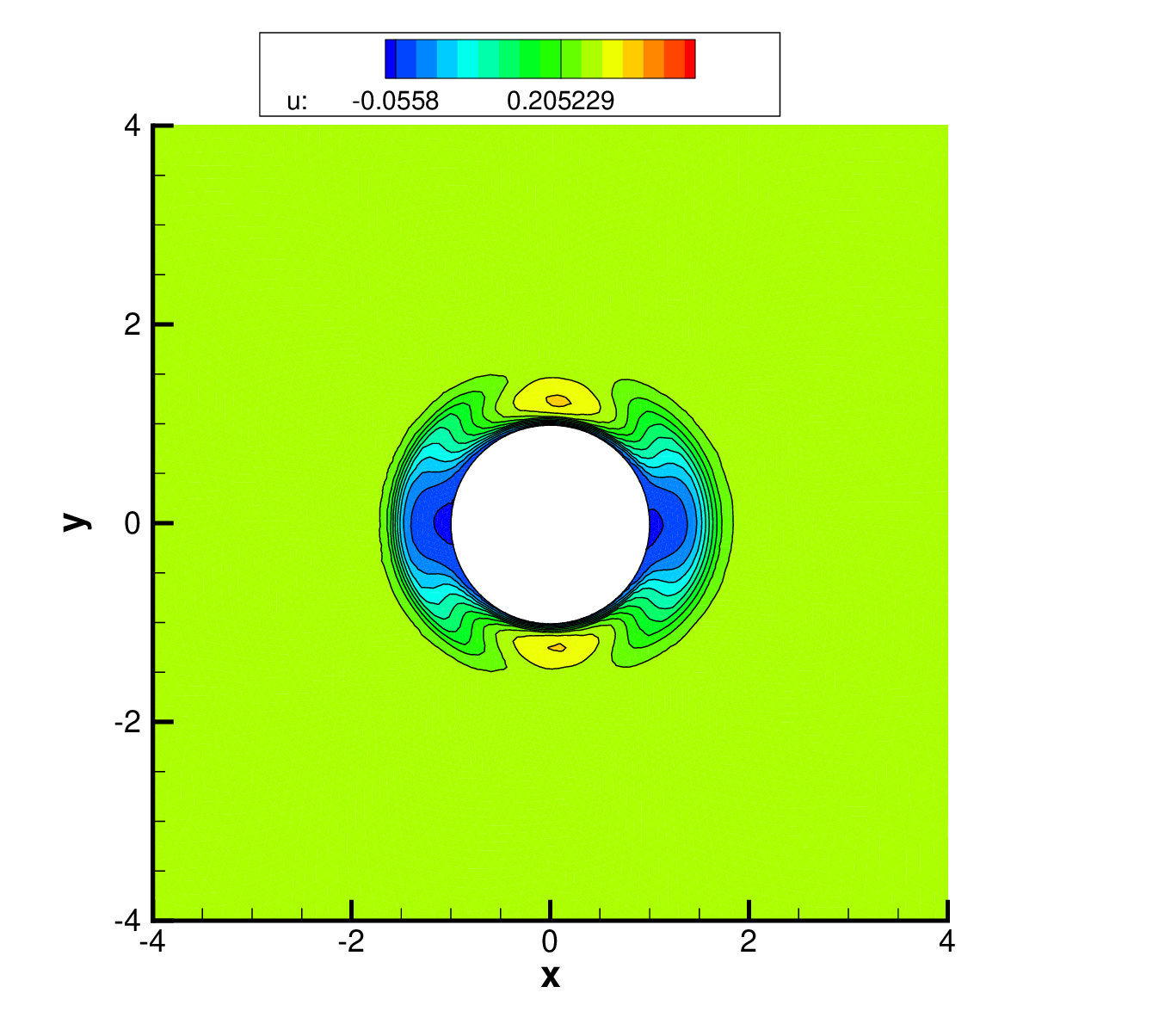}}
  \subfigure[$t=0.5$ s]{\includegraphics[width=0.4\textwidth]{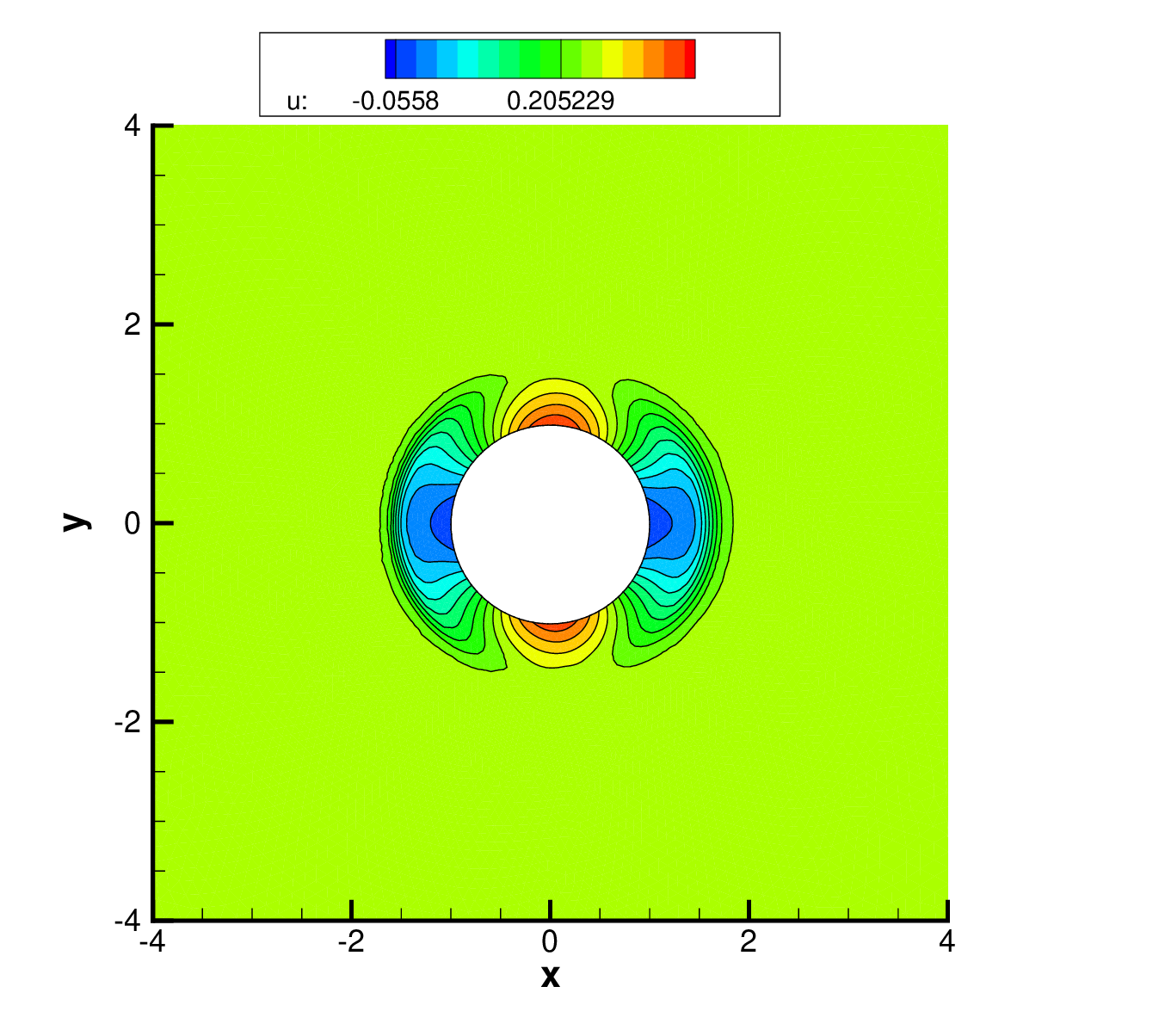}}
  \subfigure[$t=4$ s]{\includegraphics[width=0.4\textwidth]{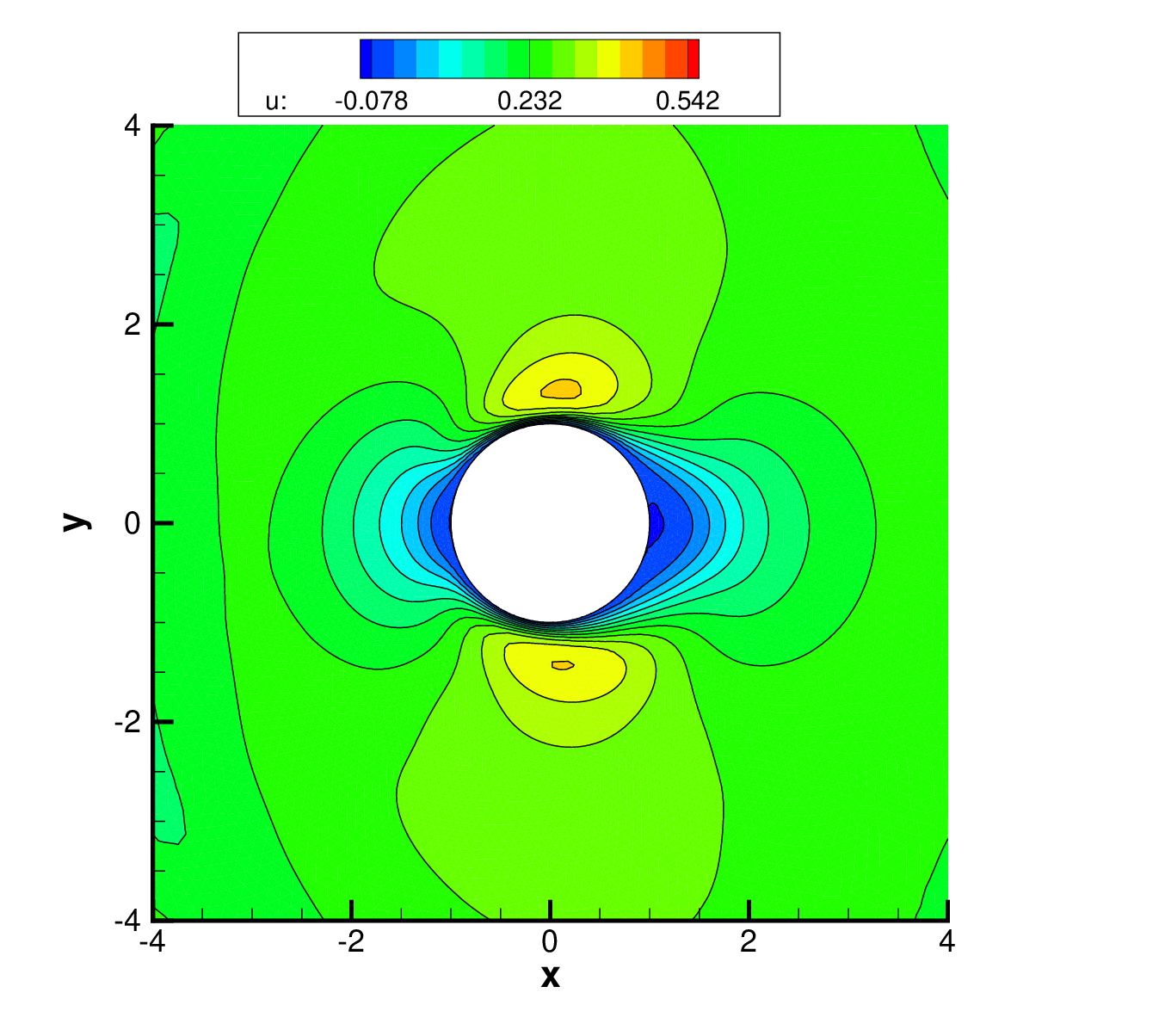}}
  \subfigure[$t=4$ s]{\includegraphics[width=0.4\textwidth]{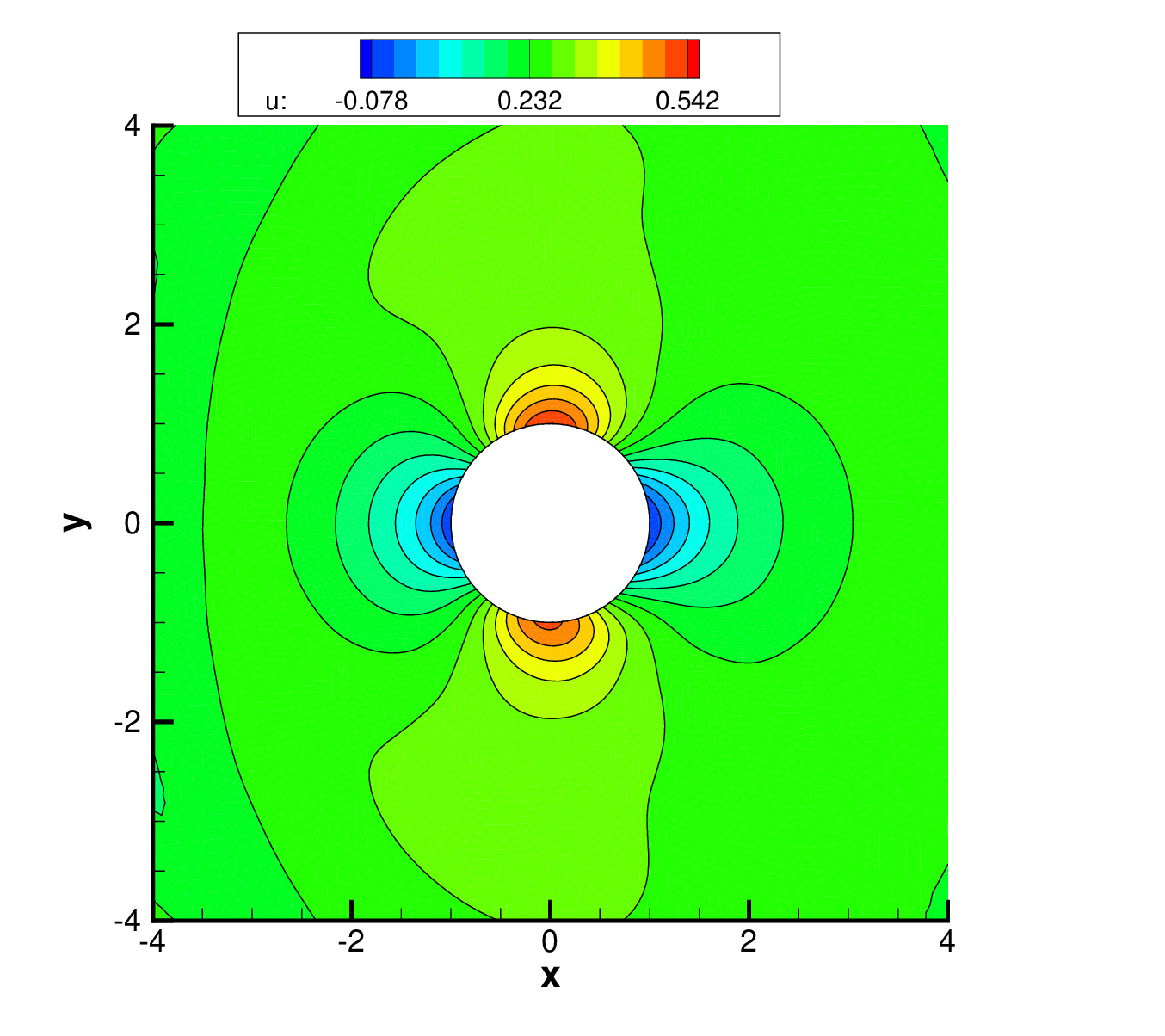}}
  \subfigure[$t=8$ s]{\includegraphics[width=0.4\textwidth]{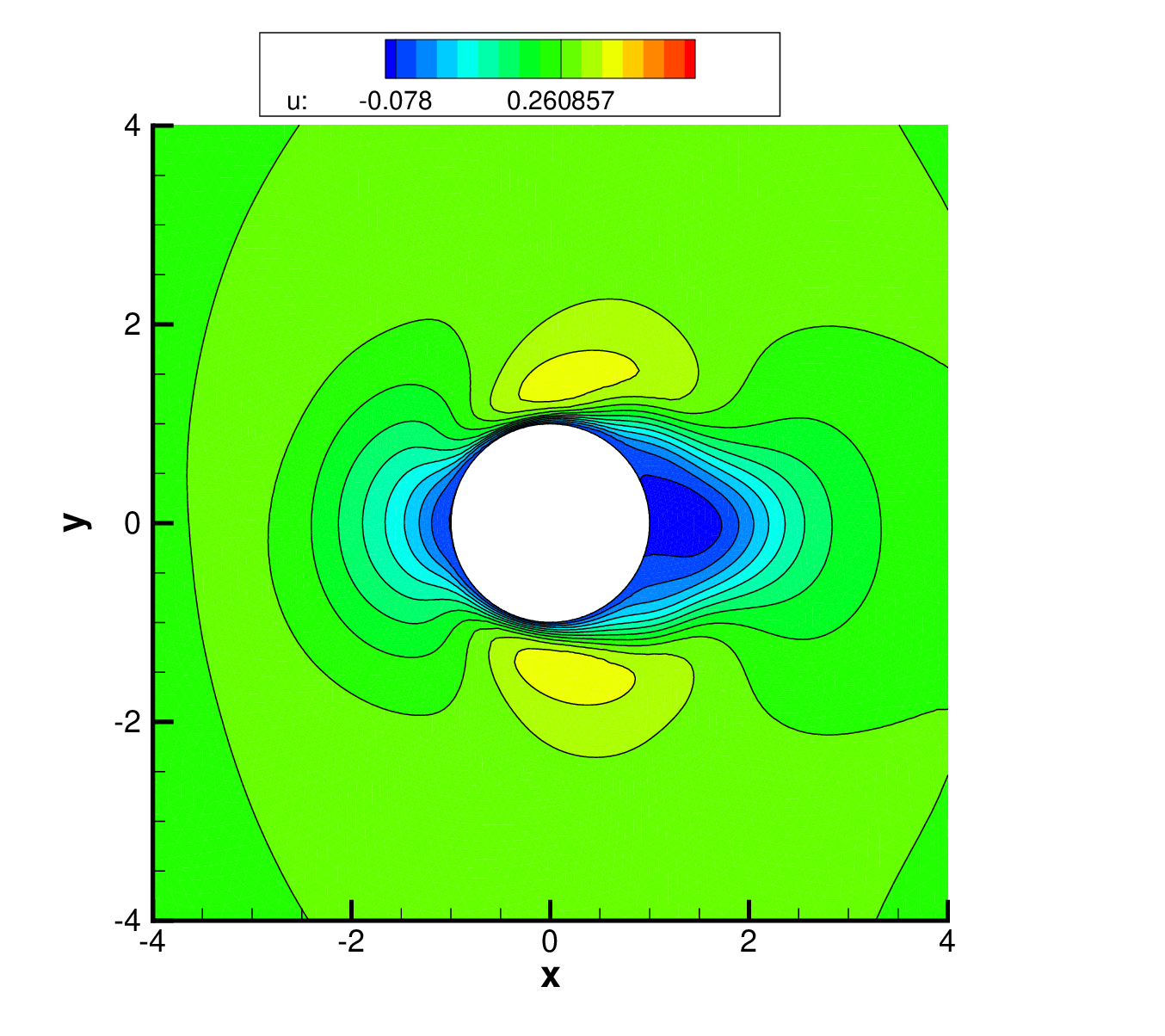}}
  \subfigure[$t=8$ s]{\includegraphics[width=0.4\textwidth]{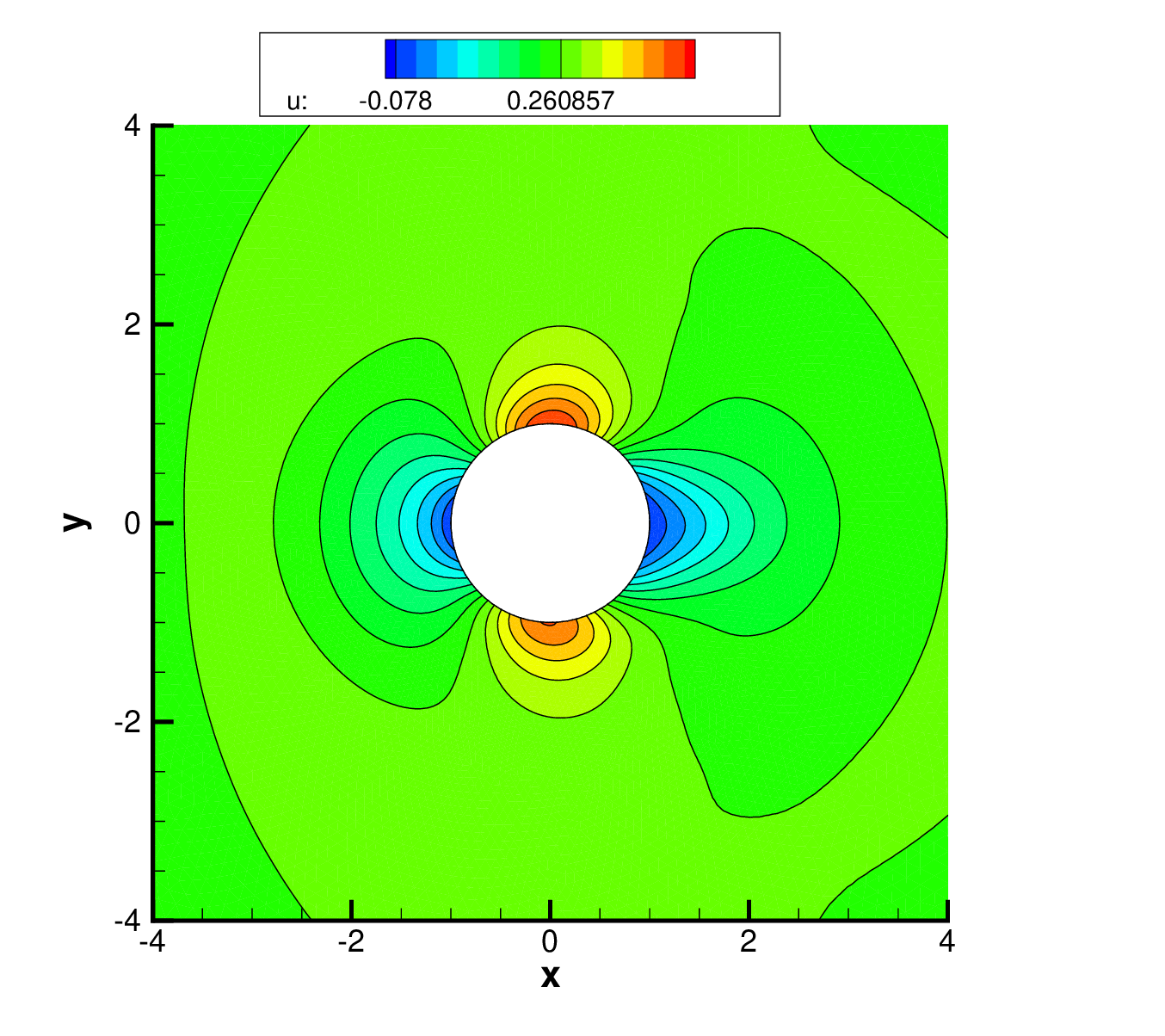}}
  \caption{Vertical oscillating cylinder: velocity component $u$ distribution without correction (left) and with correction (right).
  The solution is plot at different times: $t=0.5$ s (top), $t=4$ s (middle) and $t=8$ s (bottom).}
  \label{fig:OC-V-velU}
\end{figure}

\begin{figure}
  \centering
  \subfigure[$t=0.5$ s]{\includegraphics[width=0.4\textwidth]{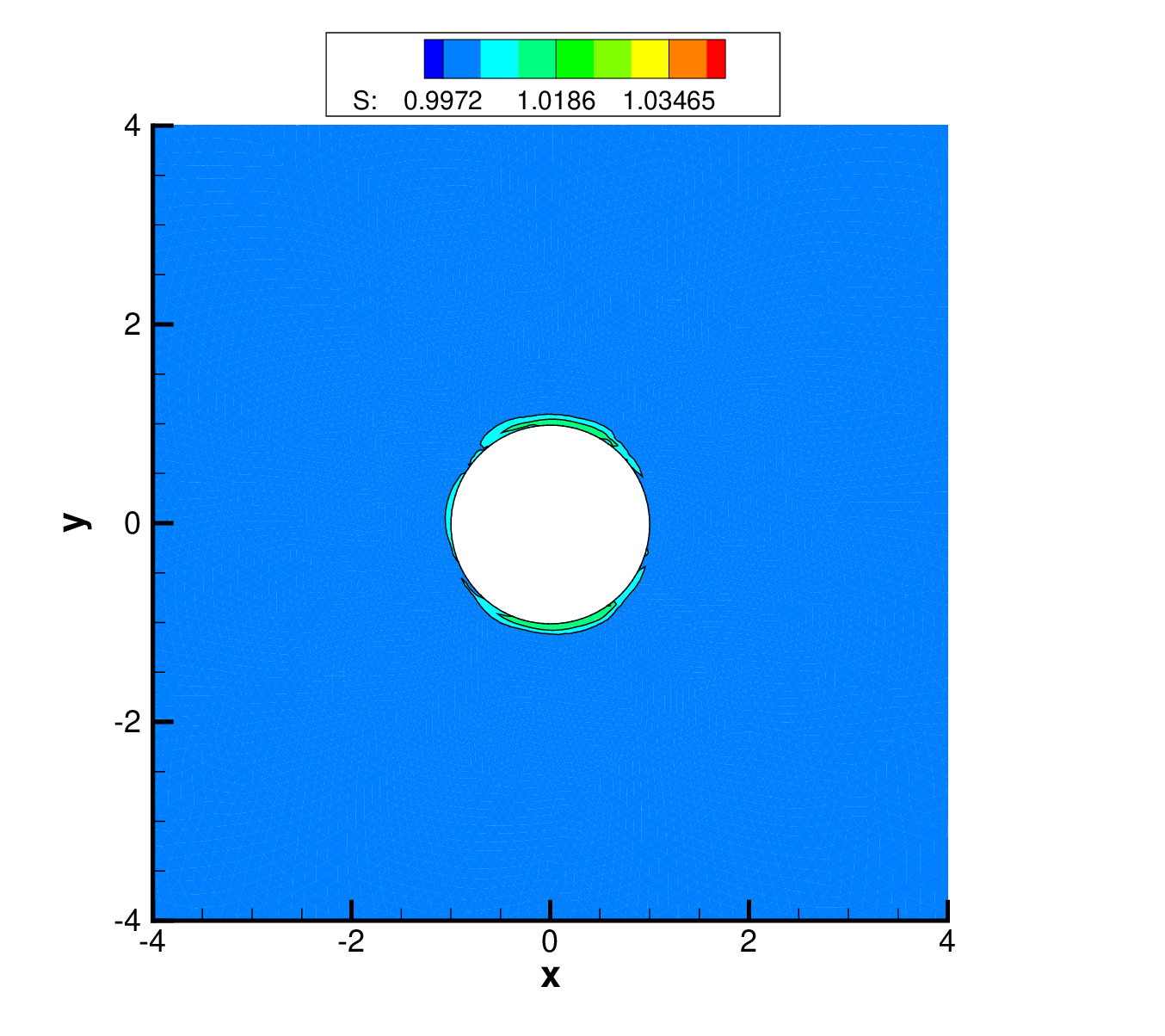}}
  \subfigure[$t=0.5$ s]{\includegraphics[width=0.4\textwidth]{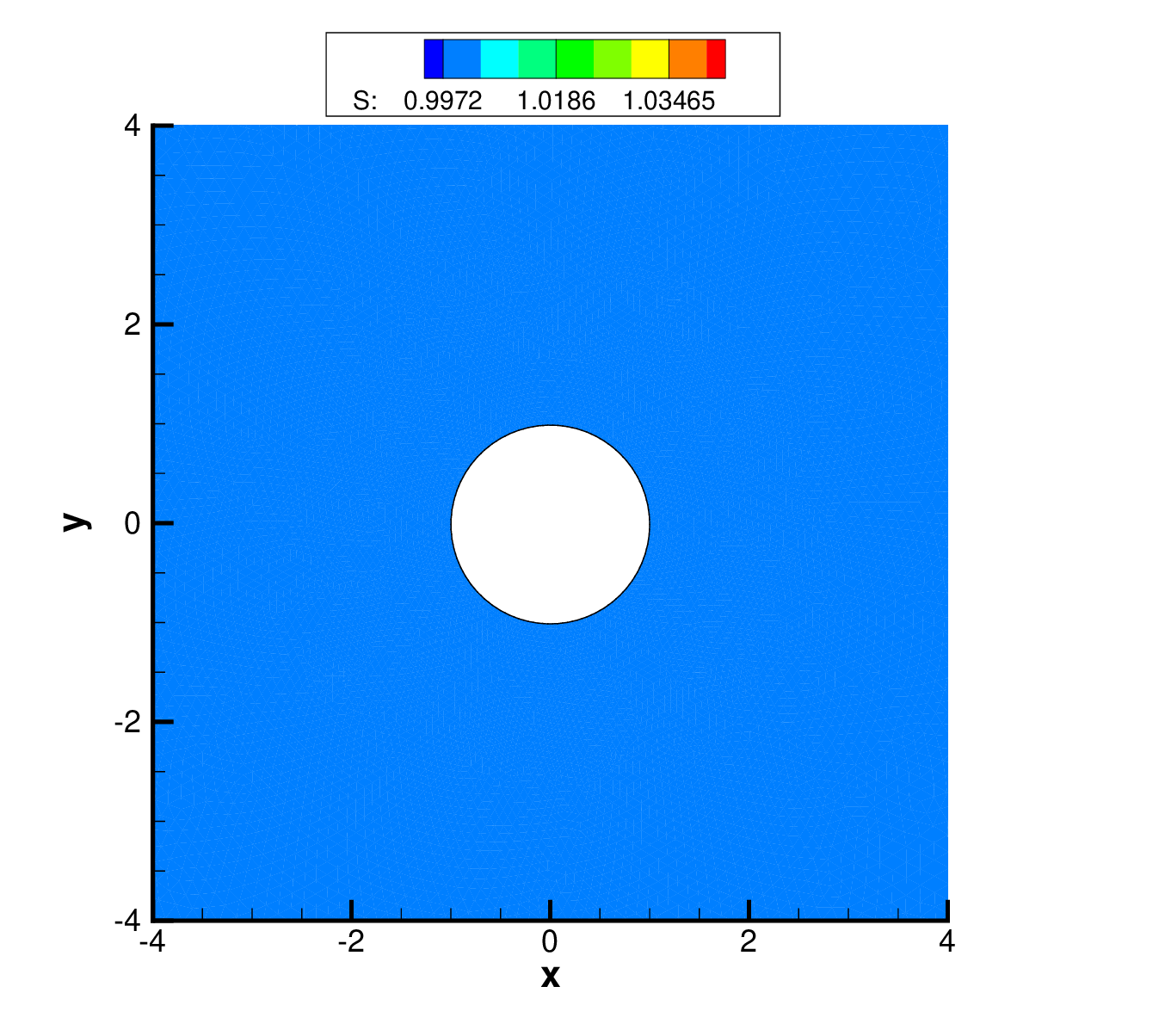}}
  \subfigure[$t=4$ s]{\includegraphics[width=0.4\textwidth]{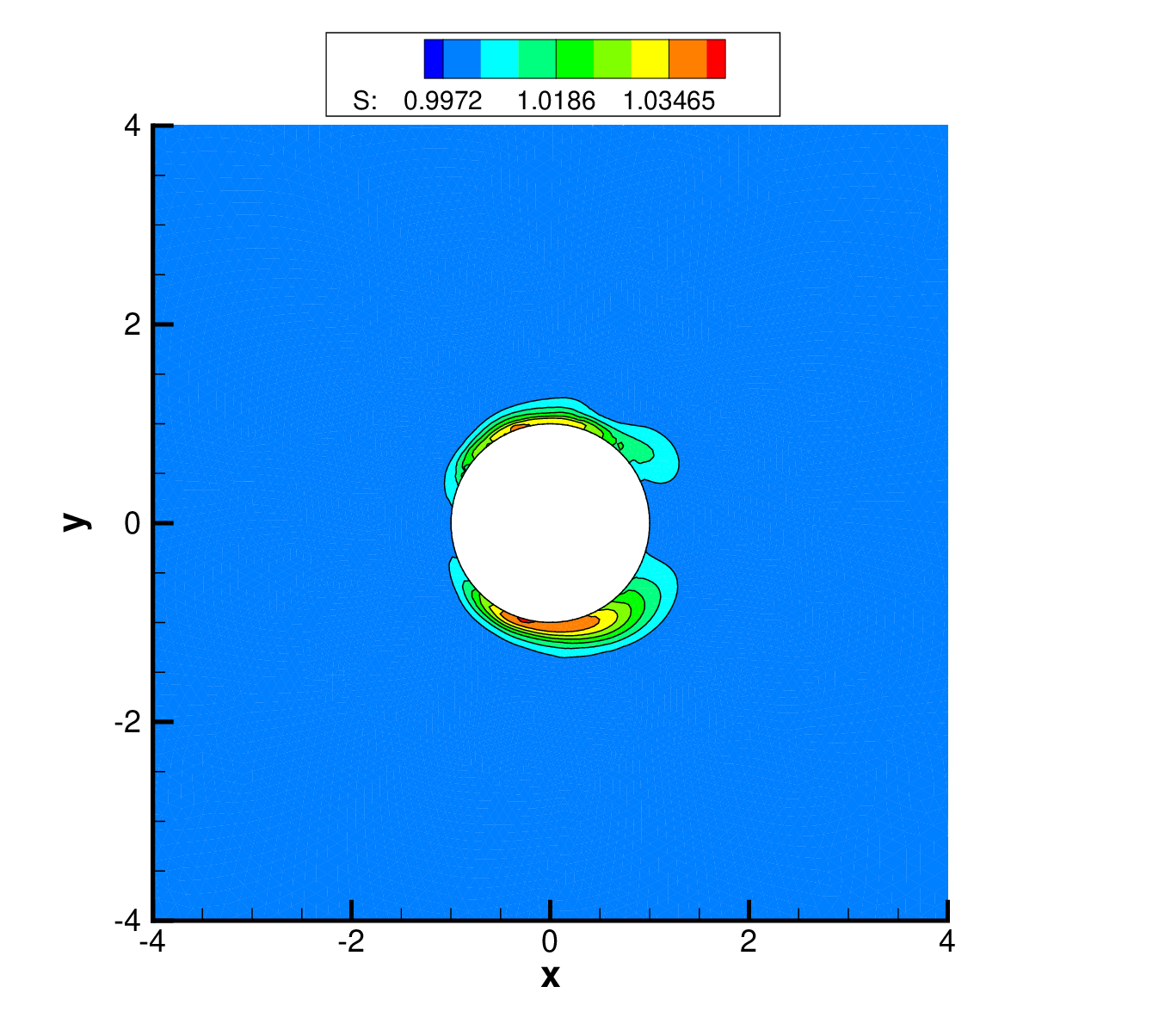}}
  \subfigure[$t=4$ s]{\includegraphics[width=0.4\textwidth]{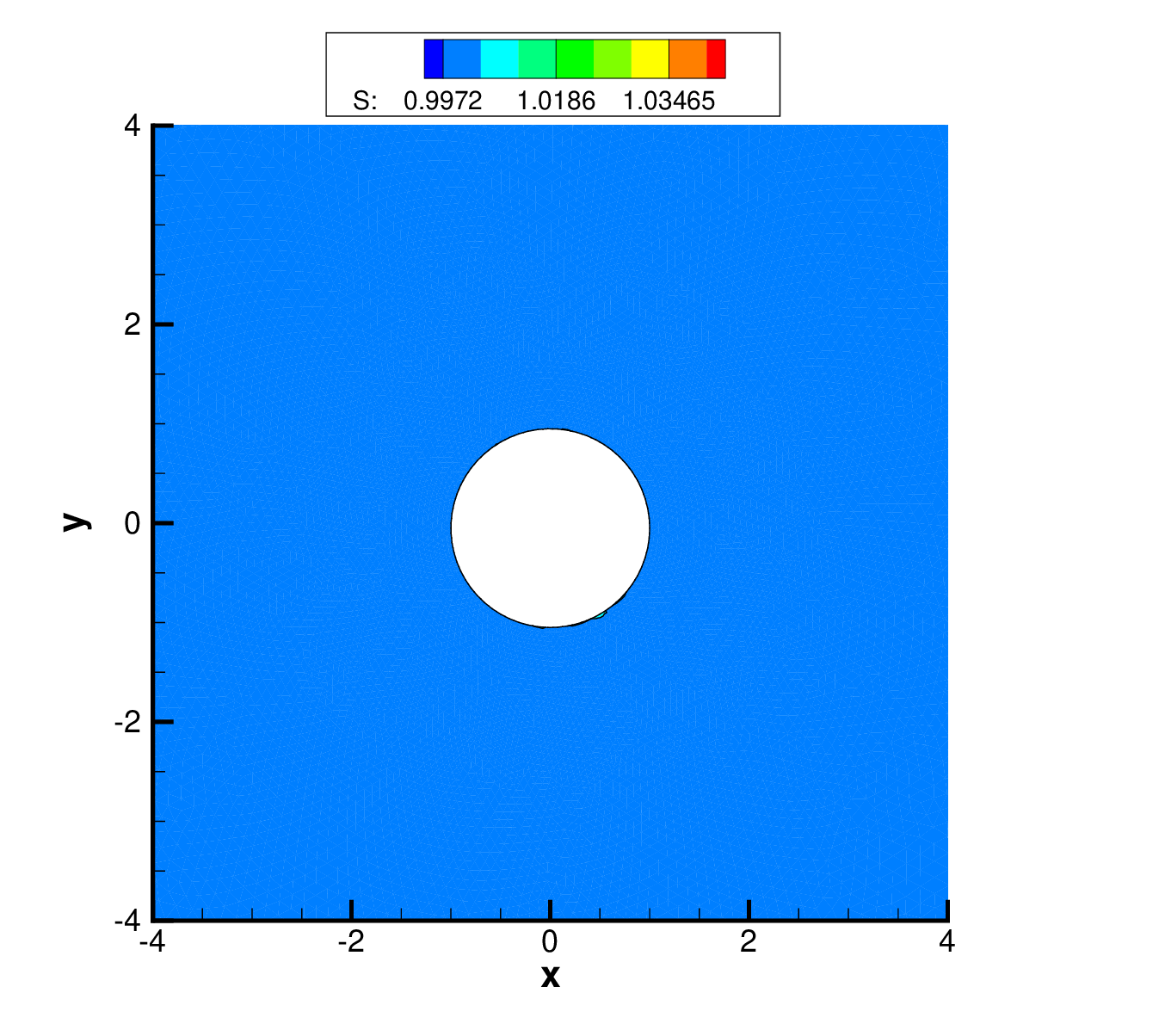}}
  \subfigure[$t=8$ s]{\includegraphics[width=0.4\textwidth]{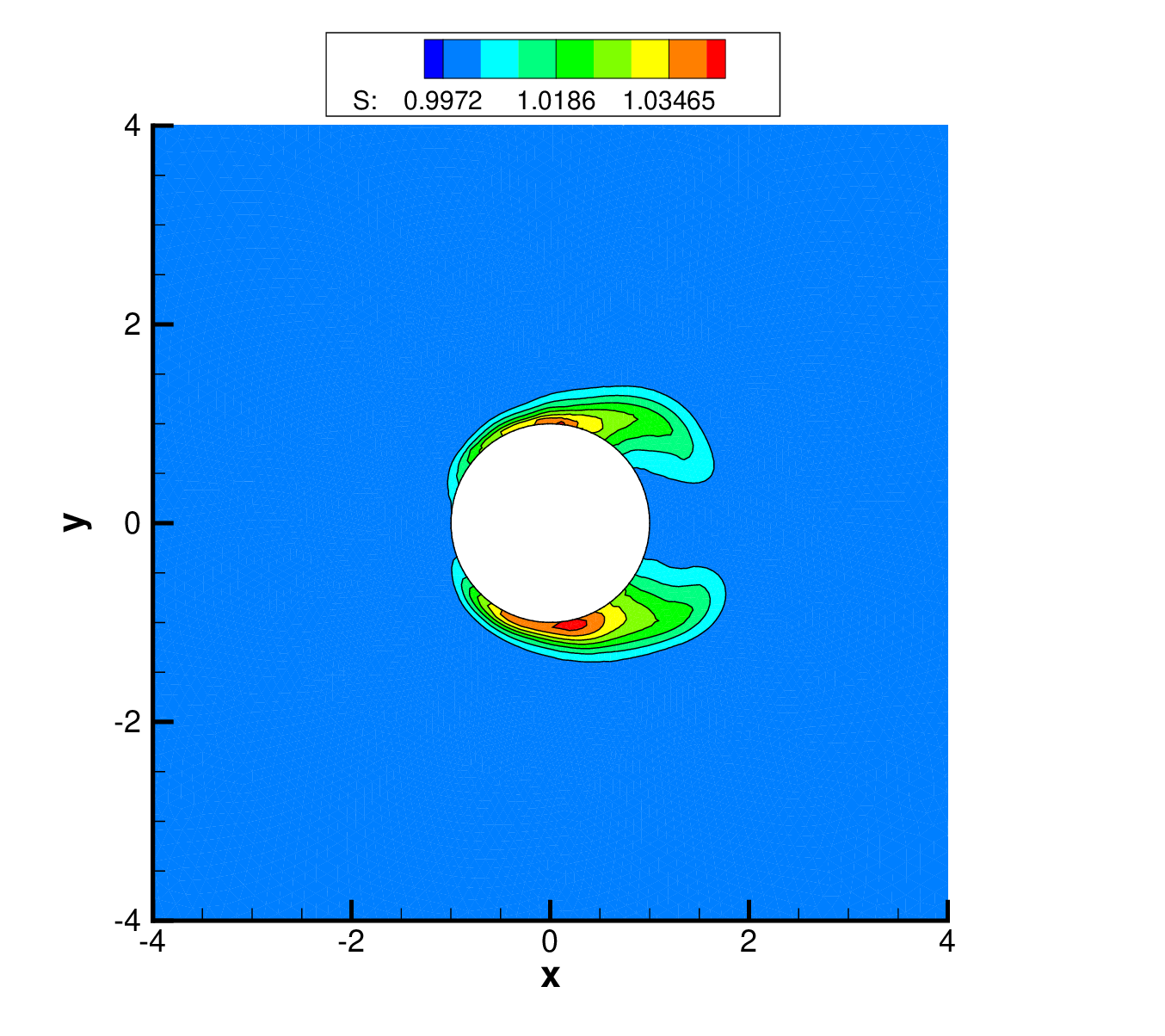}}
  \subfigure[$t=8$ s]{\includegraphics[width=0.4\textwidth]{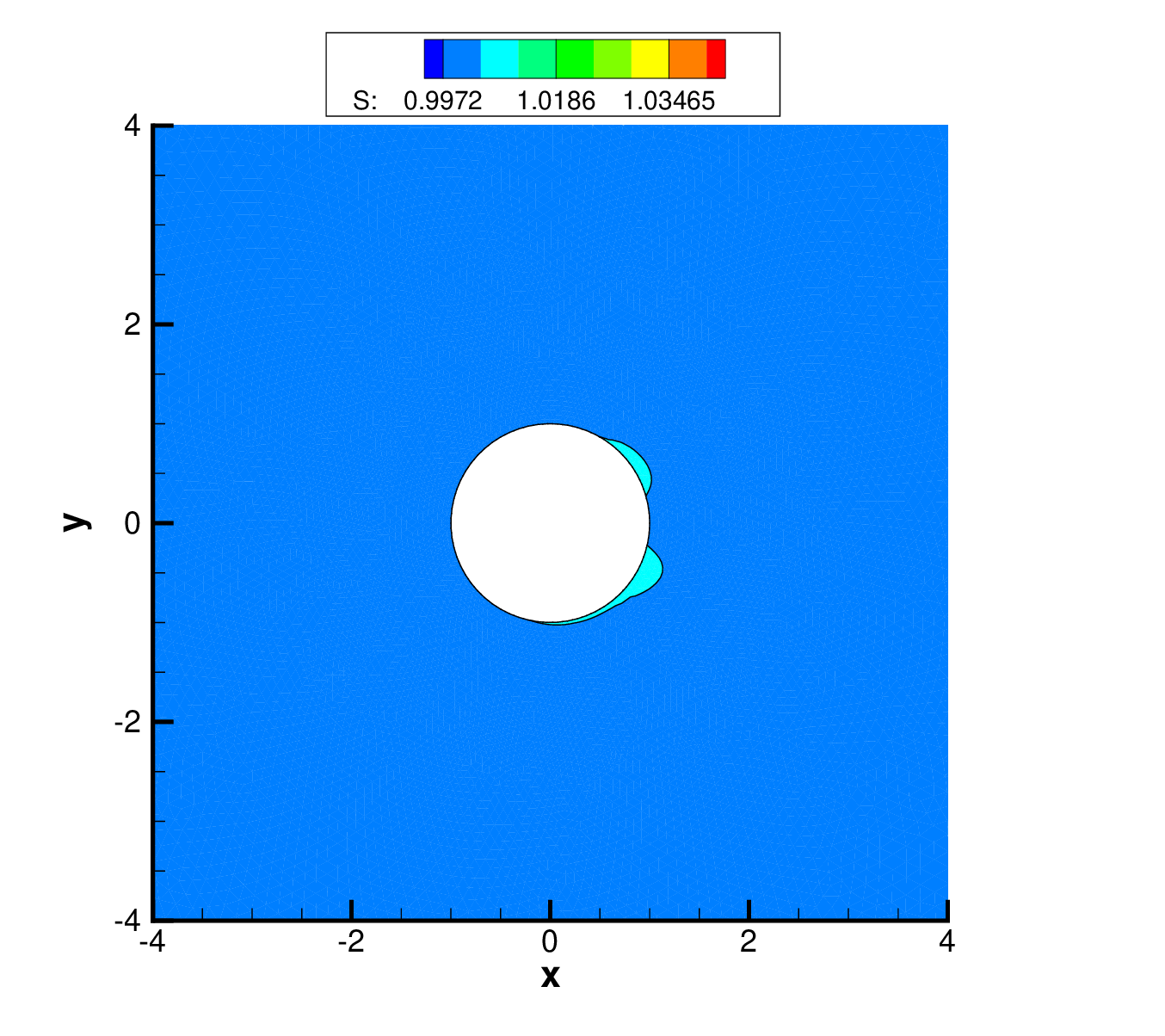}}
  \caption{Vertical oscillating cylinder: entropy $S$ distribution without correction (left) and with correction (right).
  The solution is plot at different times: $t=0.5$ s (top), $t=4$ s (middle) and $t=8$ s (bottom).}
  \label{fig:OC-V-entropy}
\end{figure}

In figure \ref{fig:OC-V-velU} we plot the velocity field $u$ at different times for the vertical oscillating cylinder, and
compare the results obtained with and without the polynomial correction.
The left column shows the results obtained without the correction, while the right column shows the results obtained with the correction.
Once again, the results confirm that the polynomial correction allows us to recover a much smoother velocity field,
while the results obtained without the correction exhibit a thick layer of spurious errors close to the boundary.
In figure \ref{fig:OC-V-entropy}, we also depict the entropy distribution $S=p/\rho^\gamma$ for the two simulations,
where the entropy production of the simulation without the correction is much greater with respect to the solution with the correction.

\section{Conclusions and future perspectives}\label{section:Conclusions}

In this paper, we have presented a new high order method to handle moving curved domains, 
considering both external and internal boundaries, by means of a direct ALE framework for 
the discretization of the hyperbolic Euler system for compressible gas dynamics on triangular and tetrahedral meshes. 
The proposed strategy relies on embedding a shifted boundary polynomial correction \cite{ciallella2023shifted} to 
overcome the second-order geometrical error due to the inconsistent treatment of moving curved 
boundaries discretized with unstructured control volumes, without being obliged to manage 2D and 3D 
moving curvilinear meshes and the challenges that come with that:
mesh generation process, isoparametric transformations and special quadrature formulas.

The idea of extrapolating the boundary conditions through an off-element evaluation of the 
cell polynomial greatly simplifies the high order Taylor expansions, which no longer need the evaluation
of high order derivative terms, and allows us to recover the expected order of accuracy of the scheme
even for complex configurations with moving curved boundaries.
The numerical results presented have shown that this can have a huge impact for the high order numerical simulation
of fluid-structure interaction problems, where boundaries would need the special treatment
provided by the more cumbersome and algorithmically complex 3D moving curvilinear meshes \cite{boscheri2016high}.  

The presented method is very versatile and its perspectives are very promising.
In particular, towards the direction of simulating real compressible fluid-structure interaction problems, 
more complex fluid models will be implemented to simulate the full Navier-Stokes system, 
and coupled with real structure models in Lagrangian formalism to describe with high order of accuracy 
moving deformable bodies. For the moment, the mesh topology is considered to not vary in time.
However, in order to deal with large mesh deformations given by structure large displacements, 
the considered ALE framework will be also extended to deal with topology changes, 
as for example done in \cite{gaburro2020high}.

\section*{Acknowledgments}

WB received financial support by the Italian Ministry of University and Research (MUR) with the PRIN Project 2022 No. 2022N9BM3N and by the "Agence Nationale de la Recherche" (ANR) with project No. ANR-23-EXMA-0004.

\bibliographystyle{elsarticle-num} 
\bibliography{literature}

\end{document}